\input amssym.def
\input amssym.tex

\headline={\ifnum \pageno=1 {\hfill} 
\else{\hss \tenrm -- \folio\ -- \hss}\fi}
\footline={\hfil}

\def\dater{\vglue-10mm\rightline{(\the\day/\the\month/\the\year)}}

\hsize 146mm
\vsize 224mm
\hoffset=6mm
\voffset=8mm
\baselineskip=5mm
\overfullrule =0pt

\font\ninerm=cmr9
 at 10,5pt

\font\GGtitre=cmbx10 at 16pt


\def\og{\leavevmode\raise.30ex
\hbox{$\scriptscriptstyle\langle\!\langle\>$}}    
\def\bigog{\leavevmode\raise.30ex
\hbox{$\langle\!\langle\>$}}  
\def\fg{\leavevmode\raise.30ex
\hbox{$\scriptscriptstyle\>\rangle\!\rangle$}}    
\def\bigfg{\leavevmode\raise.30ex
\hbox{$\>\rangle\!\rangle$}}    

\catcode`\@=11

\font\author=cmcsc10
\font\pauthor=cmcsc10 at 8pt
\font\tenmsx=msam10
\font\sevenmsx=msam10 scaled 700
\font\fivemsx=msam10 scaled 500
\font\tenmsy=msbm10
\font\sevenmsy=msbm10 scaled 700
\font\fivemsy=msbm10 scaled 500
\newfam\msxfam
\newfam\msyfam
\textfont\msxfam=\tenmsx  \scriptfont\msxfam=\sevenmsx
\scriptscriptfont\msxfam=\fivemsx
\textfont\msyfam=\tenmsy  \scriptfont\msyfam=\sevenmsy
\scriptscriptfont\msyfam=\fivemsy

\def\hexnumber@#1{\ifnum#1<10 \number#1\else
\ifnum#1=10 A\else\ifnum#1=11 B\else\ifnum#1=12 C\else
\ifnum#1=13 D\else\ifnum#1=14 E\else\ifnum#1=15 F\fi\fi\fi\fi\fi\fi\fi}

\def\msx@{\hexnumber@\msxfam}
\def\msy@{\hexnumber@\msyfam}
\mathchardef\nmid="3\msy@2D
\mathchardef\varnothing="0\msy@3F
\mathchardef\nexists="0\msy@40
\mathchardef\smallsetminus="2\msy@72
\def\Bbb{\ifmmode\let\next\Bbb@\else
\def\next{\errmessage{Use \string\Bbb\space only in math mode}}\fi\next}
\def\Bbb@#1{{\Bbb@@{#1}}}
\def\Bbb@@#1{\fam\msyfam#1}

\font\tentbl=cmr10 scaled 900
\font\seventbl=cmr7 scaled 900
\font\fivetbl=cmr5 scaled 900

\newfam\tblfam

\textfont\tblfam=\tentbl
\scriptfont\tblfam=\seventbl
\scriptscriptfont\tblfam=\fivetbl


\def \D {{\Bbb D}}

\def \Z {{\Bbb Z}}

\font\ci= eufm10

\def\ciU{\hbox{\ci U}}

\font\cccm=cmmi10 at 7pt

\font\cccm=cmmi10 at 9,5pt

\font\f=cmr10 at 7pt

\def\f1{\hbox{\f 1}}
\def\f2{\hbox{\f 2}}
\def\f3{\hbox{\f 3}}
\def\f4{\hbox{\f 4}}
\def\f5{\hbox{\f 5}}
\def\f6{\hbox{\f 6}}
\def\f7{\hbox{\f 7}}
\def\f8{\hbox{\f 8}}
\def\f9{\hbox{\f 9}}

\def \d {\,{\rm d}}

\def \dm {{\hbox {${1\over 2}$}}}

\def\le{\leqslant}
\def\ge{\geqslant}

\topskip=10pt
\font\sept=cmti9

\def\rightheadline{\ifnum\pageno=\chstart{\hfill}
          \else{\centerline{\sept Chen's double sieve, Goldbach's
conjecture and the twin
prime problem}}
          \hfill \hskip -6mm \tenrm\folio\fi}
\def\leftheadline{\ifnum\pageno=\chstart{\hfill}
          \else\tenrm\folio \hskip -3,5mm \hfill{\centerline{\pauthor 
J. Wu}}\fi}
\headline={\ifnum\pageno=\chstart{\hfill}
\else{\ifodd\pageno\rightheadline\else\leftheadline\fi}\fi}
\footline={\hfill}

\pageno=1
\newcount\chstart
\chstart=\pageno

\noindent{{\it Acta Arithmetica} {\bf 114} (2004), 215--273\hfill}

\vglue 5mm

\centerline{\GGtitre Chen's double sieve,}

\bigskip

\centerline{\GGtitre Goldbach's conjecture and the twin prime problem}

\vskip 8mm

\centerline{\author J. Wu}

\vskip 8mm

{\leftskip=1cm
\rightskip=1cm
{{\bf Abstract}.
For every even integer $N$,
denote by $D(N)$ and $D_{1,2}(N)$ the number of representations of
$N$ as a sum of two primes
and as a sum of a prime and an integer having at most two prime
factors, respectively.
In this paper, we give a new upper bound for $D(N)$ and a new lower
bound for $D_{1,2}(N)$,
which improve the corresponding results of Chen.
We also obtain similar results for the twin prime problem
\par}}

\footnote{{}}{2000 {\it Mathematics Subject Classification}: 
11P32, 11N35, 11N05.}

\vskip 10mm

{\hsize 122mm
\baselineskip=4mm
\overfullrule =0pt
\leftskip=22mm
\font\cccm=cmmi10 at 8pt

\def\cccmh{\hbox{\cccm h}}
\def\cccmH{\hbox{\cccm H}}

\hskip 36mm{\bf Contents}

\vskip 5mm

\ninerm
\S\ 1. Introduction
\dotfill \hskip 6mm 1

\vskip 0,3mm

\S\ 2. Preliminary lemmas
\dotfill \hskip 6mm 5

\vskip 0,3mm

\S\ 3. Chen's double sieve
\dotfill \hskip 6mm 9

\vskip 0,3mm

\S\ 4. Weighted inequalities for sieve function
\dotfill \hskip 4,34mm 15

\vskip 0,3mm

\S\ 5. Functional inequalities between $\cccmH$ and $\cccmh$
\dotfill \hskip 4,34mm 22

\S\ 6. Proofs of Propositions 3 and 4
\dotfill \hskip 4,34mm 28

\vskip 0,3mm

\S\ 7. Proof of Theorem 1
\dotfill \hskip 4,34mm 30

\vskip 0,3mm

\S\ 8. Proof of Theorem 3
\dotfill \hskip 4,34mm 32

\vskip 0,3mm

\S\ 9. Chen's system of weights
\dotfill \hskip 4,34mm 33

\vskip 0,3mm

\S\ 10. Proofs of Theorems 2 and 5
\dotfill \hskip 4,34mm 38

\vskip 0,3mm

\S\ 11. Proof of Theorem 4
\dotfill \hskip 4,34mm 42

\vskip 0,3mm

References
\dotfill \hskip 4,34mm 47
\par}

\vskip 10mm

\noindent{\bf \S\ 1. Introduction}

\medskip

Let $\Omega(n)$ be the number of all prime factors of the integer $n$ with the
convention $\Omega(1)=0$.
For an even integer $N\ge 4$,
we define $D(N)$ as the number of representations of $N$ as a sum of
two primes:
$$D(N) := |\{p\le N :\Omega(N - p) = 1\}|,$$
where and in what follows,
the letter $p$, with or without subscript, denotes a prime number.
The well known Goldbach conjecture can be stated as $D(N)\ge 1$ for
every even integer $N\ge 4$.
A more precise version of this conjecture was proposed by  Hardy \&
Littlewood [15]:
$$D(N) \sim 2 \Theta(N)
\qquad(N\to \infty),
\leqno(1.1)$$
where
$$\Theta(N) := {C_N N\over (\log N)^2}
\qquad{\rm and}\qquad
C_N := \prod_{p\mid N, \, p>2} {p-1\over p-2} \prod_{p>2}
\bigg(1-{1\over (p-1)^2}\bigg).$$
Certainly, the asymptotic formula (1.1) is extremely difficult.
Although the lower bound problem remains open,
the upper bound problem has a rich history.
In 1949 Selberg [25] proved
$$D(N)\le \{16 + o(1)\} \Theta(N)
\leqno(1.2)$$
with the help of his well known $\lambda^2$-upper bound sieve.
By applying Linnik's large sieve method,
C.D. Pan [20] in 1964 improved 16 to 12.
In 1966, Bombieri \& Davenport [1] obtained 8 instead.
Their proof is based on the linear sieve formulas
and the mean value theorem of Bombieri--Vinogradov.
It seems very difficult to prove (1.2) with a constant strictly less
than 8 by the method in [1].
Firstly the linear sieve formulas (see Lemma~2.2 below)
$$X V(z) f\bigg({\log Q\over \log z}\bigg) + {\rm error}
\le S({\cal A}; {\cal P}, z)\le X V(z) F\bigg({\log Q\over \log z}\bigg) + {\rm error}
\leqno(1.3)$$
are the best possible in the sense that taking
$${\cal A}
= {\cal B}_\nu
:= \{n : 1\le n\le x, \Omega(n)\equiv \nu\,({\rm mod} 2)\}
\qquad(\nu=1, 2),$$
the upper and  lower bounds in (1.3) are respectively attained
by $\nu=1$ and $\nu=2$ (see [14], page 239).
Secondly it is hopeless to try to improve the level of distribution
${1\over 2}$ in Bombieri--Vinogradov's theorem.

In 1978, Chen [10] introduced a new idea in Selberg's sieve and proved
$$D(N)\le 7.8342 \, \Theta(N)
\qquad(N\ge N_0).
\leqno(1.4)$$
His sieve machine involves two variables and is quite complicated.
Roughly speaking, for the sequence
$${\cal A} = \{N-p : p\le N\}$$
he introduced two new functions $h(s)$ and $H(s)$ such that (1.3) holds
with $f(s)+h(s)$ and $F(s)-H(s)$ in place of $f(s)$ and $F(s)$, respectively.
The key innovation is to prove $h(s)>0$ and $H(s)>0$ via three weighted
inequalities (see [10], (23), (47), (64), (90), and (91)).
It is worth  pointing out that he did not give complete proofs for
these three inequalities.
Among the three inequalities, the third one is the most complicated
(with 43 terms) and
it seems quite difficult to reconstruct a proof.
Indeed, combining any one of these three inequalities
with the Chen--Iwaniec switching principle (see [7] and [16])
leads to a constant less than 8.
In order to derive a better result,
Chen further introduced a very complicated iterative method.
In 1980, C.B. Pan [19] applied essentially the first weighted
inequality of Chen to get 7.988.
According to [22], Chen's proof is very long and somewhat difficult
to follow, but his idea is clear.

In this paper, inspired by the ideas in [26]
we shall first try to give a more comprehensive treatment
on Chen's double sieve and prove an upper bound sharper than (1.4).

\proclaim Theorem 1.
For sufficiently large $N$, we have
$$D(N) \le 7.8209 \,  \Theta(N).$$

The improvement comes from a new weighted inequality (see Lemma 4.2 below),
which is still quite complicated with 21 terms, but much simpler than Chen's third and
more powerful than his second and third inequality.
Recently Cai \& Lu [6] give another weighted inequality (with 31 terms),
which is simpler but slightly weaker than Chen's third.

\medskip

One way of approaching the lower bound problem in (1.1) is
to give a non trivial lower bound for the quantity
$$D_{1,2}(N) := |\{p\le N : \Omega(N-p)\le 2\}|.$$
In this direction, Chen [7] proved, by his system of
weights and switching principle,
the following famous theorem:
{\sl Every sufficiently large even integer can be written as sum of a
prime and an integer
having at most two prime factors}.
More precisely he established
$$D_{1,2}(N)\ge 0.67 \, \Theta(N)
\qquad
(N\ge N_0).
\leqno(1.5)$$
Then Halberstam \& Richert [14] obtained a better constant 0.689 in
place of 0.67
by a careful numerical calculation.
As they indicated in [14], it would be interesting to know whether a
more elaborate weighting procedure
could be adapted to the purpose of (1.5).
This might lead to numerical improvements and could be important.
In 1978 Chen improved the constant 0.689 of Halberstam \& Richert
to 0.7544 and to 0.81 by two more elaborate systems of weights ([8], [9]).
Very recently by improving Chen's weighting device Cai and Lu [5] obtained $0.8285$, 
which they described as being near to the limit of what could be obtained by the method employed. 

The second aim of this paper is to propose a larger constant.

\proclaim Theorem 2.
For sufficiently large $N$, we have
$$D_{1,2}(N)\ge 0.836 \, \Theta(N).$$

The proof of Theorem 2 is based on a modified version of Chen's
weights (see Lemma 9.2 below),
the linear sieve and the mean value theorems of Pan \& Ding [21] and
of Fouvry [11].

\medskip

A conjecture of the same nature is the twin prime problem,
which can be stated as
$$\pi_2(x) := |\{p\le x : \Omega(p+2)=1\}|\to\infty
\quad(x\to\infty).$$
Similar to (1.1), Hardy \& Littlewood [15] conjectured
$$\pi_2(x)\sim \Pi(x)
\qquad
(x\to \infty),
\leqno(1.6)$$
where
$$\Pi(x) := {Cx\over (\log x)^2}
\qquad{\rm and}\qquad
C := 2 \prod_{p>2} \bigg(1-{1\over (p-1)^2}\bigg).$$
The methods of Selberg, Pan, Bombieri \& Davenport and Chen work in a
similar way
and give upper bounds of this type
$$\pi_2(x)\le \{a + o(1)\} \Pi(x),
\leqno(1.7)$$
where the constant $a$ is half of the corresponding constant in the
Goldbach problem.
Due to the sieve of Rosser--Iwaniec and
mean value theorems of Bombieri, Fouvry, Friedlander and Iwaniec,
the history of (1.7) is much richer than that of (1.2).
We refer the reader to [26] and [6] for a detailed historical description of this problem.
In particular Wu [26] obtained 3.418 in place of $a+o(1)$
by placing these new mean value theorems in Chen's method.
The main difficulty for applying these mean value theorems in [26] is to not 
destroy the fact that the error terms are affected by well factorisable coefficients.
Recently Cai \& Lu [6] improved the constant 3.418 to 3.406.
Our argument in proving Theorem 1 allows us to give a better result.

\proclaim Theorem 3.
For sufficiently large $x$, we have
$$\pi_2(x) \le 3.3996 \, \Pi(x).$$

\smallskip

As an analogue of Theorem 2, Chen [7] proved that
$$\pi_{1,2}(x)\ge 0.335 \, \Pi(x)
\qquad(x\ge x_0),
\leqno(1.8)$$
where
$$\pi_{1,2}(x) := |\{p\le x : \Omega(p + 2)\le 2\}|.$$
The constant 0.335 was improved by many mathematicians.
Like (1.7), the history of (1.8) is much richer than that of (1.5).
A detailed historical description on this problem can be found
in the recent paper of Cai [3].
In particular he obtained 1.0974 in place of 0.335,
which is an improvement of Wu's constant 1.05 [26].
Here we can propose a slightly better result.

\proclaim Theorem 4.
For sufficiently large $x$, we have
$$\pi_{1,2}(x)\ge 1.104 \, \Pi(x).$$

\smallskip

{\bf Remark 1.}
(i)
Theorems 1 and 3 show that the principal terms in the linear sieve
formulas can be improved
in the special cases ${\cal A} = \{N-p : p\le N\}$ or ${\cal A} =
\{p+2 : p\le x\}$ (see the end of Section 3).
This seems to be interesting and important.
Our argument is quite general,
which works for all sequences satisfying the Chen--Iwaniec
switching principle.

(ii)
Certainly we could obtain a better constant than 3.3996 in Theorem 3
if we used mean value theorems of ([11], Corollary 2),
([12], Lemma 6) and ([18], Proposition) as in the proof of Theorem 4.
But the numerical computation involved would be quite complicated.

\medskip

The Chen theorem in short intervals was first studied by Ross [23].
Let $\alpha\in (0, 1)$ be a fixed constant and
define, for $\theta\in(0, 1)$, $x\ge 2$ and even integer $N\ge 4$,
$$\eqalign{
D_{1,2}(N, \theta)
& := |\{\alpha N\le p\le \alpha N+N^\theta : \Omega(N-p)\le 2\}|,
\cr
\pi_{1,2}(x, \theta)
& := |\{x\le p\le x+x^\theta : \Omega(p + 2)\le 2\}|.
\cr}$$
He proved (see [28]) that for $\theta\ge 0.98$, $N\ge N_0(\theta)$
and $x\ge x_0(\theta)$,
$$D_{1,2}(N, \theta)\gg \Xi(N,\theta),
\qquad
\pi_{1,2}(x, \theta)\gg \Pi(x,\theta),$$
where
$$\Xi(N,\theta)
:= {N^\theta\over (\log N)^2}
\prod_{p\mid N, \, p>2} {p-1\over p-2} \prod_{p>2} \bigg(1-{1\over
(p-1)^2}\bigg)$$
and
$$\Pi(x,\theta)
:= {2 x^\theta\over (\log x)^2} \prod_{p>2} \bigg(1-{1\over (p-1)^2}\bigg).$$
The constant 0.98 was further improved to 0.973 by Wu [28],
to 0.9729 by Salerno \& Vitolo [24]
and to 0.972 by Cai \& Lu [4].

Our method allows us to take a smaller exponent.

\proclaim Theorem 5.
For every $\theta\ge 0.971$, $N\ge N_0(\theta)$ and $x\ge x_0(\theta)$, we have
$$D_{1,2}(N, \theta)\ge 0.012 \, \Xi(N,\theta),
\qquad
\pi_{1,2}(x, \theta)\ge 0.006 \, \Pi(x,\theta).$$

\noindent{\bf Acknowledgement.}
The author would like to thank E. Fouvry for his generous help in writing this article, 
and the referee for his very careful reading of the manuscript.

\vskip 5mm

\noindent{\bf \S\ 2. Preliminary lemmas}

\medskip

This section is devoted to present the formula of the Rosser--Iwaniec linear sieve and
some mean value theorem on the distribution of primes in arithmetic progressions, 
which will be needed later.
Before stating these results, it is necessary to recall some definitions.

Let $k$ be a positive integer and
$\tau_k(n)$ the number of ways of writing $n$ as the product of $k$ 
positive integers.
An arithmetical function $\lambda(q)$ is of {\it level} $Q$ and of 
{\it order} $k$ if
$$\lambda(q) = 0
\quad{\rm for}\quad q>Q
\qquad{\rm and}\qquad
|\lambda(q)|\le \tau_k(q)
\quad{\rm for}\quad q\ge 1.$$
We say that $\lambda$ is {\it well factorable}
if for every decomposition $Q = Q_1 Q_2$ $(Q_1, Q_2\ge 1)$ there 
exist two arithmetical functions
$\lambda_1$ and $\lambda_2$ of level $Q_1, Q_2$ and of order $k$
such that $\lambda = \lambda_1*\lambda_2$.

\proclaim Lemma 2.1.
If $\lambda'$ is an arithmetical function of level $Q'$ $(\le Q)$ and 
of order $k'$,
then $\lambda*\lambda'$ is well factorable of level $Q Q'$ and of order $k+k'$.

Let ${\cal A}$ be a finite sequence of integers and
${\cal P}$ a set of prime numbers.
For $z\ge 2$, we put $P(z) := \prod_{p<z, \, p\in {\cal P}} p$ and 
define the sieve function
$$S({\cal A}; {\cal P}, z) := |\{a\in {\cal A} : (a, P(z))=1\}|.$$
If $d$ is a square-free integer with all its prime factors belonging 
to ${\cal P}$,
we denote by ${\cal A}_d$ the set of elements of ${\cal A}$ divisible 
by $d$ and
we write the following approximate formula
$$|{\cal A}_d|
={w(d)\over d}X+r({\cal A},d),
\leqno(2.1)$$
where $X>1$ is independent of $d$, and $w(d)$ is a multiplicative 
function satisfying
$$0\le w(p)< p
\qquad{\rm for}\qquad
p\in {\cal P}.
\leqno(2.2)$$
We also define
$$V(z) := \prod_{p<z} \bigg(1-{w(p)\over p}\bigg)$$
and suppose that there exists an absolute constant $K>1$ such that
$${V(z_1)\over V(z_2)}\le {\log z_2\over \log z_1}\bigg(1+{K\over 
\log z_1}\bigg)
\qquad
(z_2\ge z_1\ge 2).
\leqno(2.3)$$

The formula of the Rosser--Iwaniec linear sieve [17] is stated as follows.

\proclaim Lemma 2.2.
Let $0<\varepsilon<{1\over 8}$ and $2\le z\le Q^{1/2}$.
Under the assumptions $(2.1)$, $(2.2)$ and $(2.3)$, we have
$$S({\cal A}; {\cal P}, z)
\le X V(z)
\bigg\{F\bigg({\log Q\over \log z}\bigg) + E\bigg\}
+ \sum_{l<L} \sum_{q\mid P(z)} \lambda_l^+(q) r({\cal A},q)
\leqno(2.4)$$
and
$$S({\cal A}; {\cal P}, z)
\ge X V(z)
\bigg\{f\bigg({\log Q\over \log z}\bigg) + E\bigg\}
- \sum_{l<L} \sum_{q\mid P(z)} \lambda_l^-(q) r({\cal A},q).
\leqno(2.5)$$
In these formulas, $L$ depends only on $\varepsilon$,
the $\lambda_l^\pm$ are well factorable functions of order 1 and of level $Q$,
and $E\ll \varepsilon + \varepsilon^{-8} e^K/(\log Q)^{1/3}$.
The functions $F, f$ are defined by
$$\eqalign{F(u) = 2e^\gamma/u,&
\quad\qquad
f(u)=0
\qquad
(0<u\le 2),
\cr
(uF(u))' = f(u-1),&
\qquad\quad
(uf(u))' = F(u-1)
\quad
(u>2),
\cr}
\leqno(2.6)$$
where $\gamma$ is Euler's constant.

As usual, we denote by $\mu(q)$ M\"obius' function, $\varphi(q)$ 
Euler's function
and $\nu(q)$ the number of distinct prime factors of $q$.
Define
$$\pi(y; q, a, m)
:= \sum_{\scriptstyle mp\le y\atop\scriptstyle mp\equiv a ({\rm mod} \, q)} 1,
\qquad
{\rm li}(y) := \int_2^y {\d t\over \log t}$$
and
$$\overline E_0(y; q, a, m)
:= \pi(y; q, a, m)-{{\rm li}(y/m)\over \varphi(q)}.$$

The next lemma is due to Pan \& Ding [21], which implies 
Bombieri--Vinogradov's theorem.
Here we state it in the form of ([22], Corollary 8.12).

\proclaim Lemma 2.3.
Let $f(m)\ll 1$ and $\alpha\in (0,1]$.
Let $r_1(y)$ be a positive function depending on $x$ and satisfying
$$r_1(y)\ll x^\alpha,
\qquad
y\le x.$$
Let $r_2(m)$ be a positive function depending on $x$ and $y$, and satisfying
$$mr_2(m)\ll x,
\qquad
m\le x^\alpha,
\qquad
y\le x.$$
Then for every $A>0$, there exists a constant $B=B(A)>0$ such that
$$\eqalign{
& \sum_{q\le \sqrt x/(\log x)^B}\mu(q)^2 3^{\nu(q)}
\max_{y\le x} \max_{(a,q)=1}
\Big|\sum_{\scriptstyle m\le x^{1-\alpha}\atop\scriptstyle (m, q)=1}
f(m) \overline E_0(y; q, a, m)\Big|
\ll {x\over (\log x)^A},
\cr
& \sum_{q\le \sqrt x/(\log x)^B}\mu(q)^2 3^{\nu(q)}
\max_{y\le x} \max_{(a,q)=1}
\Big|\sum_{\scriptstyle m\le x^{1-\alpha}\atop\scriptstyle (m, q)=1}
f(m) \overline E_0(mr_1(y); q, a, m)\Big|
\ll {x\over (\log x)^A},
\cr
& \sum_{q\le \sqrt x/(\log x)^B}\mu(q)^2 3^{\nu(q)}
\max_{y\le x} \max_{(a,q)=1}
\Big|\sum_{\scriptstyle m\le x^{1-\alpha}\atop\scriptstyle (m, q)=1}
f(m) \overline E_0(mr_2(m); q, a, m)\Big|
\ll {x\over (\log x)^A}.
\cr}$$

\smallskip

In order to prove Theorem 5,
it is necessary to generalize the mean value theorem of Pan \& Ding
in short intervals.
Such a result was established by Wu ([27], theorem 2).

\proclaim Lemma 2.4.
Let $f(m)\ll 1$, $\varepsilon$ be an arbitrarily small positive 
number and define
$$H(y, h, q, a, m)
:= \pi(y+h; q, a, m) - \pi(y; q, a, m)
- {{\rm li}((y+h)/m)-{\rm li}(y/m)\over\varphi(q)}.$$
Then for any $A>0$, there exists a constant $B=B(A)>0$
such that
$$\sum_{q\le Q} \mu(q)^2 3^{\nu(q)}
\max_{(a,q)=1} \max_{h\le x^\theta} \max_{x/2<y\le x}
\big|\sum_{m\le M, \, (m,q)=1} f(m) H(y, h, q, a, m)\big|
\ll {x^\theta\over (\log x)^A}$$
for $x\ge 10, \, {3\over 5} + \varepsilon\le\theta\le 1$,
$Q = x^{\theta-1/2} /(\log x)^B$ and $M = x^{(5\theta-3)/2-\varepsilon}$.

In the proofs of Theorem 3 and 4,
we shall need some mean value theorems with well factorable or almost 
well factorable coefficients.

Let $M\ge 1$, $N\ge 1$ and $X:=MN$.
Let $\{\alpha_m\}$ and $\{\beta_n\}$ be two sequences of order $k$
supported in $[M, 2M]$ and $[N, 2N]$ respectively.
We also suppose the conditions below:

(i)
For all $B$, the equality
$$\sum_{n\equiv n_0 ({\rm mod} \, k), \, (n,d)=1} \beta_n
= {1\over \varphi(k)} \sum_{(n,dk)=1} \beta_n
+ O_{B,k}\big(N \tau_k(d)/(\log 2N)^B\big)$$
holds for $d\ge 1$, $k\ge 1$ and $(k,n_0)=1$.

(ii)
If $n$ has a prime factor $p$ with $p<\exp\{(\log\log n)^2\}$, then 
$\beta_n=0$.

\smallskip

The following result is an immediate consequence of Corollary 2 of [11],
Lemma 6 of [12] and the proposition of [18].

\proclaim Lemma 2.5.
Under the conditions (i) and (ii) above, for any $A$ and for any 
$\varepsilon>0$ we have
$$\sum_{(q,a)=1} \lambda(q)
\Big(\sum_{mn\equiv a ({\rm mod} \, q)} \alpha_m \beta_n
- {1\over \varphi(q)} \sum_{(mn,q)=1} \alpha_m \beta_n\Big)
\ll_{\varepsilon, A} {X\over (\log X)^A}$$
uniformly for $|a|\le (\log X)^A$ and $\nu := \log N/\log X$ 
$(\varepsilon\le \nu\le 1-\varepsilon)$.
Here $\lambda(q)$ is a well factorisable function of order 1 and of 
level $Q:=X^{\theta(\nu)-\varepsilon}$,
where $\theta(\nu)$ is given by
$$\theta(\nu) = \cases{
{6-5\nu\over 10}               & for $\,\,\varepsilon\le \nu\le {1\over 15}$,
\cr\noalign{\medskip}
{1\over 2}+\nu                 & for ${1\over 15}\le \nu\le {1\over 10}$,
\cr\noalign{\medskip}
{5-2\nu\over 8}                & for ${1\over 10}\le \nu\le {3\over 14}$,
\cr\noalign{\medskip}
{3+2\nu\over 6}                & for ${3\over 14}\le \nu\le {1\over 4}$,
\cr\noalign{\medskip}
{2-\nu\over 3}                 & for ${1\over 4}\le \nu\le {2\over 7}$,
\cr\noalign{\medskip}
{2+\nu\over 4}                 & for ${2\over 7}\le \nu\le {2\over 5}$,
\cr\noalign{\medskip}
1 - \nu                        & for ${2\over 5}\le \nu\le {1\over 2}$,
\cr\noalign{\medskip}
{1\over 2}                     & for ${1\over 2}\le \nu\le 1-\varepsilon$.
\cr}$$

\noindent{\sl Proof}.
The value $(6-5\nu)/10$ in $[\varepsilon, 1/15]$ comes from 
Corollary 2 (ii) of [11].
The intervals $[1/15, 1/10]$ and $[1/10, 3/14]$ follow from 
the proposition of [18] by decomposing $\lambda = \lambda_1*\lambda_2$ with
$$Q_1=Q=x^{1/2-\varepsilon},
\qquad
Q_2=R=Nx^{-\varepsilon}$$
and
$$Q_1=Q=x^{5/8-\varepsilon} N^{-5/4},
\qquad
Q_2=R=Nx^{-\varepsilon},$$
respectively.
The remaining case is Lemma 6 of [12].
\hfill
$\square$

\goodbreak
\smallskip

The next lemma is Corollary 2 (i) of [11].
This is the first result,
which is valid uniformly for $|a|\le X$ and has the level of 
distribution $>\dm$.

\proclaim Lemma 2.6.
Under the conditions (i) and (ii) above, for any $A$ and for any 
$\varepsilon>0$ we have
$$\sum_{(q,a)=1} \lambda(q)
\Big(\sum_{mn\equiv a ({\rm mod} \, q)} \alpha_m \beta_n
- {1\over \varphi(q)} \sum_{(mn,q)=1} \alpha_m \beta_n\Big)
\ll_{\varepsilon, A} {X\over (\log X)^A}$$
uniformly for $|a|\le X$ and $\varepsilon\le \nu := \log N/\log X\le 
{1\over 10}$.
Here $\lambda(q)$ is a well factorisable function of order 1
and of level $Q:=X^{5(1-\nu)9-\varepsilon}$.

\smallskip

As usual define
$$\pi(y; q,a)
:= \sum_{p\le y, \, p\equiv a ({\rm mod} \, q)} 1.$$
The following result is due to 
Bombieri, Friedlander \& Iwaniec ([2], theorem 10).

\proclaim Lemma 2.7.
Let $\lambda$ be a well factorable function of order $k$ and of level 
$Q=x^{4/7-\varepsilon}$.
For any $\varepsilon>0$ and any $A$,
we have uniformly for $x\ge 3$ and $|a|\le (\log x)^A$,
$$\sum_{(q,a)=1} \lambda(q) 
\bigg(\pi(x; q,a)-{{\rm li}(x)\over \varphi(q)}\bigg)
\ll_{\varepsilon,k,A} {x\over (\log x)^A}.$$

When we use the weighted inequality,
some coefficients are merely ``almost well factorable''.
So we need the following results, 
due to Fouvry \& Grupp ([13], theorem 2 and the corollary).

\proclaim Lemma 2.8.
Let $\lambda$ be a well factorable function of level $Q_1$ and of order $k$,
$\xi$ an arithmetical function satisfying the conditions
$|\xi(q_2)|\le \log x$
and $\xi(q_2)=0$ $(q_2>Q_2)$ and let $\Lambda$ be the von Mangoldt function.
Then we have for any integer $a$, any $\varepsilon>0$ and any $A>0$,
$$\sum_{(q_1q_2,a)=1} \lambda(q_1) \xi(q_2)
\bigg(\pi(x; q_1q_2,a)-{{\rm li}(x)\over \varphi(q_1q_2)}\bigg)
\ll_{a,\varepsilon,k,A} {x\over (\log x)^A},$$
so long as one of the following three conditions is true:
$$\leqalignno{
& Q_2\le Q_1,
\qquad
Q_1 Q_2\le x^{4/7-\varepsilon},
& ({\rm C}.1)
\cr\noalign{\smallskip}
& Q_2\ge Q_1,
\qquad
Q_1 Q_2^6\le x^{2-\varepsilon},
& ({\rm C}.2)
\cr\noalign{\smallskip}
& \xi(q)=\Lambda(q),
\qquad
Q_1 Q_2\le x^{11/20-\varepsilon},
\qquad
Q_2\le x^{1/3-\varepsilon}.
& ({\rm C}.3)
\cr}$$

The next two lemmas also are useful when we apply the switching principle.

\proclaim Lemma 2.9 ([26], Lemma 7). 
Let $\lambda$ be a well factorable function of level $Q := x^{4/7-\varepsilon}$ and of order $k$.
Let $\eta>0$ and $\{\varepsilon_i\}_{1\le i\le r}$ be real numbers such that 
$$\varepsilon_i\ge \eta,
\qquad
\varepsilon_1+\varepsilon_2+\cdots+\varepsilon_r=1.$$
Then for any integer $a$, any $\varepsilon>0$ and any $A>0$, we have
$$\sum_{(q,a)=1} \lambda(q)
\Big(
\sum_{\scriptstyle p_1\cdots p_r\equiv a ({\rm mod}\, q)
\atop{\atop\scriptstyle x^{\varepsilon_i}<p_i\le 2x^{\varepsilon_i}\,(1\le i\le r)}} 1
-{1\over \varphi(q)}
\sum_{\scriptstyle (p_1\cdots p_r,q)=1
\atop{\atop\scriptstyle (x^{\varepsilon_i}<p_i\le 2x^{\varepsilon_i}\,(1\le i\le r)}} 1
\Big) 
\ll_{a,\varepsilon,k,A} {x\over (\log x)^A}.$$ 

\proclaim Lemma 2.10 ([26], Lemma 12).
Let $x\ge 2$ and $y=x^{1/u}$. Then
$$\sum_{\scriptstyle n\le x\atop\scriptstyle p\mid n\Rightarrow p\ge y} 1
= {x\over \log y} \omega(u) + O\bigg({x\over (\log y)^2}\bigg),$$
where $\omega(u)$ is Buchstab's function defined by
$$\omega(u) = 1/u \quad (1\le u\le 2)
\quad{\rm and}\quad
\big(u\omega(u)\big)' = \omega(u-1) \quad (u\ge 2).$$
Moreover we have
$\omega(u)\le 0.561522$
$(u\ge 3.5)$
and
$\omega(u)\le 0.567144$
$(u\ge 2)$.

\vskip 5mm

\noindent{\bf \S\ 3. Chen's double sieve}

\medskip

We shall sieve the sequence
$${\cal A} := \{N - p : p\le N\}.$$
Let $\delta>0$ be a sufficiently small number and $k\in \Z$.
Put 
$$Q := N^{1/2-\delta},
\qquad
\underline d := Q/d,
\qquad
{\cal L} := \log N,
\qquad
W_k := N^{\delta^{1+k}}.$$
Let $\Delta$ be a real number with
$1 + {\cal L}^{-4}\le \Delta<1 + 2 {\cal L}^{-4}$.
We put ${\cal P}(N) := \{p : (p, N) = 1\}$
and denote by $\pi_{[Y, Z)}$ the characteristic function of the set
${\cal P}(N) \cap [Y, Z)$.
For $k\in \Z^+$ and $N\ge 2$,
let $\ciU_k(N)$ be the set of all arithmetical functions $\sigma$ which 
can be written as the form
$$\sigma = \pi_{[V_1/\Delta, V_1)}*\cdots*\pi_{[V_i/\Delta, V_i)},$$
where $i$ is an integer with $0\le i\le k$, and $V_1, \dots, V_i$ are 
real numbers satisfying
$$\cases{
V_1^2\le Q,                            & {}
\cr\noalign{\smallskip}
V_1 V_2^2\le Q,                        & {}
\cr\noalign{\smallskip}
\cdots\cdots\cdots\cdots\cdot          & {}
\cr\noalign{\smallskip}
V_1\cdots V_{i-1} V_i^2\le Q,          & {}
\cr\noalign{\smallskip}
V_1\ge V_2\ge \cdots\ge V_i\ge W_k.    & {}
\cr}
\leqno(3.1)$$
By convention, $\sigma$ is the characteristic function of the set 
$\{1\}$ if $i=0$.
 From this definition and Lemma 2.1, we see immediately the following result.

\proclaim Lemma 3.1.
{\rm (i)}
We have $\ciU_k(N)\subset \ciU_{k+1}(N)$ for $k\in \Z^+$.

\vskip -2mm

{\sl {\rm (ii)}
Let $\sigma = \pi_{[V_1/\Delta, V_1)}*\cdots*\pi_{[V_i/\Delta, 
V_i)}\in \ciU_k(N)$.
Then $\sigma$ is well factorable of level $V := V_1 \cdots V_i$ and 
of order $i$.
If $\lambda$ is well factorable of level $Q/V$ and of order 1,
then $\sigma * \lambda$ is well factorable of level $Q$ and of order $k+1$.}

\medskip

Let $F$ and $f$ be defined as in (2.6) and let
$$A(s) := sF(s)/2e^\gamma
\qquad{\rm and}\qquad
a(s) := sf(s)/2e^\gamma,
\leqno(3.2)$$
We introduce the notation
$$\Phi(N, \sigma, s) := \sum_d \sigma(d) S({\cal A}_d; {\cal P}(d N), 
\underline d^{1/s}),
\qquad
\Theta(N, \sigma)
:= 4 {\rm li}(N)
\sum_d {\sigma(d) C_{d N}\over \varphi(d)\log \underline d}.$$

For $k\in \Z^+$, $N_0\ge 2$ and $s\in [1, 10]$, we define
$H_{k,N_0}(s)$ and $h_{k,N_0}(s)$ as the supremum of $h\ge -\infty$
such that for all $N\ge N_0$ and $\sigma\in \ciU_k(N)$ one has
the following inequalities
$$\Phi(N, \sigma, s)
\le \{A(s) - h\} \, \Theta(N, \sigma)$$
and
$$\Phi(N, \sigma, s)
\ge \{a(s) + h\} \, \Theta(N, \sigma)$$
respectively.

From this definition, we deduce immediately the following result.

\proclaim Lemma 3.2.
For $k\in \Z^+, N\ge N_0, s\in [1, 10]$ and $\sigma\in \ciU_k(N)$,
we have
$$\Phi(N, \sigma, s)
\le \{A(s) - H_{k,N_0}(s)\} \Theta(N, \sigma)
\leqno(3.3)$$
and
$$\Phi(N, \sigma, s)
\ge \{a(s) + h_{k,N_0}(s)\} \Theta(N, \sigma).
\leqno(3.4)$$

Obviously $H_{k,N_0}(s), h_{k,N_0}(s)$ are decreasing on $N_0$,
and they are also decreasing on $k$ by Lemma 3.1.
Hence we can write
$$\eqalign{
H_k(s)
& := \lim_{N_0\rightarrow \infty} H_{k,N_0}(s),
\cr
H(s)
& := \lim_{k\rightarrow \infty} H_k(s),
\cr}
\qquad
\eqalign{
h_k(s)
& := \lim_{N_0\rightarrow \infty} h_{k,N_0}(s),
\cr
h(s)
& := \lim_{k\rightarrow \infty} h_k(s).
\cr}$$

\proclaim Lemma 3.3.
For $N\ge 2$ and $\sigma= \pi_{[V_1/\Delta, 
V_1)}*\cdots*\pi_{[V_i/\Delta, V_i)}\in \ciU_k(N)$,
we have
$$\leqalignno{
& {\cal L}^{-5k}\ll_{\delta, k} \sum_d \sigma(d)/d\ll_{\delta, k} 1,
& (3.5)
\cr
& \sum_d \sigma(d)\ll_{\delta, k} V_1 \cdots V_i,
& (3.6)
\cr
& \Theta(N, \sigma)\gg_{\delta, k} N / {\cal L}^{5k+2}.
& (3.7)
\cr}$$

\noindent{\sl Proof}.
Let $\sigma = \pi_{[V_1/\Delta, V_1)}*\cdots*\pi_{[V_i/\Delta, 
V_i)}\in \ciU_k(N)$.
We have
$$\sum_d {\sigma(d)\over d}
= \prod_{1\le j\le i}
\sum_{p_j\in {\cal P}(N)\cap [V_j/\Delta, V_j)} {1\over p_j}.
\leqno(3.8)$$
The prime number theorem of the form
$\sum_{p\le x} 1 = {\rm li}(x) + O(x \, e^{-2(\log x)^{1/2}})$
implies
$$\eqalign{
\sum_{p_j\in {\cal P}(N)\cap [V_j/\Delta, V_j)} {1\over p_j}
& = \sum_{V_j/\Delta\le p_j<V_j} {1\over p_j}
- \sum_{V_j/\Delta\le p_j<V_j, \, p_j\mid N} {1\over p_j}
\cr
& = \log\bigg({\log V_j\over \log(V_j/\Delta)}\bigg)
+ O\bigg(e^{-\log^{1/2}(V_j/\Delta)}
+ {{\cal L}\over V_j \log {\cal L}}\bigg).
\cr}$$
Therefore our choice of $\Delta$ and (3.1) give us
$$\sum_{p_j\in {\cal P}(N)\cap [V_j/\Delta, V_j)} 1/ p_j
\asymp_{\delta, k} {\cal L}^{-5}.
\leqno(3.9)$$
Now (3.5) follows from (3.8) and (3.9).

Since $\sigma(d)\not = 0$ implies $d\le V_1 \cdots V_i$,
the second inequality in (3.5) implies (3.6).
Noticing $\Theta(N, \sigma)\gg N{\cal L}^{-2} \sum_d \sigma(d)/d$,
we obtain (3.7) by the first inequality in (3.5).
$\hfill\square$

\proclaim Proposition 1.
For $k\in \Z^+$ and $s\in [1,10]$, we have 
$$H_k(s)\ge 0
\qquad\hbox{and}\qquad
h_k(s)\ge 0.$$

\noindent{\sl Proof}.
We shall prove only the first inequality.
The second one can be treated similarly.
Let $\sigma = \pi_{[V_1/\Delta, V_1)}*\cdots*\pi_{[V_i/\Delta, 
V_i)}\in \ciU_k(N)$.
We use Lemma 2.2 with
$$X = {{\rm li}(N)\over \varphi(d)},
\qquad
w(p) = \cases{
p/(p-1) & if $p\in {\cal P}(N)$,
\cr\noalign{\smallskip}
0       & otherwise
\cr}$$
to estimate $\sigma(d) S({\cal A}_d; {\cal P}(d N), \underline d^{1/s})$.
By Merten's formula and (3.1), we can infer that for any $\varepsilon>0$
$$V(\underline d^{1/s})
= \{1 + O_{\delta, k}(\varepsilon)\}
{2sC_{d N}\over e^\gamma \log \underline d}.
\leqno(3.10)$$
By using Lemma 2.2 and (3.10), we deduce
$$\leqalignno{\qquad
\sigma(d) S({\cal A}_d; {\cal P}(d N), \underline d^{1/s})
& \le 4 {\rm li}(N) {\sigma(d) C_{d N}\over \varphi(d)\log \underline d}
\bigg\{A\bigg({\log(Q/V)\over \log \underline d^{1/s}}\bigg) + 
O_{\delta, k}(\varepsilon)\bigg\}
& (3.11)
\cr
& \quad
+ \sum_{l<L} \sigma(d)
\sum_{q\mid P(\underline d^{1/s})} \lambda_l^+(q) r({\cal A}_d, q),
\cr}$$
where $\lambda_l^+(q)$ is well factorable of level $Q/V$ with
$V := V_1\cdots V_i$ and of order 1.

If $\sigma(d)\not =0$, we have $d\in [V/\Delta^i, V]$, which implies
$0\le \log V-\log d\le i\log \Delta\le 2k{\cal L}^{-4}$.
 From this we deduce that
$A\big(\log(Q/V)/\log \underline d^{1/s}\big) = A(s) + O_{\delta, 
k}(\varepsilon)$.
Inserting (3.11) and summing over $d$, we obtain
$$\Phi(N, \sigma, s)\le \{A(s) + O_{\delta, k}(\varepsilon)\} 
\Theta(N, \sigma) + R,
\leqno(3.12)$$
where
$$R := \sum_{l<L} \sum_d\sigma(d)
\sum_{q\mid P(\underline d^{1/s})} \lambda_l^+(q) r({\cal A}_d, q).$$
Let $q\mid P(\underline d^{1/s})$. It is clear that $\mu(q)^2=1$ and 
$(Nd,q)=1$.
Thus we have
$$\eqalign{r({\cal A}_d, q)
& = |{\cal A}_{dq}| - {\rm li}(N)/\varphi(dq)
\cr
& = \pi(N; dq, N) - {\rm li}(N)/\varphi(dq).
\cr}$$
Hence we can see, by using Lemmas 3.1(ii) and 2.3, that
$$\leqalignno{R
& \ll_\varepsilon
\sum_{q\le Q} \tau_{k+1}(q) |\pi(N; dq, N) - {\rm li}(N)/\varphi(dq)|
& (3.13)
\cr
& \ll_{\delta, k, \varepsilon} N/{\cal L}^{5k+3}.
\cr}$$
 From (3.7), (3.12) and (3.13), we deduce
$$\Phi(N, \sigma, s)
\le \{A(s) + O_{\delta, k}(\varepsilon)\} \Theta(N, \sigma),$$
which implies, by the definition of $H_{k,N_0}(s)$,
for any $\varepsilon>0$ and sufficiently large $N_0$
$$H_{k,N_0}(s)\ge - O_{\delta, k}(\varepsilon).$$
First making $N_0\to \infty$ and then $\varepsilon\to 0$, we 
obtain $H_k(s)\ge 0$.
$\hfill\square$

\proclaim Proposition 2.
For $2\le s\le s'\le 10$, we have
$$h(s)\ge h(s') + \int_{s-1}^{s'-1} {H(t)\over t} \d t$$
and
$$H(s)\ge H(s') + \int_{s-1}^{s'-1} {h(t)\over t} \d t.$$

\noindent{\sl Proof}.
We shall only prove the first inequality as the second one can be 
established in the same way.

Let $k\ge 0$ and $\sigma=\pi_{[V_1/\Delta, V_1)}*\cdots 
*\pi_{[V_i/\Delta, V_i)}\in \ciU_k(N)$.
By Buchstab's identity, we write
$$\Phi(N, \sigma, s)
= \Phi(N, \sigma, s')
- \sum_d \sigma(d) \sum_{\underline d^{1/s'}\le p<\underline d^{1/s}}
S({\cal A}_{dp}; {\cal P}(d N), p).
\leqno(3.14)$$

Next we shall give an upper bound for the last double sums $S$.
The idea is to prove that the characteristic function of $dp$ belongs 
to $\ciU_{k+1}(N)$.
Thus $S$ can be estimated by a function $H_{k+1, N_0}$.
We put $V := V_1\cdots V_i$, $\underline V := Q/V$ and $\alpha_j := 
\underline V^{1/s'}\Delta^j$.
Let $r$ be the integer satisfying $\alpha_r\le \underline V^{1/s}< 
\alpha_{r+1}$.
Noticing that
$\sigma(d)\not=0
\Rightarrow
\underline V^{1/s'}\le \underline d^{1/s'}$ and
$\underline V^{1/s}\le \underline d^{1/s},$
we deduce
$$\leqalignno{S
& \le \sum_d \sigma(d) \sum_{\alpha_0\le p<\alpha_r}
S({\cal A}_{dp}; {\cal P}(d N), p)
+ R_1
& (3.15)
\cr
& = \sum_{1\le j\le r} \sum_{d,\,p}
\sigma(d)
\pi_{[\alpha_{j-1}, \alpha_j)}(p)
S({\cal A}_{dp}; {\cal P}(d p N), (\underline{dp})^{1/s^*})
+ R_1,
\cr}$$
where $s^* := \log \underline d/\log p-1$ and
$$R_1 :=
\sum_d \sigma(d) \!\! \sum_{\alpha_r\le p<\underline d^{1/s}} \!\!
S({\cal A}_{dp}; {\cal P}(d N), p).$$

We would prove that $\sigma*\pi_{[\alpha_{j-1}, \alpha_j)}\in \ciU_{k+1}(N)$.
It suffices to verify that $V_1, V_2, \dots, V_i, \alpha_j$ satisfy 
(3.1) for $j\le r$.
If $V_i\ge \alpha_j$, then
$V_1 V_2 \cdots V_i \alpha_j^2
\le V \underline V^{2/s}
= Q^{2/s}V^{1-2/s}
\le Q$
and
$\alpha_j
\ge \underline V^{1/s'}
\ge V_i^{1/s'}
\ge W_k^{1/s'}\ge W_{k+1}.$
If $V_1\ge \cdots \ge V_l\ge \alpha_j\ge V_{l+1}\ge \cdots \ge V_i$,
we have
$V_1 \cdots V_l \alpha_j V_{l+1} \cdots V_n^2
\le V \alpha_j^2
\le V \underline V^{2/s}
\le Q$
for $l<n\le i$.
Thus $\sigma*\pi_{[\alpha_{j-1}, \alpha_j)}\in \ciU_{k+1}(N)$.

Since $s^*$ depends on $d$ and $p$,
we replace it by a suitable quantity independent of $d$ and $p$
such that we can use (3.3) with $H_{k+1, N_0}$.
For this we introduce
$s_1 := \log (\underline V/\alpha_j)/\log \alpha_j$,
$s_2 := \log (\underline V/\alpha_{j-i-1})/\log \alpha_{j-1}$.
Noticing that
$\sigma(d)
\pi_{[\alpha_{j-1}, \alpha_j)}(p)
\not=0
\Rightarrow
s_1\le s^*\le s_2$,
we deduce from (3.15) that
$$S\le
\sum_{1\le j\le r} \sum_{d,\,p}
\sigma(d)
\pi_{[\alpha_{j-1}, \alpha_j)}(p)
S({\cal A}_{dp}; {\cal P}(d p N), (\underline{dp})^{1/s_1})
+ R_1 + R_2$$
where
$$R_2 :=
\sum_{1\le j\le r} \sum_{d, \, p}
\sigma(d)
\pi_{[\alpha_{j-1}, \alpha_j)}(p)
\big\{
S({\cal A}_{dp}; {\cal P}(d p N), (\underline{dp})^{1/s^*})
- S({\cal A}_{dp}; {\cal P}(d p N), (\underline{dp})^{1/s_1})
\big\}.$$
Now we can use (3.3) in Lemma 3.2 to write
$$\eqalign{S
& \le \sum_{1\le j\le r}
\{A(s_1) - H_{k+1, N_0}(s_1)\}
\Theta(N, \sigma*\pi_{[\alpha_{j-1}, \alpha_j)}) + R_1 + R_2
\cr
& \le 4 {\rm li}(N)
\sum_d {\sigma(d) C_{d N}\over \varphi(d) \log \underline d}
\sum_{\alpha_0\le p<\alpha_r}
{A(s^*) - H_{k+1, N_0}(s^*)\over
(p - 2) (1 - \log p/\log \underline d)}
+ R_1 + R_2
\cr
& \le 4 {\rm li}(N)
\sum_d {\sigma(d) C_{d N}\over \varphi(d) \log \underline d}
\sum_{\underline d^{1/s'}\le p<\underline d^{1/s}}
{A(s^*) - H_{k+1, N_0}(s^*)\over
(p - 2) (1 - \log p/\log \underline d)}
+ R_1 + R_2 + R_3,
\cr}$$
where we have used the fact that $A(s)-H_{k+1, N_0}(s)$ is increasing 
on $s$, and the notation
$$R_3 := 4 {\rm li}(N)
\sum_d {\sigma(d) C_{d N}\over \varphi(d) \log \underline d}
\sum_{\underline V^{1/s'}\le p<\underline d^{1/s'}}
{A(s^*) - H_{k+1, N_0}(s^*)\over
\varphi(p)(1-\log p/\log \underline d)}.$$
Applying the prime number theorem, an integration by parts shows that
$$\sum_{\underline d^{1/s'}\le p<\underline d^{1/s}}
{A(s^*) - H_{k+1, N_0}(s^*)\over
(p - 2) (1 - \log p/\log \underline d)}
= \int_{s-1}^{s'-1} {A(t) - H_{k+1, N_0}(t)\over t}\d t + O_{\delta, 
k}(\varepsilon).$$
Hence
$$S
\le \bigg\{\int_{s-1}^{s'-1} {A(t) - H_{k+1, N_0}(t)\over t}\d t + 
O_{\delta, k}(\varepsilon)\bigg\}
\Theta(N, \sigma)
+ R_1 + R_2 + R_3.
\leqno(3.16)$$

It remains to estimate $R_1, R_2, R_3$.
Observing that $\sigma(d)\not= 0\Rightarrow V/\Delta^i\le d<V$,
we have $\underline d^{1/s}\le \underline V^{1/s}\Delta^{i/s}$.
Thus
$\log ({\log \underline d^{1/s}/\log \alpha_r})
\le \log (1 + \log \Delta^{1+i/s}/\log(\underline V^{1/s}/\Delta))
\ll_{\delta, k} {\cal L}^{-5}$.
By using the prime number theorem and the previous estimate, we have
$$\leqalignno{R_1
& \ll \sum_d \sigma(d)
\sum_{\alpha_r\le p<\underline d^{1/s}} N/dp
& (3.17)
\cr
& \ll N {\cal L}^{-5} \sum_d \sigma(d)/\varphi(d)
\cr
& \ll_{\delta, k} \Theta(N, \sigma)/{\cal L}^3.
\cr}$$
Similarly we can show that
$$R_3
\ll_{\delta, k} \Theta(N, \sigma)/{\cal L}^3.
\leqno(3.18)$$
By the definition of $R_2$, we easily see that
$$R_2
\ll \sum_d \sigma(d)
\sum_{\alpha_0\le p<\alpha_r}
\sum_{(\underline{dp})^{1/s^*}\le p'<(\underline{dp})^{1/s_1}} N/dpp'.$$
Using a similar preceding argument, we can show that
$$\leqalignno{R_2
& \ll_{\delta, k} {N\over {\cal L}^4}
\sum_d {\sigma(d)\over \varphi(d)\log \underline d}
\sum_{\underline V^{1/s'}\le p<\underline V^{1/s}} {1\over p}
& (3.19)
\cr
& \ll_{\delta, k} {\Theta(N, \sigma)\over {\cal L}^3}.
\cr}$$
Combining (3.16)--(3.19), we obtain the desired upper bound, for $N\ge N_0$,
$$S
\le \bigg\{\int_{s-1}^{s'-1} {A(t) - H_{k+1, N_0}(t)\over t}\d t
+ O_{\delta, k}(\varepsilon)\bigg\}
\Theta(N, \sigma).
\leqno(3.20)$$

Inserting it in (3.14), estimating the first sum on the right-hand 
side of (3.14) by (3.9)
and noticing the relation
$$a(s') - a(s) = \int_s^{s'} A(t-1) \d t,$$
we find that, for $N\ge N_0(\varepsilon,\delta,k)$,
$$\Phi(N, \sigma, s)
\ge \Big\{a(s) + h_{k, N_0}(s')
+ \int_{s-1}^{s'-1} {H_{k+1, N_0}(t)\over t} \d t + O_{\delta, 
k}(\varepsilon)\Big\}
\Theta(N, \sigma).$$
 From the definition of $h_{k, N_0}(s)$,
we deduce that, for any $\varepsilon>0$ and for sufficiently large $N_0$,
$$h_{k, N_0}(s)
\ge h_{k, N_0}(s')
+ \int_{s-1}^{s'-1} {H_{k+1, N_0}(t)\over t} \d t
+ O_{\delta, k}(\varepsilon).$$
Taking $N_0\to \infty$ and then $\varepsilon\to 0$, we obtain
$$h_k(s)
\ge h_k(s')
+\int_{s-1}^{s'-1} {H_{k+1}(t)\over t} \d t$$
which implies the required inequality.
This completes the proof.
\hfill $\square$

\proclaim Corollary 1.
The function $H(s)$ is decreasing on $[1, 10]$.
The function $h(s)$ is increasing on $[1, 2]$ and is decreasing on $[2, 10]$.

\noindent{\sl Proof}.
According to the definition, we easily see that $H_{k, N_0}(s)$ is 
decreasing on $[1, 3]$
since $A(s) = 1$ for $1\le s\le 3$.
Thus $H(s)$ is also decreasing on $[1, 3]$.
When $3\le s\le 10$, the required result follows immediately from 
Propositions 2 and 1.

Similarly the definition of $h_{k, N_0}(s)$ and the fact that $a(s) = 
0$ for $1\le s\le 2$
show that $h(s)$ is increasing on $[1, 2]$.
Propositions 2 and 1 imply that $h(s)$ is decreasing on $[2, 10]$.
This concludes the proof.
\hfill $\square$

\smallskip

The central results in this section are Propositions 3 and 4 below.
Before stating it, it is necessary to introduce some notation.

Let $1\le s\le 3\le s'\le 5$ and $s\le \kappa_3\le \kappa_2\le \kappa_1\le s'$.
Define
$$\eqalign{
& \alpha_1 := \kappa_1-2,
\cr
& \alpha_4 := s'-s'/\kappa_2-1,
\cr
& \alpha_7 := s'-s'/\kappa_1-s'/\kappa_3,
\cr}
\qquad
\eqalign{
& \alpha_2 := s'-2,
\cr
& \alpha_5 := s'-s'/\kappa_3-1,
\cr
& \alpha_8 := s'-s'/\kappa_1-s'/\kappa_2,
\cr}
\qquad
\eqalign{
& \alpha_3 := s'-s'/s-1,
\cr
& \alpha_6 := s'-2s'/\kappa_2,
\cr
& \alpha_9 := \kappa_1-\kappa_1/\kappa_2-1.
\cr}$$
Let ${\bf 1}_{[a, b]}(t)$ be the characteristic function of the 
interval $[a, b]$.
We put
$$\qquad
\sigma(a, b, c)
:= \int_a^b \log{c\over t-1} {\d t\over t},
\qquad
\sigma_0(t) := {\sigma(3, t+2, t+1)\over 1-\sigma(3, 5, 4)}.$$

\smallskip

We can prove that $H(s)$ satisfies some functional inequalities.

\proclaim Proposition 3.
For $5\ge s'\ge 3\ge s\ge 2$ and $s' - s'/s\ge 2$, we have
$$H(s)\ge \Psi_1(s) + \int_1^3 H(t) \Xi_1(t; s) \d t,
\leqno(3.21)$$
where $\Psi_1(s)$ is defined as in Lemma 5.1 below and 
$\Xi_1(t; s) = \Xi_1(t; s, s')$ is given by
$$\eqalign{\Xi_1(t; s)
& := {\sigma_0(t)\over 2t} \log\bigg({16\over (s-1)(s'-1)}\bigg)
\cr
& \quad
+ {{\bf 1}_{[\alpha_2, 3]}(t)\over 2t}\log\bigg({(t+1)^2\over 
(s-1)(s'-1)}\bigg)
\cr
& \quad
+ {{\bf 1}_{[\alpha_3, \alpha_2]}(t)\over 2t}\log\bigg({t+1\over 
(s-1)(s'-1-t)}\bigg).
\cr}$$

\proclaim Proposition 4.
Let $5\ge s'\ge 3\ge s\ge 2$ and $s\le \kappa_3<\kappa_2<\kappa_1\le s'$ satisfy
$$s' - s'/s\ge 2,
\qquad
1\le \alpha_i\le 3\quad(1\le i\le 9),
\qquad
\alpha_1<\alpha_4,
\qquad
\alpha_5<\alpha_8.$$
Then we have
$$H(s)\ge \Psi_2(s) + \int_1^3 H(t) \Xi_2(t; s) \d t,
\leqno(3.22)$$
where
$\Psi_2(s)$ is defined as in Lemma 5.2 below, and
$\Xi_2(t; s) = \Xi_2(t; s, s', \kappa_1, \kappa_2, \kappa_3)$ is given by
$$\eqalign{\Xi_2(t; s)
& := {\sigma_0(t)\over 5t}
\log\bigg({1024\over (s-1)(s'-1)(\kappa_1-1)(\kappa_2-1)(\kappa_3-1)}\bigg)
\cr
& \quad
+ {{\bf 1}_{[\alpha_2, 3]}(t)\over 5t}
\log\bigg({(t+1)^5\over (s-1)(s'-1)(\kappa_1-1)(\kappa_2-1)(\kappa_3-1)}\bigg)
\cr
& \quad
+ {{\bf 1}_{[\alpha_9, \alpha_1]}(t)\over 5t}
\log\bigg({t+1\over (\kappa_2-1)(\kappa_1-1-t)}\bigg)
\cr
& \quad
+ {{\bf 1}_{[\alpha_5, \alpha_2]}(t)\over 5t}
\log\bigg({t+1\over (\kappa_3-1)(s'-1-t)}\bigg)
\cr
& \quad
+ {{\bf 1}_{[\alpha_3, \alpha_2]}(t)\over 5t} \log\bigg({t+1\over 
(s-1)(s'-1-t)}\bigg)
\cr
& \quad
+ {{\bf 1}_{[\alpha_1, \alpha_2]}(t)\over 5t}
\log\bigg({(t+1)^2\over (\kappa_1-1)(\kappa_2-1)}\bigg)
\cr
& \quad
+ {{\bf 1}_{[\alpha_7, \alpha_5]}(t)\over 5t(1-t/s')}
\log\bigg({s'^2\over (\kappa_1s'-s'-\kappa_1t)(\kappa_3s'-s'-\kappa_3t)}\bigg)
\cr
& \quad
+ {{\bf 1}_{[\alpha_5, \alpha_8]}(t)\over 5t(1-t/s')}
\log\bigg({s'(s'-1-t)\over \kappa_1s'-s'-\kappa_1t}\bigg)
\cr
& \quad
+ {{\bf 1}_{[\alpha_6, \alpha_8]}(t)\over 5t(1-t/s')}
\log\bigg({s'\over \kappa_2s'-s'-\kappa_2t}\bigg)
\cr
& \quad
+ {{\bf 1}_{[\alpha_8, \alpha_2]}(t)\over 5t(1-t/s')} \log(s'-1-t).
\cr}$$

We shall prove these two propositions in Section 6.
It is easy to see that $\Xi_i(t;s)$ is positive and
that for $s\in [1, 3)$ there exist parameters $s', \kappa_i$ such 
that $\Psi_i(s)>0$.
Therefore $H(s)>0$ for $s\in [1, 3)$ and 
then Proposition 2 implies that $h(s)>0$ for $s\in [1, 3)$.
In Sections 7 and 8, we shall give numeric solution of (3.20) and (3.21),
and prove Theorems 1 and 3.

\vskip 5mm

\noindent{\bf \S\ 4. Weighted inequalities for sieve function}

\medskip

The aim of this section is to present two weighted inequalities for 
sieve function.
The first is essentially due to Chen ([10], (23)).
The second is new,
which is not only much simpler than the third weighted inequality of Chen
([10], (64), (90) and (91)) but also more powerful.

\proclaim Lemma 4.1.
Let $1\le s<s'\le 10$.
For $N\ge 2$, $k\ge 0$ and $\sigma\in \ciU_k(N)$,
we have
$$\eqalign{2 \Phi(N, \sigma, s)
& \le \sum_d \sigma(d) (\Omega_1 - \Omega_2 + \Omega_3) + 
O_{\delta,k}(N^{1-\eta}),
\cr}$$
where $\eta = \eta(\delta,k)>0$ and $\Omega_i = \Omega_i(d)$ is given by
$$\eqalign{
& \Omega_1
:= 2 S({\cal A}_d; {\cal P}(d N), \underline d^{1/s'}),
\cr\noalign{\bigskip}
& \Omega_2
:= \sum_{\scriptstyle \underline d^{1/s'}\le p<\underline 
d^{1/s}\atop\scriptstyle (p, N)=1}
S({\cal A}_{dp}; {\cal P}(d N), \underline d^{1/s'}),
\cr
& \Omega_3
:= \mathop{\sum \, \sum \, \sum}_{\scriptstyle
\underline d^{1/s'}\le p_1<p_2<p_3<\underline 
d^{1/s}\atop\scriptstyle (p_1p_2p_3, N)=1}
S({\cal A}_{dp_1p_2p_3}; {\cal P}(d p_1 N), p_2).
\cr}$$

\noindent{\sl Proof}.
By the Buchstab identity, we have
$$\leqalignno{\qquad
2 S({\cal A}_d;  {\cal P}(d N), \underline d^{1/s})
& = \Omega_1
- 2 \sum_{\scriptstyle \underline d^{1/s'}\le p<\underline 
d^{1/s}\atop\scriptstyle (p, N)=1}
S({\cal A}_{dp}; {\cal P}(d N),  p),
& (4.1)
\cr
\sum_{\scriptstyle \underline d^{1/s'}\le p_1<\underline 
d^{1/s}\atop\scriptstyle (p_1, N)=1}
S({\cal A}_{dp_1}; {\cal P}(d N), p_1)
& = \Omega_0
+ \mathop{\sum \, \sum}_{\underline
d^{1/s'}\le p_1\le p_3<\underline d^{1/s}\atop\scriptstyle (p_1p_3, N)=1}
S({\cal A}_{dp_1p_3}; {\cal P}(d p_1 N), p_3),
& (4.2)
\cr
\sum_{\scriptstyle \underline d^{1/s'}\le p_3<\underline 
d^{1/s}\atop\scriptstyle (p_3, N)=1}
S({\cal A}_{dp_3}; {\cal P}(d N), p_3)
& = \Omega_2
- \mathop{\sum \, \sum}_{\underline
d^{1/s'}\le p_1<p_3<\underline d^{1/s}\atop\scriptstyle (p_1p_3, N)=1}
S({\cal A}_{dp_1p_3}; {\cal P}(d N), p_1),
& (4.3)
\cr}$$
where
$$\Omega_0
:= \sum_{\scriptstyle \underline d^{1/s'}\le p_1<\underline 
d^{1/s}\atop\scriptstyle (p_1, N)=1}
S({\cal A}_{dp_1}; {\cal P}(d N), \underline d^{1/s}).$$
Inserting (4.2)--(4.3) into (4.1),
dropping the term $\Omega_0$ (which is non-negative) 
and replacing $p_1\le p_3$ by $p_1<p_3$, we find that
$$2 S({\cal A}_d; {\cal P}(d N), \underline d^{1/s})
\le \Omega_1 - \Omega_2 + \Delta_1,$$
where
$$\leqalignno{\Delta_1
& := \mathop{\sum \, \sum}_{\scriptstyle
\underline d^{1/s'}\le p_1<p_3<\underline d^{1/s}\atop\scriptstyle 
(p_1p_3, N)=1}
\big\{S({\cal A}_{dp_1p_3}; {\cal P}(d N), p_1)
- S({\cal A}_{dp_1p_3}; {\cal P}(d p_1 N), p_3)\big\}
& (4.4)
\cr
& \,\, = \mathop{\sum \, \sum \, \sum}_{\underline
d^{1/s'}\le p_1\le p_2<p_3<\underline d^{1/s}\atop\scriptstyle 
(p_1p_2p_3, N)=1}
S({\cal A}_{dp_1p_2p_3}; {\cal P}(d p_1 N), p_2).
\cr}$$
By the inequality $S({\cal A}_{dp_1^2p_3}; {\cal P}(d N), p_1)\ll N/dp_1^2p_3$
and the fact that $\underline d\ge W_k$,
we easily see
$$\mathop{\sum \, \sum}_{\underline d^{1/s'}\le p_1<p_3<\underline d^{1/s}}
S({\cal A}_{dp_1^2p_3}; {\cal P}(d N), p_1)
\ll \mathop{\sum \, \sum}_{\underline d^{1/s'}\le p_1<p_3<\underline 
d^{1/s}} N/dp_1^2p_3
\ll_{\delta,k} N^{1-\eta}/d $$
for some $\eta = \eta(\delta,k)>0$.
Inserting it in (4.4), we obtain that
$$\Delta_1
= \Omega_3 + O_{\delta,k}(N^{1-\eta}/d).
\leqno(4.5)$$
Finally  we complete the proof with (3.5).
\hfill
$\square$

\proclaim Lemma 4.2.
Let $1\le s\le \kappa_3<\kappa_2<\kappa_1\le s'\le 10$.
For $N\ge 2$, $k\ge 0$ and $\sigma\in \ciU_k(N)$,
we have
$$\eqalign{5 \Phi(N, \sigma, s)
& \le \sum_d \sigma(d)
( \Gamma_1
- \Gamma_2
- \Gamma_3
- \Gamma_4
+ \Gamma_5
+ \cdots
+ \Gamma_{21})
+ O_{\delta, k} (N^{1-\eta}),
\cr}$$
where $\eta = \eta(\delta, k)>0$ and $\Gamma_i = \Gamma_i(d)$ is given by
$$\eqalign{
\Gamma_1
& := 4 S({\cal A}_d; {\cal P}(d N), \underline d^{1/s'})
+ S({\cal A}_d; {\cal P}(d N), \underline d^{1/\kappa_1}),
\cr\noalign{\bigskip}
\Gamma_2
& := \sum_{\scriptstyle \underline d^{1/s'}\le p<\underline d^{1/s}
\atop\scriptstyle (p, N)=1}
S({\cal A}_{dp}; {\cal P}(d N), \underline d^{1/s'}),
\cr
\Gamma_3
& := \sum_{\scriptstyle \underline d^{1/s'}\le p<\underline d^{1/\kappa_2}
\atop\scriptstyle (p, N)=1}
S({\cal A}_{dp}; {\cal P}(d N), \underline d^{1/s'}),
\cr
\Gamma_4
& := \sum_{\scriptstyle \underline d^{1/s'}\le p<\underline d^{1/\kappa_3}
\atop\scriptstyle (p, N)=1}
S({\cal A}_{dp}; {\cal P}(d N), \underline d^{1/s'}),
\cr
\Gamma_5
& := \mathop{\sum \, \sum}_{\scriptstyle \underline d^{1/s'}\le 
p_1<p_2<\underline d^{1/\kappa_2}
\atop\scriptstyle (p_1p_2, N)=1}
S({\cal A}_{dp_1p_2}; {\cal P}(d N), \underline d^{1/s'}),
\cr
\Gamma_6
& := \mathop{\sum\,\,\sum}_{\scriptstyle \underline d^{1/s'}\le 
p_1<\underline d^{1/\kappa_1},\,
\underline d^{1/\kappa_2}\le p_2<\underline d^{1/\kappa_3}
\atop\scriptstyle (p_1p_2, N)=1}
S({\cal A}_{dp_1p_2}; {\cal P}(d N), \underline d^{1/s'}),
\cr
\Gamma_7
& :=  \mathop{\sum \, \sum}_{\scriptstyle \underline d^{1/s'}\le 
p_1<p_2<\underline d^{1/\kappa_1}
\atop\scriptstyle (p_1p_2, N)=1}
S({\cal A}_{dp_1p_2}; {\cal P}(d N), p_1),
\cr
\Gamma_8
& := \mathop{\sum\,\,\sum}_{\scriptstyle
\underline d^{1/s'}\le p_1<\underline d^{1/\kappa_1}\le 
p_2<\underline d^{1/\kappa_2}
\atop\scriptstyle (p_1p_2, N)=1}
S({\cal A}_{dp_1p_2}; {\cal P}(d N), p_1),
\cr
\Gamma_9
& := \mathop{\sum\,\,\sum\,\,\sum}_{\scriptstyle
\underline d^{1/\kappa_1}\le p_1<p_2<p_3<\underline d^{1/\kappa_3}
\atop\scriptstyle (p_1p_2p_3, N)=1}
S({\cal A}_{dp_1p_2p_3}; {\cal P}(d N),  p_2),
\cr
\Gamma_{10}
& := \mathop{\sum\,\,\sum\,\,\sum}_{\scriptstyle
\underline d^{1/\kappa_1}\le p_1<p_2<\underline d^{1/\kappa_2}\le 
p_3<\underline d^{1/s}
\atop\scriptstyle (p_1p_2p_3, N)=1}
S({\cal A}_{dp_1p_2p_3}; {\cal P}(d N),  p_2),
\cr
\Gamma_{11}
& := \mathop{\sum\,\,\sum\,\,\sum}_{\scriptstyle
\underline d^{1/\kappa_1}\le p_1<\underline d^{1/\kappa_2}\le 
p_2<p_3<\underline d^{1/\kappa_3}
\atop\scriptstyle (p_1p_2p_3, N)=1}
S({\cal A}_{dp_1p_2p_3}; {\cal P}(d N),  p_2),
\cr
\Gamma_{12}
& := \mathop{\sum\,\,\sum\,\,\sum}_{\scriptstyle
\underline d^{1/s'}\le p_1<p_2<\underline d^{1/\kappa_1},\,
\underline d^{1/\kappa_3}\le p_3<\underline d^{1/s}
\atop\scriptstyle (p_1p_2p_3, N)=1}
S({\cal A}_{dp_1p_2p_3}; {\cal P}(d N),  p_2),
\cr
\Gamma_{13}
& := \mathop{\sum\,\,\sum\,\,\sum}_{\scriptstyle
\underline d^{1/s'}\le p_1<\underline d^{1/\kappa_1}\le 
p_2<\underline d^{1/\kappa_2}
\le p_3<\underline d^{1/s}
\atop\scriptstyle (p_1p_2p_3, N)=1}
S({\cal A}_{dp_1p_2p_3}; {\cal P}(d N),  p_2),
\cr
\Gamma_{14}
& := \mathop{\sum\,\,\sum\,\,\sum}_{\scriptstyle
\underline d^{1/s'}\le p_1<\underline d^{1/\kappa_1}, \,
\underline d^{1/\kappa_2}\le p_2<p_3<\underline d^{1/s}
\atop\scriptstyle (p_1p_2p_3, N)=1}
S({\cal A}_{dp_1p_2p_3}; {\cal P}(d N),  p_2),
\cr}$$
$$\eqalign{
\Gamma_{15}
& := \mathop{\sum\,\,\sum\,\,\sum}_{\scriptstyle \underline d^{1/\kappa_1}
\le p_1<\underline d^{1/\kappa_2}\le p_2<\underline d^{1/\kappa_3}\le 
p_3<\underline d^{1/s}
\atop\scriptstyle (p_1p_2p_3, N)=1}
S({\cal A}_{dp_1p_2p_3}; {\cal P}(d N),  p_2),
\cr
\Gamma_{16}
& := \mathop{\sum\,\,\sum\,\,\sum\,\,\sum}_{\scriptstyle
\underline d^{1/\kappa_2}\le p_1<p_2<p_3<p_4<\underline d^{1/\kappa_3}
\atop\scriptstyle (p_1p_2p_3p_4, N)=1}
S({\cal A}_{dp_1p_2p_3p_4}; {\cal P}(d N), p_3),
\cr
\Gamma_{17}
& := \mathop{\sum\,\,\sum\,\,\sum\,\,\sum}_{\scriptstyle
\underline d^{1/\kappa_2}\le p_1<p_2<p_3<\underline d^{1/\kappa_3}\le 
p_4<\underline d^{1/s}
\atop\scriptstyle (p_1p_2p_3p_4, N)=1}
S({\cal A}_{dp_1p_2p_3p_4}; {\cal P}(d N), p_3),
\cr
\Gamma_{18}
& := \mathop{\sum\,\,\sum\,\,\sum\,\,\sum}_{\scriptstyle
\underline d^{1/\kappa_2}\le p_1<p_2<\underline d^{1/\kappa_3}\le 
p_3<p_4<\underline d^{1/s}
\atop\scriptstyle (p_1p_2p_3p_4, N)=1}
S({\cal A}_{dp_1p_2p_3p_4}; {\cal P}(d N), p_3),
\cr
\Gamma_{19}
& := \mathop{\sum\,\,\sum\,\,\sum\,\,\sum}_{\scriptstyle
\underline d^{1/\kappa_1}\le p_1<\underline d^{1/\kappa_2}, \,
\underline d^{1/\kappa_3}\le p_2<p_3<p_4<\underline d^{1/s}
\atop\scriptstyle (p_1p_2p_3p_4, N)=1}
S({\cal A}_{dp_1p_2p_3p_4}; {\cal P}(d N),  p_3),
\cr
\Gamma_{20}
& := \mathop{\sum\,\,\sum\,\,\sum\,\,\sum\,\,\sum}_{\scriptstyle
\underline d^{1/\kappa_2}\le p_1<\underline d^{1/\kappa_3}\le 
p_2<p_3<p_4<p_5<\underline d^{1/s}
\atop\scriptstyle (p_1p_2p_3p_4p_5, N)=1}
S({\cal A}_{dp_1p_2p_3p_4p_5}; {\cal P}(d N), p_4),
\cr
\Gamma_{21}
& := \mathop{\sum\,\,\sum\,\,\sum\,\,\sum\,\,\sum\,\,\sum}_{\scriptstyle
\underline d^{1/\kappa_3}\le p_1<p_2<p_3<p_4<p_5<p_6<\underline d^{1/s}
\atop\scriptstyle (p_1p_2p_3p_4p_5p_6, N)=1}
S({\cal A}_{dp_1p_2p_3p_4p_5p_6}; {\cal P}(d N), p_5).
\cr}$$

\noindent{\sl Proof}.
Let $S := S({\cal A}_d;  {\cal P}(d N), \underline d^{1/s})$.
By using the Buchstab identity, we have
$$\leqalignno{\qquad
2 S
& = 2 S({\cal A}_d;  {\cal P}(d N), \underline d^{1/s'})
- \sum_{\scriptstyle \underline d^{1/s'}\le p<\underline 
d^{1/s}\atop\scriptstyle (p, N)=1}
S({\cal A}_{dp}; {\cal P}(d N),  p)
& (4.6)
\cr
& \quad
- \Gamma_3
+ \mathop{\sum \, \sum}_{\scriptstyle
\underline d^{1/s'}\le p_1<p_2<\underline 
d^{1/\kappa_2}\atop\scriptstyle (p_1p_2, N)=1}
S({\cal A}_{dp_1p_2}; {\cal P}(d N), p_1)
- \sum_{\scriptstyle \underline d^{1/\kappa_2}\le p<\underline 
d^{1/s}\atop\scriptstyle (p, N)=1}
S({\cal A}_{dp}; {\cal P}(d N),  p)
\cr
& =: 2 S({\cal A}_d;  {\cal P}(d N), \underline d^{1/s'})
- E_1 - \Gamma_3 + D_1' - E_2.
\cr}$$

We can also write, always by Buchstab's identity,
$$\leqalignno{S
& = S({\cal A}_d;  {\cal P}(d N), \underline d^{1/s'})
- \sum_{\scriptstyle \underline d^{1/s'}\le p<\underline d^{1/\kappa_3}
\atop\scriptstyle (p, N)=1}
S({\cal A}_{dp}; {\cal P}(d N),  p)
& (4.7)
\cr
& \quad
- \sum_{\scriptstyle \underline d^{1/\kappa_3}\le p<\underline d^{1/s}
\atop\scriptstyle (p, N)=1}
S({\cal A}_{dp}; {\cal P}(d N),  p).
\cr}$$
But we have
$$\eqalign{
\sum_{\scriptstyle \underline d^{1/s'}\le p<\underline 
d^{1/\kappa_3}\atop\scriptstyle (p, N)=1}
S({\cal A}_{dp}; {\cal P}(d N),  p)
& = \Gamma_4
- \Gamma_7 - \Gamma_8
- \mathop{\sum\,\, \sum}_{\scriptstyle \underline d^{1/\kappa_1}\le 
p_1<p_2<\underline d^{1/\kappa_3}
\atop\scriptstyle (p_1p_2, N)=1}
S({\cal A}_{dp_1p_2}; {\cal P}(d N),  p_1)
\cr
& - \Gamma_6
+ \mathop{\sum\,\,\sum\,\,\sum}_{\scriptstyle
\underline d^{1/s'}\le p_1<p_2<\underline d^{1/\kappa_1}, \,
\underline d^{1/\kappa_2}\le p_3<\underline d^{1/\kappa_3}
\atop\scriptstyle (p_1p_2p_3, N)=1}
S({\cal A}_{dp_1p_2p_3}; {\cal P}(d N),  p_1).
\cr}$$
Inserting these relations into (4.7), it yields that
$$\leqalignno{\qquad\qquad
S
& = S({\cal A}_d;  {\cal P}(d N), \underline d^{1/s'})
- \Gamma_4
+ \Gamma_6
+ \Gamma_7
+ \Gamma_8
& (4.8)
\cr\noalign{\medskip}
& \quad
- \sum_{\scriptstyle
\underline d^{1/\kappa_3}\le p<\underline d^{1/s}\atop\scriptstyle (p, N)=1}
S({\cal A}_{dp}; {\cal P}(d N),  p)
+ \mathop{\sum\,\, \sum}_{\scriptstyle \underline d^{1/\kappa_1}\le 
p_1<p_2<\underline d^{1/\kappa_3}
\atop\scriptstyle (p_1p_2, N)=1}
S({\cal A}_{dp_1p_2}; {\cal P}(d N),  p_1)
\cr
& \quad
- \mathop{\sum\,\,\sum\,\,\sum}_{\scriptstyle
\underline d^{1/s'}\le p_1<p_2<\underline d^{1/\kappa_1}, \,
\underline d^{1/\kappa_2}\le p_3<\underline d^{1/\kappa_3}
\atop\scriptstyle (p_1p_2p_3, N)=1}
S({\cal A}_{dp_1p_2p_3}; {\cal P}(d N), p_1)
\cr
& =: S({\cal A}_d;  {\cal P}(d N), \underline d^{1/s'})
- \Gamma_4
+ \Gamma_6
+ \Gamma_7
+ \Gamma_8
- E_3
+ D_2
- E_4.
\cr}$$

Similar to (4.8), we can prove that
$$\leqalignno{\qquad
S
& = S({\cal A}_d;  {\cal P}(d N), \underline d^{1/s'})
- \Gamma_2
+ \Gamma_5
- \mathop{\sum\,\,\sum\,\,\sum}_{\scriptstyle
\underline d^{1/s'}\le p_1<p_2<p_3<\underline d^{1/\kappa_2}
\atop\scriptstyle (p_1p_2p_3, N)=1}
S({\cal A}_{dp_1p_2p_3}; {\cal P}(d N),  p_1)
& (4.9)
\cr
& \quad
+ \Big\{
\mathop{\sum\,\,\sum}_{\scriptstyle
\underline d^{1/s'}\le p_1<\underline d^{1/\kappa_2}\le p_2<\underline d^{1/s}
\atop\scriptstyle (p_1p_2, N)=1}
+ \mathop{\sum\,\,\sum}_{\scriptstyle
\underline d^{1/\kappa_2}\le p_1<p_2<\underline d^{1/s}
\atop\scriptstyle (p_1p_2, N)=1}
\Big\}
S({\cal A}_{dp_1p_2}; {\cal P}(d N),  p_1)
\cr
& =: S({\cal A}_d;  {\cal P}(d N), \underline d^{1/s'})
- \Gamma_2
+ \Gamma_5
- E_5
+ D_1''
+ D_1'''.
\cr}$$

Finally we write
$$S
= S({\cal A}_d;  {\cal P}(d N), \underline d^{1/\kappa_1})
-  \Big\{
\sum_{\scriptstyle \underline d^{1/\kappa_1}\le p<\underline d^{1/\kappa_3}
\atop\scriptstyle (p, N)=1}
+  \sum_{\scriptstyle \underline d^{1/\kappa_3}\le p<\underline 
d^{1/s}\atop\scriptstyle (p, N)=1}
\Big\}
S({\cal A}_{dp}; {\cal P}(d N),  p).
\leqno(4.10)$$
For $p_1<\underline d^{1/\kappa_3}<\underline d^{1/s}$, we have
$$S({\cal A}_{dp}; {\cal P}(d N),  p)
\ge \sum_{\scriptstyle p\le p_1<\underline 
d^{1/\kappa_3}\atop\scriptstyle (p_1, N)=1}
S({\cal A}_{dpp_1}; {\cal P}(d N),  p_1)
+ \sum_{\scriptstyle \underline d^{1/\kappa_3}\le p_1<\underline 
d^{1/s}\atop\scriptstyle (p_1, N)=1}
S({\cal A}_{dpp_1}; {\cal P}(d N), p_1).$$
This implies
$$\eqalign{\sum_{\scriptstyle \underline d^{1/\kappa_1}\le 
p<\underline d^{1/\kappa_3}
\atop\scriptstyle (p, N)=1}
S({\cal A}_{dp}; {\cal P}(d N),  p)
& \ge \mathop{\sum\,\,\sum}_{\scriptstyle \underline 
d^{1/\kappa_1}\le p_1<p_2<\underline d^{1/\kappa_3}
\atop\scriptstyle (p_1p_2, N)=1}
S({\cal A}_{dp_1p_2}; {\cal P}(d N),  p_2)
\cr
& \quad
+ \mathop{\sum\,\,\sum}_{\scriptstyle
\underline d^{1/\kappa_1}\le p_1<\underline d^{1/\kappa_3}\le 
p_2<\underline d^{1/s}
\atop\scriptstyle (p_1p_2, N)=1}
S({\cal A}_{dp_1p_2}; {\cal P}(d N),  p_2).
\cr}$$
Inserting it into (4.10), we obtain
$$\leqalignno{\qquad
S
& \le S({\cal A}_d;  {\cal P}(d N), \underline d^{1/\kappa_1})
- \mathop{\sum\,\,\sum}_{\scriptstyle \underline d^{1/\kappa_1}\le 
p_1<p_2<\underline d^{1/\kappa_3}
\atop\scriptstyle (p_1p_2, N)=1}
S({\cal A}_{dp_1p_2}; {\cal P}(d N),  p_2)
& (4.11)
\cr
& \quad
- \mathop{\sum\,\,\sum}_{\scriptstyle
\underline d^{1/\kappa_1}\le p_1<\underline d^{1/\kappa_3}\le 
p_2<\underline d^{1/s}
\atop\scriptstyle (p_1p_2, N)=1}
S({\cal A}_{dp_1p_2}; {\cal P}(d N),  p_2)
-  \sum_{\scriptstyle \underline d^{1/\kappa_3}\le p<\underline 
d^{1/s}\atop\scriptstyle (p, N)=1}
S({\cal A}_{dp}; {\cal P}(d N),  p)
\cr
& =: S({\cal A}_d;  {\cal P}(d N), \underline d^{1/\kappa_1}) - E_6 - 
E_7 - E_3.
\cr}$$

Now by adding up the inequalities (4.6), (4.8), (4.9) and (4.11)
and by noticing the estimate $D_2 - E_6\le \Gamma_9 
+ O_{\delta, k}(N^{1-\eta}/d)$,
we get
$$5 S
\le \Gamma_1
- \Gamma_2
- \Gamma_3
- \Gamma_4
+ \Gamma_5
+ \cdots
+ \Gamma_9
+ \Delta_2
+ O_{\delta, k}(N^{1-\eta}/d)
\leqno(4.12)$$
where
$$\Delta_2 := D_1 - E_1 - E_2 - 2E_3 - E_4 - E_5 - E_7$$
and
$$D_1 := D_1' + D_1'' + D_1'''
= \mathop{\sum\,\,\sum}_{\scriptstyle
\underline d^{1/s'}\le p_1<p_2<\underline d^{1/s}
\atop\scriptstyle (p_1p_2, N)=1}
S({\cal A}_{dp_1p_2}; {\cal P}(d N),  p_1).$$

Clearly we have
$$E_1
\ge \mathop{\sum\sum}_{\scriptstyle \underline d^{1/s'}\le p_1<p_2<\underline 
d^{1/s}\atop\scriptstyle (p_1p_2, N)=1}
S({\cal A}_{dp_1p_2}; {\cal P}(d N), p_2).$$
Thus an application of Bechstab's identity gives us
$$D_1-E_1
\le \mathop{\sum\,\,\sum\,\,\sum}_{\scriptstyle
\underline d^{1/s'}\le p_1<p_2<p_3<\underline 
d^{1/s}\atop\scriptstyle (p_1p_2p_3, N)=1}
S({\cal A}_{dp_1p_2p_3}; {\cal P}(d N),  p_2) + O_{\delta, k}(N^{1-\eta}/d).$$
From this, we can deduce
$$\eqalign{D_1 - E_1 - E_5
& \le \mathop{\sum\,\,\sum\,\,\sum}_{\scriptstyle
\underline d^{1/s'}\le p_1<p_2<p_3<\underline 
d^{1/s}\atop\scriptstyle (p_1p_2p_3, N)=1}
S({\cal A}_{dp_1p_2p_3}; {\cal P}(d N),  p_2)
\cr
& \quad
- \mathop{\sum\,\,\sum\,\,\sum}_{\scriptstyle
\underline d^{1/s'}\le p_1<p_2<p_3<\underline d^{1/\kappa_2}
\atop\scriptstyle (p_1p_2p_3, N)=1}
S({\cal A}_{dp_1p_2p_3}; {\cal P}(d N),  p_2)
+ O_{\delta, k}(N^{1-\eta}/d)
\cr
& = \mathop{\sum\,\,\sum\,\,\sum}_{\scriptstyle
\underline d^{1/s'}\le p_1<p_2<\underline d^{1/\kappa_2}\le 
p_3<\underline d^{1/s}
\atop\scriptstyle (p_1p_2p_3, N)=1}
S({\cal A}_{dp_1p_2p_3}; {\cal P}(d N),  p_2)
\cr
& \quad
+ \mathop{\sum\,\,\sum\,\,\sum}_{\scriptstyle
\underline d^{1/s'}\le p_1<\underline d^{1/\kappa_2}\le 
p_2<p_3<\underline d^{1/s}
\atop\scriptstyle (p_1p_2p_3, N)=1}
S({\cal A}_{dp_1p_2p_3}; {\cal P}(d N),  p_2)
\cr
& \quad
+ \mathop{\sum\,\,\sum\,\,\sum}_{\scriptstyle
\underline d^{1/\kappa_2}\le p_1<p_2<p_3<\underline d^{1/s}
\atop\scriptstyle (p_1p_2p_3, N)=1}
S({\cal A}_{dp_1p_2p_3}; {\cal P}(d N),  p_2)
 + O_{\delta, k}(N^{1-\eta}/d)
\cr
& =: D_3 + D_4 + D_5 + O_{\delta, k}(N^{1-\eta}/d).
\cr}$$

We have
$$D_3 - E_4
\le \Gamma_{10} + \Gamma_{12} + \Gamma_{13}.$$

By splitting $D_4$ into 4 subsums, we have
$$D_4 = \Gamma_{11} + \Gamma_{14} + \Gamma_{15} 
+ \mathop{\sum\,\,\sum\,\,\sum}_{\scriptstyle
\underline d^{1/\kappa_1}\le p_1<\underline d^{1/\kappa_2}, \,
\underline d^{1/\kappa_3}\le p_2<p_3<\underline d^{1/s}
\atop\scriptstyle (p_1p_2p_3, N)=1}
S({\cal A}_{dp_1p_2p_3}; {\cal P}(d N),  p_2).$$
Similarly by splitting $E_7$ into 2 subsums, we have
$$\eqalign{E_7 
& = E_8 
+ \mathop{\sum\,\,\sum}_{\scriptstyle
\underline d^{1/\kappa_1}\le p_1<\underline d^{1/\kappa_2}, \,
\underline d^{1/\kappa_3}\le p_2<\underline d^{1/s}
\atop\scriptstyle (p_1p_2, N)=1}
S({\cal A}_{dp_1p_2}; {\cal P}(d N),  p_2)
\cr
& \ge E_8 
+ \mathop{\sum\,\,\sum\,\,\sum}_{\scriptstyle
\underline d^{1/\kappa_1}\le p_1<\underline d^{1/\kappa_2}, \,
\underline d^{1/\kappa_3}\le p_2<p_3<\underline d^{1/s}
\atop\scriptstyle (p_1p_2p_3, N)=1}
S({\cal A}_{dp_1p_2p_3}; {\cal P}(d N),  p_3),
\cr}$$
where
$$E_8 := \mathop{\sum\,\,\sum}_{\scriptstyle
\underline d^{1/\kappa_2}\le p_1<\underline d^{1/\kappa_3}\le 
p_2<\underline d^{1/s}
\atop\scriptstyle (p_1p_2, N)=1}
S({\cal A}_{dp_1p_2}; {\cal P}(d N),  p_2).$$
By noticing that
$$\eqalign{
& \mathop{\sum\,\,\sum\,\,\sum}_{\scriptstyle
\underline d^{1/\kappa_1}\le p_1<\underline d^{1/\kappa_2}, \,
\underline d^{1/\kappa_3}\le p_2<p_3<\underline d^{1/s}
\atop\scriptstyle (p_1p_2p_3, N)=1}
\big\{
S({\cal A}_{dp_1p_2p_3}; {\cal P}(d N),  p_2)
- S({\cal A}_{dp_1p_2p_3}; {\cal P}(d N),  p_3)
\big\}
\cr
& \qquad\qquad
= \Gamma_{19} + O_{\delta, k}(N^{1-\eta}/d),
\cr}$$
we can deduce
$$D_4 - E_7
\le \Gamma_{11} + \Gamma_{14} + \Gamma_{15} + \Gamma_{19} - E_8 + 
O_{\delta, k}(N^{1-\eta}/d).$$

Since
$$E_2
\ge \mathop{\sum\,\,\sum\,\,\sum}_{\scriptstyle
\underline d^{1/\kappa_2}\le p_1<p_2<p_3<\underline d^{1/s}
\atop\scriptstyle (p_1p_2p_3, N)=1}
S({\cal A}_{dp_1p_2p_3}; {\cal P}(d N), p_3),$$
we have
$$\eqalign{D_5 - E_2
& \le \mathop{\sum\,\,\sum\,\,\sum\,\,\sum}_{\scriptstyle
\underline d^{1/\kappa_2}\le p_1<p_2<p_3<p_4<\underline d^{1/s}
\atop\scriptstyle (p_1p_2p_3p_4, N)=1}
S({\cal A}_{dp_1p_2p_3p_4}; {\cal P}(d N), p_3)
+ O_{\delta, k}(N^{1-\eta}/d)
\cr
& =: D_6 + O_{\delta, k}(N^{1-\eta}/d).
\cr}$$

Similarly$$\eqalign{E_3
& \ge \mathop{\sum\,\,\sum\,\,\sum\,\,\sum}_{\scriptstyle
\underline d^{1/\kappa_3}\le p_1\le p_2<p_3<p_4<\underline d^{1/s}
\atop\scriptstyle (p_1p_2p_3p_4, N)=1}
S({\cal A}_{dp_1p_2p_3p_4}; {\cal P}(d N),  p_4)
=: E_3',
\cr
E_3
& \ge \mathop{\sum\,\,\sum\,\,\sum\,\,\sum\,\,\sum}_{\scriptstyle
\underline d^{1/\kappa_3}\le p_1\le p_2<p_3<p_4<p_5<\underline d^{1/s}
\atop\scriptstyle (p_1p_2p_3p_4p_5, N)=1}
S({\cal A}_{dp_1p_2p_3p_4p_5}; {\cal P}(d N),  p_5)
=: E_3'',
\cr
E_8
& \ge \mathop{\sum\,\,\sum\,\,\sum\,\,\sum}_{\scriptstyle
\underline d^{1/\kappa_2}\le p_1<\underline d^{1/\kappa_3}\le p_2\le 
p_3<p_4<\underline d^{1/s}
\atop\scriptstyle (p_1p_2p_3p_4, N)=1}
S({\cal A}_{dp_1p_2p_3p_4}; {\cal P}(d N),  p_4)
=: E_8',
\cr
D_6
& = \Gamma_{16} + \Gamma_{17} + \Gamma_{18}
+ \mathop{\sum\,\,\sum\,\,\sum\,\,\sum}_{\scriptstyle
\underline d^{1/\kappa_2}\le p_1<\underline d^{1/\kappa_3}\le 
p_2<p_3<p_4<\underline d^{1/s}
\atop\scriptstyle (p_1p_2p_3p_4, N)=1}
S({\cal A}_{dp_1p_2p_3p_4}; {\cal P}(d N), p_3)
\cr
& \quad
+ \mathop{\sum\,\,\sum\,\,\sum\,\,\sum}_{\scriptstyle
\underline d^{1/\kappa_3}\le p_1<p_2<p_3<p_4<\underline d^{1/s}
\atop\scriptstyle (p_1p_2p_3p_4, N)=1}
S({\cal A}_{dp_1p_2p_3p_4}; {\cal P}(d N), p_3)
\cr
& =: \Gamma_{16} + \Gamma_{17} + \Gamma_{18} + D_6' + D_6''.
\cr}$$
Since
$$D_6' - E_8' = \Gamma_{20} + O_{\delta, k}(N^{1-\eta}/d)$$
and
$$\eqalign{D_6'' - E_3' - E_3''
& = \mathop{\sum\,\,\sum\,\,\sum\,\,\sum\,\,\sum}_{\scriptstyle
\underline d^{1/\kappa_3}\le p_1<p_2<p_3<p_4<p_5<\underline d^{1/s}
\atop\scriptstyle (p_1p_2p_3p_4p_5, N)=1}
S({\cal A}_{dp_1p_2p_3p_4p_5}; {\cal P}(d N), p_4)
\cr
& \quad
- \mathop{\sum\,\,\sum\,\,\sum\,\,\sum\,\,\sum}_{\scriptstyle
\underline d^{1/\kappa_3}\le p_1<p_2<p_3<p_4<p_5<\underline d^{1/s}
\atop\scriptstyle (p_1p_2p_3p_4p_5, N)=1}
S({\cal A}_{dp_1p_2p_3p_4p_5}; {\cal P}(d N), p_5)
+ O_{\delta, k}(N^{1-\eta}/d)
\cr
& = \Gamma_{21} + O_{\delta, k}(N^{1-\eta}/d),
\cr}$$
we have
$$\eqalign{D_6 - 2 E_3 - E_8
& \le \Gamma_{16} + \Gamma_{17} + \Gamma_{18} + D_6' + D_6'' - E_3' - 
E_3'' - E_8'
\cr
& = \Gamma_{16} + \Gamma_{17} + \Gamma_{18} + \Gamma_{20} + \Gamma_{21}
+ O_{\delta, k}(N^{1-\eta}/d).
\cr}$$

Combining these estimations leads to the following inequalities
$$\eqalign{\Delta_2
& \le D_3 + D_4 + D_5 - E_2 - 2E_3 - E_4 - E_7
\cr
& \le \Gamma_{10} + \cdots + \Gamma_{15} + \Gamma_{19} + D_6 - 2E_3 - E_8
+ O_{\delta, k}(N^{1-\eta}/d)
\cr
& \le \Gamma_{10} + \cdots + \Gamma_{21} + O_{\delta, k}(N^{1-\eta}/d).
\cr}
\leqno(4.13)$$
Now the desired result follows from (4.12) and (4.13).
\hfill
$\square$

\vskip 5mm

\noindent{\bf \S\ 5. Functional inequalities between $H(s)$ and $h(s)$}

\medskip

In this section, we start from two weighted inequalities for the sieve function
to deduce two functional inequalities between $H(s)$ and $h(s)$.
They will be used to prove Propositions 3 and 4 in the next section.

\proclaim Lemma 5.1.
For $5\ge s'\ge 3\ge s\ge 2$ and $s'-s'/s\ge 2$,
we have
$$H(s)
\ge \Psi_1(s)
+ {1\over 2} \int_{1-1/s}^{1-1/s'} {h(s't)\over t(1-t)} \d t + H(s'),$$
where $\Psi_1(s) = \Psi_1(s, s')$ is given by
$$\Psi_1(s)
:=
- \int_2^{s'-1} {\log(t-1)\over t}\d t
+ {1\over 2} \int_{1-1/s}^{1-1/s'} {\log(s't-1)\over t(1-t)} \d t
- I_1(s)$$
and $I_1(s) = I_1(s, s')$ is given by
$$I_1(s)
:= \max_{\phi\ge 2} \mathop{\int\int\int}_{1/s'\le t\le u\le v\le 1/s}
\omega\bigg({\phi - t - u - v\over u}\bigg)
{\d t \d u \d v\over t u^2 v}.$$

\noindent{\sl Proof}.
Our starting point is the inequality in Lemma 4.1.
We need to estimate all terms in the right-hand side of this inequality.

Firstly, (3.3) of Lemma 3.2 gives us
$$\sum_d \sigma(d) \Omega_1
\le 2 \big\{A(s') - H_{k, N_0}(s')\big\} \Theta(N, \sigma).
\leqno(5.1)$$

Secondly, by an argument similar to the proof of (3.20),
we can prove, for any $\varepsilon>0$ and $N\ge N_0(\varepsilon, \delta, k)$,
$$\sum_d \sigma(d) \Omega_2
\ge \bigg\{
\int_{1-1/s}^{1-1/s'} {a(s't) + h_{k+1, N_0}(s't)\over t(1-t)} \d t - 
\varepsilon
\bigg\}
\Theta(N, \sigma).
\leqno(5.2)$$

Finally we apply the switching principle to estimate $\sum_d 
\sigma(d) \Omega_3$.
For this, we introduce
$${\cal E} := \{e : e = d n p_1 p_2, \, \sigma(d)\not= 0, \,
(n, p_1, p_2) \,\, \hbox{satisfies} \,\, (5.3) \,\, \hbox{below}\},$$
where
$$\underline d^{1/s'}\le p_1<p_2<\underline d^{1/s},
\quad
(p_1 p_2, d N) = 1,
\quad
n\le N/d p_1 p_2^2,
\quad
(n, \, N P(p_2)/p_1) = 1;
\leqno(5.3)$$
and
$${\cal B}
:= \{b : b = N - e p_3, \, e\in {\cal E}, \, p_2 < p_3 \le \kappa(d, e)\},$$
where $\kappa(d, e) := \min\{N/e, \, \underline d^{1/s}\}$.
The set ${\cal B}$ is a multiset and an element $b$ may occur more than once.

Clearly $\sum_d \sigma(d) \Omega_3$ does not exceed the number of 
primes in the set ${\cal B}$.
Thus
$$\sum_d \sigma(d) \Omega_3
\le S({\cal B}; {\cal P}(N), Q^{1/2}) + O(Q^{1/2}).$$
In the set ${\cal E}$, $d$ is not determined uniquely by $e$.
This causes technical difficulty.
In order to avoid it, we define ${\cal E}'$ and ${\cal B}'$, similar 
to ${\cal E}$ and ${\cal B}$, with the condition
 $(n, \, N P(p_2)/p_1) = 1$ replaced by $(n, \, d N P(p_2)/p_1) = 1$ and
${\cal E}$ by ${\cal E}'$ respectively.
Obviously
the difference $S({\cal B}; {\cal P}(N), Q^{1/2}) - S({\cal B}'; 
{\cal P}(N), Q^{1/2})$ is
$$\le \sum_d \sigma(d)
\mathop{\sum \,\,\, \sum \,\,\, \sum}_{\underline d^{1/s'}\le 
p_1<p_2<p_3<\underline d^{1/s}}
\sum_{\scriptstyle n\le N/d p_1 p_2 p_3\atop{\atop\scriptstyle (n, d) = 1}} 1
\ll_{\delta, k} N^{1-\eta}$$
where $\eta = \eta(\delta, k)>0$.
Hence
$$\sum_d \sigma(d) \Omega_3
\le S({\cal B}'; {\cal P}(N), Q^{1/2}) + O_{\delta, k}(N^{1-\eta}).
\leqno(5.4)$$
In order to estimate $S({\cal B}'; {\cal P}(N), Q^{1/2})$,
we use Theorem 5.2 in [14] with
$$X :=
\sum_{e\in {\cal E}'} \sum_{p_2<p_3\le \kappa(d, e)} 1,
\qquad
w(q) = \cases{
q/\varphi(q) & if $\mu(q)^2 = (q, N) = 1$
\cr\noalign{\smallskip}
0            & otherwise
\cr}$$
to write
$$S({\cal B}; {\cal P}(N), Q^{1/2})
\le 8 (1+\varepsilon) {C_N X\over \log N}
+ O(R_1 + R_2),
\leqno(5.5)$$
where
$$\eqalign{
& R_1 := \sum_{\scriptstyle q<Q\atop\scriptstyle q\mid P(Q^{1/2})}
3^{\nu_1(q)}
\Big|\sum_{\scriptstyle e\in {\cal E}'\atop\scriptstyle (e, q)=1}
\Big(\sum_{\scriptstyle p_2<p_3\le \kappa(d, e)
\atop\scriptstyle e p_3\equiv N ({\rm mod} \, q)} 1
- {1\over \varphi(q)} \sum_{p_2<p_3\le \kappa(d, e)} 1\Big)\Big|,
\cr
& R_2 := \sum_{\scriptstyle q<Q\atop\scriptstyle q\mid P(Q^{1/2})}
3^{\nu_1(q)}/\varphi(q)
\sum_{\scriptstyle e\in {\cal E}'\atop\scriptstyle (e, q)>1}
\sum_{p_2<p_3\le \kappa(d, e)} 1.
\cr}$$

We first estimate $R_2$.
Noticing that for $e\in {\cal E}'$ we have $e\le N^{1-\eta}$ and
the smallest prime factor of $e$ is $\ge \min\{p_1, V_k/\Delta\}\ge 
W_k^{1/s'}$,
we can deduce
$$\eqalign{
R_2
& \le N \sum_{q\le Q} {3^{\nu_1(q)}\over \varphi(q)}
\sum_{\scriptstyle e\le N^{1-\eta}\atop\scriptstyle (e, q)\ge 
W_k^{1/s'}} {1\over e}
\cr
& \ll N {\cal L} \sum_{m\ge W_k^{1/s'}} {1\over m}
\sum_{\scriptstyle q\le Q\atop{\atop\scriptstyle q \equiv 0 ({\rm mod} \, m)}}
{3^{\nu_1(q)}\over \varphi(q)}
\cr
& \ll N W_k^{1/3s'} \sum_{m\ge W_k^{1/s'}} {1\over m}
\sum_{\scriptstyle q\le Q\atop{\atop\scriptstyle q \equiv 0 ({\rm 
mod} \, m)}} {1\over q}
\cr
& \ll N W_k^{1/2s'} \sum_{m\ge W_k^{1/s'}} {1\over m^2}
\cr
& \ll_{\delta, k} {\Theta(N, \sigma)\over {\cal L}^3}.
\cr}
\leqno(5.6)$$

Next we estimate $R_1$.
Let $g(a) := \sum_{e\in {\cal E}', e = a} 1$.
Obviously for each $e = d n p_1 p_2\in {\cal E}'$,
the integers $d, n, p_1, p_2$ are pairwisely coprime.
Therefore they are uniquely determined by $e$.
Thus $g(a)\le 1$ and there are some injections $r_0(e) = d$ and $r(e) = p_2$.
Then we have
$$R_1
\ll
\!\!\!\!\! \sum_{\scriptstyle q\le Q\atop{\atop\scriptstyle (q, 
N)=1}} \!\!\!\!\!
\mu(q)^2 3^{\nu_1(q)} \,
\Big|
\!\!\! \sum_{\scriptstyle a\in I(a)\atop{\atop\scriptstyle (a, q)=1}} \!\!\!
g(a)
\Big(
\!\!\! \sum_{\scriptstyle p_2<p_3\le \kappa(d, a)
\atop{\atop\scriptstyle a p_3\equiv N ({\rm mod} \, q)}} \!\!\! 1
- {1\over \varphi(q)} \sum_{p_2<p_3\le \kappa(d, a)} 1\Big)
\Big|,$$
where $I(a) := \big(r_0(a) \underline{r_0(a)}^{2/s'}, \, 
N/\underline{r_0(a)}^{1/s'}\big)$.

Since $\underline d^{1/s'}<r(e)<\underline d^{1/s}$ and $e r(e)\le N$,
we can write
$$R_1
\ll R_1^{(1)} + R_1^{(2)} + R_1^{(3)},$$
where
$$\eqalign{
& R_1^{(1)} := \!\!\!
\sum_{\scriptstyle q\le Q\atop{\atop\scriptstyle (q, N)=1}} \!\!\!
\mu(q)^2 3^{\nu_1(q)}
\Big|
\sum_{\scriptstyle a\in I_1(a)\atop{\atop\scriptstyle (a, q)=1}}
g(a)
\Big(\sum_{\scriptstyle p_3\le \underline{r_0(a)}^{1/s}
\atop{\atop\scriptstyle a p_3\equiv N ({\rm mod} \, q)}} 1
- {1\over \varphi(q)} \sum_{p_3\le \underline{r_0(a)}^{1/s}} 1\Big)
\Big|,
\cr
& R_1^{(2)} := \!\!\!
\sum_{\scriptstyle q\le Q\atop{\atop\scriptstyle (q, N)=1}} \!\!\!
\mu(q)^2 3^{\nu_1(q)}
\Big|
\sum_{\scriptstyle a\in I_2(a)\atop{\atop\scriptstyle (a, q)=1}}
g(a)
\Big(\sum_{\scriptstyle p_3\le N/a
\atop{\atop\scriptstyle a p_3\equiv N ({\rm mod} \, q)}} 1
- {1\over \varphi(q)} \sum_{p_3\le N/a} 1\Big)
\Big|,
\cr
& R_1^{(3)} := \!\!\!
\sum_{\scriptstyle q\le Q\atop{\atop\scriptstyle (q, N)=1}} \!\!\!
\mu(q)^2 3^{\nu_1(q)}
\Big|
\sum_{\scriptstyle a\in I(a)\atop{\atop\scriptstyle (a, q)=1}}
g(a)
\Big(\sum_{\scriptstyle p_3\le r(a)
\atop{\atop\scriptstyle a p_3\equiv N ({\rm mod} \, q)}} 1
- {1\over \varphi(q)} \sum_{p_3\le r(a)} 1\Big)
\Big|,
\cr}$$
and $I_1(a) := \big(r_0(a) \underline{r_0(a)}^{2/s'}, \, 
N/\underline{r_0(a)}^{1/s}\big]$,
$I_2(a) := \big(N/\underline{r_0(a)}^{1/s}, \, 
N/\underline{r_0(a)}^{1/s'}\big)$.

\smallskip

Applying Lemma 2.3 yields
$R_1^{(j)}\ll_{\delta, k} N/{\cal L}^{5k+5}$ for $j = 1, 2, 3$.
Hence
$$R_1\ll \Theta(N, \sigma)/{\cal L}^3.
\leqno(5.7)$$

Replacing $(n, \, N P(p_2)/p_1)=1$ by $(n, \, N P(p_2))=1$ in the 
definition of $X$,
the difference is
$\ll_{\delta, k} N {\cal L}^2/d \underline d^{1/s'}\ll_{\delta, k} 
N/{\cal L}^2 d$.
Thus we can obtain, by Lemma 2.10, that
$$\eqalign{X
& = \sum_d \sigma(d)
\Big\{\mathop{\sum \,\, \sum \,\, \sum}_{\scriptstyle
\underline d^{1/s'}\le p_1<p_2<p_3<\underline d^{1/s}
\atop{\atop\scriptstyle (p_1 p_2, dN)=1}}
\sum_{\scriptstyle n\le N/d p_1 p_2 p_3\atop{\atop\scriptstyle (n, \, 
N P(p_2))=1}} 1
+ O_{\delta, k} (N/{\cal L}^2 d)
\Big\}
\cr
& \le (1 + \varepsilon) \sum_d \sigma(d)
\bigg\{
\mathop{\sum \,\, \sum \,\, \sum}_{\underline d^{1/s'}\le 
p_1<p_2<p_3<\underline d^{1/s}}
\!\!\!\! {N \omega(\log (N/d p_1 p_2 p_3)/\log p_2)\over d p_1 p_2 
p_3 \log p_2}
+ O_{\delta, k}\bigg({N\over {\cal L}^2 d}\bigg)
\bigg\}.
\cr}$$
By applying the prime number theorem, we can deduce
$$\eqalign{X
& \le (1 + \varepsilon) N \sum_d {\sigma(d)\over d \log \underline d}
\mathop{\int\int\int}_{1/s'\le t\le u\le v\le 1/s} 
\omega\bigg({\phi_{d, N} - t - u - v\over u}\bigg)
{\d t \d u \d v\over t u^2 v },
\cr}$$
where $\phi_{d, N} := \log(N/d)/\log \underline d$.
Obviously $\sigma(d)\not=0$ implies $\phi_{d, N}\ge 2$.
Thus
$$X
\le (1 + \varepsilon) I_1(s) N \sum_d {\sigma(d)\over d \log \underline d}
\le (1 + \varepsilon) I_1(s) N \sum_d {\sigma(d)\over \varphi(d) \log 
\underline d}.
\leqno(5.8)$$
Combining (5.4)--(5.8) and noticing $C_N\le C_{dN}$,
we obtain, for any $\varepsilon>0$ and $N\ge N_0(\varepsilon, \delta, k)$,
$$\sum_d \sigma(d) \Omega_3
\le \{2 I_1(s) + \varepsilon\} \Theta(N, \sigma).
\leqno(5.9)$$
Inserting (5.1), (5.2) and (5.9) into the inequality of Lemma 4.1 and 
noticing that
$$A(s') = A(s) + \int_2^{s'-1} {\log(t-1)\over t} \d t,
\qquad
a(s't) = \log(s't-1),
\leqno(5.10)$$
we find that
$$\Phi(N, \sigma, s)
\le \bigg\{
A(s)
- \Psi_1(s)
- {1\over 2} \int_{1-1/s}^{1-1/s'} {h_{k+1, N_0}(s't)\over t(1-t)} \d t
- H_{k, N_0}(s')
+ \varepsilon\bigg\}
\Theta(N, \sigma).$$
By the definition of $H_{k, N_0}(s)$, we must have
$$H_{k, N_0}(s)
\ge \Psi_1(s)
+ {1\over 2} \int_{1-1/s}^{1-1/s'} {h_{k+1, N_0}(s't)\over t(1-t)} \d t
+ H_{k, N_0}(s')
- \varepsilon.$$
Making $N_0\to \infty$ and then $\varepsilon\to 0$ yields
$$H_k(s)
\ge \Psi_1(s)
+ {1\over 2} \int_{1-1/s}^{1-1/s'} {h_{k+1}(s't)\over t(1-t)} \d t
+ H_k(s').$$
Now it remains to take $k\to \infty$ to get the desired result.
\hfill
$\square$

\proclaim Lemma 5.2.
For $5\ge s'\ge 3\ge s\ge 2$, $s'-s'/s\ge 2$ and $s\le 
\kappa_3<\kappa_2<\kappa_1\le s'$,
we have
$$\eqalign{H(s)
& \ge \Psi_2(s) + {4\over 5} H(s') + {1\over 5} H(\kappa_1)
+ {1\over 5} \int_{1-1/s}^{1-1/s'} {h(s't)\over t(1-t)} \d t
\cr
& \quad
+ {1\over 5} \int_{1-1/\kappa_2}^{1-1/s'} {h(s't)\over t(1-t)} \d t
+ {1\over 5} \int_{1-1/\kappa_3}^{1-1/s'} {h(s't)\over t(1-t)} \d t
\cr
& \quad
+ {1\over 5} \int_{1/s'}^{1/\kappa_2} {\d t\over t}
\int_t^{1/\kappa_2} {H(s'-s't-s'u)\over u(1-t-u)} \d u
\cr
& \quad
+ {1\over 5} \int_{1/s'}^{1/\kappa_1} {\d t\over t}
\int_{1/\kappa_2}^{1/\kappa_3} {H(s'-s't-s'u)\over u(1-t-u)} \d u
\cr
& \quad
+ {1\over 5} \int_{1/s'}^{1/\kappa_1} {\d t\over t}
\int_t^{1/\kappa_2} {H((1-t-u)/t)\over u(1-t-u)} \d u,
\cr}$$
where $\Psi_2(s) = \Psi_2(s, s', \kappa_1, \kappa_2, \kappa_3)$ is given by
$$\eqalign{\Psi_2(s)
& := - {2\over 5} \int_2^{s'-1} {\log(t-1)\over t} \d t
- {2\over 5} \int_2^{\kappa_1-1} {\log(t-1)\over t} \d t
- {1\over 5} \int_2^{\kappa_2-1} {\log(t-1)\over t} \d t
\cr
& \quad
+ {1\over 5} \int_{1-1/s}^{1-1/s'} {\log(s't-1)\over t(1-t)} \d t
+ {1\over 5} \int_{1-1/\kappa_3}^{1-1/\kappa_1} 
{\log(\kappa_1t-1)\over t(1-t)} \d t
- {2\over 5} \sum_{i=9}^{21} I_{2, i}(s)
\cr}$$
and $I_{2, i}(s) = I_{2, i}(s, s', \kappa_1, \kappa_2, \kappa_3)$ is given by
$$\eqalign{
I_{2, i}(s)
& :=  \max_{\phi\ge 2} \int_{\D_{2, i}}
\omega\bigg({\phi - t - u - v\over u}\bigg)
{\d t \d u \d v\over t u^2 v}
\qquad
(i = 9, \dots, 15),
\cr
I_{2, i}(s)
& := \max_{\phi\ge 2} \int_{\D_{2, i}}
\omega\bigg({\phi - t - u - v - w\over v}\bigg)
{\d t \d u \d v \d w\over t u v^2 w}
\qquad
(i = 16, \dots, 19),
\cr
I_{2, 20}(s)
& :=  \max_{\phi\ge 2} \int_{\D_{2, 20}}
\omega\bigg({\phi - t - u - v - w - x\over w}\bigg)
{\d t \d u \d v \d w \d x\over t u v w^2 x},
\cr
I_{2, 21}(s)
& := \max_{\phi\ge 2} \int_{\D_{2, 21}}
\omega\bigg({\phi - t - u - v - w - x - y\over x}\bigg)
{\d t \d u \d v \d w \d x \d y\over t u v w x^2 y}.
\cr}$$
The sets $\D_{2,i}$ $(9\le i\le 21)$ are defined as follows:
$$\eqalign{
\D_{2, 9}
& := \{(t, u, v) : 1/\kappa_1\le t\le u\le v\le 1/\kappa_3\},
\cr
\D_{2, 10}
& := \{(t, u, v) : 1/\kappa_1\le t\le u\le 1/\kappa_2\le v\le 1/s\},
\cr
\D_{2, 11}
& := \{(t, u, v) : 1/\kappa_1\le t\le 1/\kappa_2\le u\le v\le 1/\kappa_3\},
\cr
\D_{2, 12}
& := \{(t, u, v) : 1/s'\le t\le u\le 1/\kappa_1, \, 1/\kappa_3\le v\le 1/s\},
\cr
\D_{2, 13}
& := \{(t, u, v) : 1/s'\le t\le 1/\kappa_1\le u\le 1/\kappa_2\le v\le 1/s\},
\cr
\D_{2, 14}
& := \{(t, u, v) : 1/s'\le t\le 1/\kappa_1, \, 1/\kappa_2\le u\le v\le 1/s\},
\cr
\D_{2, 15}
& := \{(t, u, v) : 1/\kappa_1\le t\le 1/\kappa_2\le u\le 
1/\kappa_3\le v\le 1/s\},
\cr
\D_{2, 16}
& := \{(t, u, v, w) : 1/\kappa_2\le t\le u\le v\le w\le 1/\kappa_3\},
\cr
\D_{2, 17}
& := \{(t, u, v, w) : 1/\kappa_2\le t\le u\le v\le 1/\kappa_3\le w\le 1/s\},
\cr
\D_{2, 18}
& := \{(t, u, v, w) : 1/\kappa_2\le t\le u\le 1/\kappa_3\le v\le w\le 1/s\},
\cr
\D_{2, 19}
& := \{(t, u, v, w) : 1/\kappa_1\le t\le 1/\kappa_2,
1/\kappa_3\le u\le v\le w\le 1/s\},
\cr
\D_{2, 20}
& := \{(t, u, v, w, x) : 1/\kappa_2\le t\le 1/\kappa_3\le u\le v\le 
w\le x\le 1/s\},
\cr
\D_{2, 21}
& := \{(t, u, v, w, x, y) : 1/\kappa_3\le t\le u\le v\le w\le x\le y\le 1/s\}.
\cr}$$

\noindent{\sl Proof}.
By (3.3) of Lemma 3.2, we have
$$\sum_d \sigma(d) \Gamma_1
\le \big\{4 A(s') + A(\kappa_1)
- 4 H_{k, N_0}(s') - H_{k, N_0}(\kappa_1)\big\} \Theta(N, \sigma).
\leqno(5.11)$$
Similar to (3.20), we can prove
$$\leqalignno{
\sum_d \sigma(d) \Gamma_2
& \ge \bigg\{\int_{1-1/s}^{1-1/s'} {a(s't) + h_{k+1, N_0}(s't)\over 
t(1-t)} \d t
- \varepsilon\bigg\} \Theta(N, \sigma),
& (5.12)
\cr
\sum_d \sigma(d) \Gamma_3
& \ge \bigg\{\int_{1-1/\kappa_2}^{1-1/s'} {a(s't) + h_{k+1, 
N_0}(s't)\over t(1-t)} \d t
- \varepsilon\bigg\} \Theta(N, \sigma),
& (5.13)
\cr
\sum_d \sigma(d) \Gamma_4
& \ge \bigg\{\int_{1-1/\kappa_3}^{1-1/s'} {a(s't) + h_{k+1, 
N_0}(s't)\over t(1-t)} \d t
- \varepsilon\bigg\} \Theta(N, \sigma).
& (5.14)
\cr}$$
Similar to (3.20) and
in view of $A(s'-s't-s'u) = A((1-t-u)/t) = 1$,
we have
$$\leqalignno{\qquad
\sum_d \sigma(d) \Gamma_5
& \le \bigg\{
\int_{1/s'}^{1/\kappa_2} {\d t\over t} \int_t^{1/\kappa_2}
{1 - H_{k+2, N_0}(s'-s't-s'u)\over u(1-t-u)} \d u
+ \varepsilon
\bigg\} \Theta(N, \sigma),
& (5.15)
\cr
\sum_d \sigma(d) \Gamma_6
& \le \bigg\{
\int_{1/s'}^{1/\kappa_1} {\d t\over t} \int_{1/\kappa_2}^{1/\kappa_3}
{1 - H_{k+2, N_0}(s'-s't-s'u)\over u(1-t-u)} \d u
+ \varepsilon
\bigg\} \Theta(N, \sigma),
& (5.16)
\cr
\sum_d \sigma(d) \Gamma_7
& \le \bigg\{
\int_{1/s'}^{1/\kappa_1} {\d t\over t} \int_t^{1/\kappa_1}
{1 - H_{k+2, N_0}((1-t-u)/t)\over u(1-t-u)} \d u
+ \varepsilon\bigg\} \Theta(N, \sigma),
& (5.17)
\cr
\cr
\sum_d \sigma(d) \Gamma_8
& \le \bigg\{
\int_{1/s'}^{1/\kappa_1} {\d t\over t} \int_{1/\kappa_1}^{1/\kappa_2}
{1 - H_{k+2, N_0}((1-t-u)/t)\over u(1-t-u)} \d u
+ \varepsilon\bigg\} \Theta(N, \sigma). 
& (5.18)
\cr}$$

We have also, for $i = 9, \dots, 21$,
$$\sum_d \sigma(d) \Gamma_i
\le \{2I_{2, i}(s) + \varepsilon\} \Theta(N, \sigma).
\leqno(5.19)$$
As before, inserting (5.11)--(5.19) into the inequality of Lemma 4.2 and
using the definition of $H_{k, N_0}(s)$, we can deduce
$$5 H_{k, N_0}(s)
\ge A(s, s') + B(s, s') - \varepsilon,
\leqno(5.20)$$
where
$$\eqalign{A(s, s')
& := 5 A(s) - 4 A(s') - A(\kappa_1) + 4 H_{k, N_0}(s') + H_{k, N_0}(\kappa_1)
\cr\noalign{\medskip}
& \quad
+ \int_{1-1/s}^{1-1/s'}  {a(s't)\over t(1-t)} \d t
+ \int_{1-1/\kappa_2}^{1-1/s'} {a(s't)\over t(1-t)} \d t
+ \int_{1-1/\kappa_3}^{1-1/s'} {a(s't)\over t(1-t)} \d t
\cr
& \quad
- \int_{1/s'}^{1/\kappa_2} \int_t^{1/\kappa_2} {\d t \d u\over t u(1-t-u)}
- \int_{1/s'}^{1/\kappa_1} \int_t^{1/\kappa_3} {\d t \d u\over t u(1-t-u)}
- 2 \sum_{i=9}^{21} I_{2, i}(s)
\cr}$$
and
$$\eqalign{B(s, s')
& := \bigg(\int_{1-1/s}^{1-1/s'}
+ \int_{1-1/\kappa_2}^{1-1/s'}
+ \int_{1-1/\kappa_3}^{1-1/s'} \bigg) {h_{k+1, N_0}(s't)\over t(1-t)} \d t
\cr
& \quad
+ \int_{1/s'}^{1/\kappa_2} {\d t\over t}
\int_t^{1/\kappa_2} {H_{k+2, N_0}(s'-s't-s'u)\over u(1-t-u)} \d u
\cr
& \quad
+ \int_{1/s'}^{1/\kappa_1} {\d t\over t}
\int_{1/\kappa_2}^{1/\kappa_3} {H_{k+2, N_0}(s'-s't-s'u)\over u(1-t-u)} \d u
\cr
& \quad
+ \int_{1/s'}^{1/\kappa_1} {\d t\over t}
\int_t^{1/\kappa_2} {H_{k+2, N_0}((1-t-u)/t)\over u(1-t-u)} \d u.
\cr}$$

For $a\ge b>2$, we have
$$\eqalign{
\int_{1/a}^{1/b} {\d t\over t} \int_t^{1/b} {\d u\over u(1-t-u)}
& = \int_{1/a}^{1/b} {\d u\over u} \int_{1/a}^u {\d t\over t (1-t-u)}
\cr
& = \int_{1/a}^{1/b} {\log(a-1-au)- \log(1/u-2)\over u (1-u)} \d u
\cr
& = \int_{1-1/b}^{1-1/a} {\log(at-1)\over t (1-t)} \d t
- \int_{b-1}^{a-1} {\log(t-1)\over t} \d t,
\cr}$$
where we have used the change of variables $t = 1 - u$ and $t = 1/u - 1$ respectively.

Similarly for $a\ge b\ge c\ge d>2$, we have
$$\int_{1/a}^{1/b} {\d t\over t} \int_{1/c}^{1/d} {\d u\over u(1-t-u)}
= \int_{1-1/d}^{1-1/c} {\log(at-1)\over t(1-t)} \d t
- \int_{1-1/d}^{1-1/c} {\log(bt-1)\over t(1-t)} \d t.$$

By using these two relations and (5.10),
a simple calculation shows
$$A(s, s')
= 5 \Psi_2(s) + 4 H_{k, N_0}(s') + H_{k, N_0}(\kappa_1).$$
Inserting this into (5.20) and making $N\to\infty$, $\varepsilon\to 
0$ and $k\to\infty$,
we obtain the desired inequality.
This completes the proof.
\hfill
$\square$

\vskip 5mm

\noindent{\bf \S\ 6. Proofs of Propositions 3 and 4}

\medskip

We first prove a preliminary lemma.
Let ${\bf 1}_{[a, b]}(t)$ be the characteristic function of the 
interval $[a, b]$.
We put
$$\qquad
\sigma(a, b, c)
:= \int_a^b \log{c\over t-1} {\d t\over t},
\qquad
\sigma_0(t) := {\sigma(3, t+2, t+1)\over 1-\sigma(3, 5, 4)}.$$

\proclaim Lemma 6.1.
Let $3\le s'\le 5$, $0<a<b<1$ and $2\le ac<bc\le 4$.
Then we have
$$\leqalignno{\qquad\qquad
& h(4)
\ge \int_1^3 H(t) {\sigma_0(t)\over t} \d t.
& (6.1)
\cr
& H(s')
\ge \int_1^3 H(t)
\bigg\{{\sigma_0(t)\over t}\log\bigg({4\over s'-1}\bigg)
+ {{\bf 1}_{[s'-2,3]}(t)\over t}\log\bigg({t+1\over s'-1}\bigg)\bigg\} \d t,
& (6.2)
\cr
& \int_a^b {h(c t)\over t(1-t)} \d t
\ge \log\bigg({b-ab\over a-ab}\bigg) \int_1^3 H(t) 
{\sigma_0(t)+{\bf 1}_{[bc-1, 3]}(t)\over t} \d t
& (6.3)
\cr
& \qquad\qquad\qquad\qquad
+ \int_1^3 H(t) {{\bf 1}_{[ac-1, bc-1]}(t)\over 
t}\log\bigg({(1-a)(t+1)\over a(c-1-t)}\bigg) \d t.
\cr}$$

\noindent{\sl Proof}.
By Proposition 2, we have
$$\leqalignno{
H(s')
& \ge \int_{s'-1}^4 {h(u)\over u} \d u
\ge \int_{s'-1}^4 \bigg(h(4) + \int_{u-1}^3 H(t) {\d t\over t}\bigg) 
{\d u\over u}
& (6.4)
\cr
& = h(4) \log\bigg({4\over s'-1}\bigg)
+ \int_1^3 H(t) {{\bf 1}_{[s'-2,3]}(t)\over t} \log\bigg({t+1\over 
s'-1}\bigg) \d t.
\cr}$$
From Proposition 2 and (6.4), we deduce
$$h(4)
\ge \int_3^5 {H(v)\over v} \d v
\ge h(4) \sigma(3, 5, 4)
+ \int_1^3 H(t) {\sigma(3, t+2, t+1)\over t} \d t,$$
which implies the inequality (6.1).

The inequality (6.2) follows immediately from (6.4) and (6.1).

By using Proposition 2, we have
$$\eqalign{\int_a^b {h(c t)\over t(1-t)} \d t
& = c \int_{ac}^{bc} { h(u)\over u(c-u)} \d u
\ge c \int_{ac}^{bc} {\d u\over u(c-u)}
\bigg(h(4) + \int_{u-1}^3 H(t) {\d t\over t}\bigg)
\cr
& = \bigg\{h(4) + \int_1^3 H(t) {{\bf 1}_{[bc-1, 3]}(t)\over t} \d t\bigg\}
\log\bigg({b-ab\over a-ab}\bigg)
\cr
& \quad
+ \int_1^3 H(t) {{\bf 1}_{[ac-1, bc-1]}(t)\over t} 
\log\bigg({(1-a)(t+1)\over a(c-1-t)}\bigg) \d t,
\cr}$$
which combines (6.1) to give (6.3).
This completes the proof.
\hfill
$\square$

\goodbreak
\medskip

{\sl Proof of Proposition 3}.
By using Lemma 6.1, a simple calculation shows
$${1\over 2} \int_{1-1/s}^{1-1/s'} {h(s' t)\over t(1-t)} \d t + H(s')
\ge \int_1^3 H(t) \Xi_1(t; s) \d t,$$
which, together with Lemma 4.1, implies the desired result.
\hfill
$\square$

\smallskip

{\sl Proof of Proposition 4}.
From (6.1)--(6.3), we can deduce
$$\leqalignno{\qquad
\int_{1-1/s}^{1-1/s'} 
& {h(s't)\over t(1-t)} \d t
+ \int_{1-1/\kappa_2}^{1-1/s'} {h(s't)\over t(1-t)} \d t
+ \int_{1-1/\kappa_3}^{1-1/s'} {h(s't)\over t(1-t)} \d t
+ 4 H(s') + H(\kappa_1)
& (6.5)
\cr
& \ge \int_1^3 H(t) {{\bf 1}_{[\alpha_2, 3]}(t)\over t}
\log\bigg({(t+1)^5\over 
(s-1)(s'-1)(\kappa_1-1)(\kappa_2-1)(\kappa_3-1)}\bigg) \d t
\cr
& \quad
+ \int_1^3 H(t) {\sigma_0(t)\over t}
\log\bigg({1024\over 
(s-1)(s'-1)(\kappa_1-1)(\kappa_2-1)(\kappa_3-1)}\bigg) \d t
\cr
& \quad
+ \int_1^3 H(t)
{{\bf 1}_{[\alpha_5, \alpha_2]}(t)\over t} \log\bigg({t+1\over 
(\kappa_3-1)(s'-1-t)}\bigg) \d t
\cr
& \quad
+ \int_1^3 H(t)
{{\bf 1}_{[\alpha_4, \alpha_2]}(t)\over t} \log\bigg({t+1\over 
(\kappa_2-1)(s'-1-t)}\bigg) \d t
\cr
& \quad
+ \int_1^3 H(t)
{{\bf 1}_{[\alpha_3, \alpha_2]}(t)\over t} \log\bigg({t+1\over 
(s-1)(s'-1-t)}\bigg) \d t
\cr
& \quad
+ \int_1^3 H(t) {{\bf 1}_{[\alpha_1, \alpha_2]}(t)\over t}
\log\bigg({t+1\over \kappa_1-1}\bigg) \d t.
\cr}$$

By the change of variable $v = s'(1 - t - u)$, we have
$$\int_{1/s'}^{1/\kappa_2} {\d t\over t} \int_t^{1/\kappa_2} 
{H(s'-s't-s'u)\over u (1-t-u)} \d u
= \int_{1/s'}^{1/\kappa_2} {\d t\over t}
\int_{s'(1-1/\kappa_2-t)}^{s'(1-2t)} {s' H(v)\over v (s'-s't-v)} \d v.$$
Interchanging the order of integration and a simple calculation show that
$$\leqalignno{\qquad
\int_{1/s'}^{1/\kappa_2} {\d t\over t} \int_t^{1/\kappa_2} 
{H(s'-s't-s'u)\over u (1-t-u)} \d u
& = \int_1^3 H(t)
\bigg\{
{{\bf 1}_{[\alpha_6, \alpha_4]}(t)\over t(1-t/s')}
\log\bigg({s'\over \kappa_2s'-s'-\kappa_2t}\bigg)
& (6.6)
\cr
& \qquad\qquad\qquad
+ {{\bf 1}_{[\alpha_4, \alpha_2]}(t)\over t(1-t/s')} \log(s'-1-t)\bigg\} \d t.
\cr}$$

Similarly we can prove
$$\int_{1/s'}^{1/\kappa_1} {\d t\over t}
\int_{1/\kappa_2}^{1/\kappa_3} {H(s'-s't-s'u)\over u (1-t-u)} \d u
= \int_1^3 H(t) L_1(t) \d t
\leqno(6.7)$$
where
$$\eqalign{L_1(t)
& := {{\bf 1}_{[\alpha_7, \alpha_5]}(t)\over t(1-t/s')}
\log\bigg({s'^2\over (\kappa_1s'-s'-\kappa_1t)(\kappa_3s'-s'-\kappa_3t)}\bigg)
\cr
& \quad
+ {{\bf 1}_{[\alpha_5, \alpha_8]}(t)\over t(1-t/s')}
\log\bigg({s'(s'-1-t)\over \kappa_1s'-s'-\kappa_1t}\bigg)
\cr
& \quad
+ {{\bf 1}_{[\alpha_8, \alpha_4]}(t)\over t(1-t/s')}
\log\bigg({(s'-1-t)(\kappa_2s'-s'-\kappa_2t)\over s'}\bigg),
\cr}$$
and
$$\int_{1/s'}^{1/\kappa_1} {\d t\over t} \int_t^{1/\kappa_2} 
{H((1-t-u)/t)\over u (1-t-u)} \d u
= \int_1^3 H(t) L_2(t) \d t
\leqno(6.8)$$
with
$$\eqalign{L_2(t)
& := {{\bf 1}_{[\alpha_9, \alpha_1]}(t)\over t}
\log\bigg({t+1\over (\kappa_2-1)(\kappa_1-1-t)}\bigg)
\cr
& \quad
+ {{\bf 1}_{[\alpha_1, \alpha_4]}(t)\over t}
\log\bigg({t+1\over \kappa_2-1}\bigg)
+ {{\bf 1}_{[\alpha_4, \alpha_2]}(t)\over t}
\log (s'-1-t).
\cr}$$

Now by inserting (6.5)--(6.8) into Lemma 5.2,
we easily deduce the required result.
\hfill $\square$

\vskip 5mm

\noindent{\bf \S\ 7. Proof of Theorem 1}

\medskip

We need to resolve the functional inequalities (3.21) and (3.22).
It seems very difficult to give the exact solutions,
because we only know that $H(s)$ is decreasing.
Next we shall give a numeric lower bound for solution by using discretion,
which is sufficient to prove Theorem 1.

Put $s_0:= 1$ and $s_i := 2 + 0.1\times (i+1)$ for $i=1, \dots, 9$.
Since $H(s)$ is decreasing on $[1, 10]$, Proposition 4 allows us to deduce
$$H(s_i)\ge \Psi_2(s_i) + \sum_{j=1}^9 a_{i, j} H(s_j),
\leqno(7.1)$$
where
$$a_{i, j} := \int_{s_{j-1}}^{s_j} \Xi_2(t, s_i) \d t
\qquad
(i=1, \dots, 4; \, \, j=1, \dots, 9).$$
Similarly Proposition 3 implies
$$H(s_i)\ge \Psi_1(s_i) + \sum_{j=1}^9 a_{i, j} H(s_j),
\leqno(7.2)$$
where
$$a_{i, j} := \int_{s_{j-1}}^{s_j} \Xi_1(t, s_i) \d t
\qquad
(i=5, \dots, 9; \, \, j=1, \dots, 9)$$

\bigskip

\centerline{{\bf Table 1.} Choice of parameters}
\vskip -4mm
$$\displaylines{
\vbox{\tabskip = 0pt\offinterlineskip
\halign{
\vrule # & &\hfil$ $ $#$ $ \!\! $ \hfil & \vrule #\cr
\noalign{\hrule}
height 2mm
&&&&&&&&&&&&&&&&
\cr
& \,\,     i              \,\,        &
& \,\,\,   s_i            \,\,\,      &
& \,\,\,\, s_i'           \,\,\,\,    &
& \,\,\,   \kappa_{1,i}   \,\,\,      &
& \,\,\,   \kappa_{2,i}   \,\,\,      &
& \,\,\,   \kappa_{3,i}   \,\,\,      &
& \quad\,\,\, \Psi_1(s_i) \quad\,\,\, &
& \quad\,\,\, \Psi_2(s_i) \quad\,\,\, &
\cr
height 2mm
&&&&&&&&&&&&&&&&
\cr
\noalign{\hrule}
height 2mm
&&&&&&&&&&&&&&&&
\cr
& 1           &
& 2.2         &
& 4.54        &
& 3.53        &
& 2.90        &
& 2.44        &
&             &
& 0.015826357 &
\cr
height 2mm
&&&&&&&&&&&&&&&&
\cr
\noalign{\hrule}
height 2mm
&&&&&&&&&&&&&&&&
\cr
& 2           &
& 2.3         &
& 4.50        &
& 3.54        &
& 2.88        &
& 2.43        &
&             &
& 0.015247971 &
\cr
height 2mm
&&&&&&&&&&&&&&&&
\cr
\noalign{\hrule}
height 2mm
&&&&&&&&&&&&&&&&
\cr
& 3           &
& 2.4         &
& 4.46        &
& 3.57        &
& 2.87        &
& 2.40        &
&             &
& 0.013898757 &
\cr
height 2mm
&&&&&&&&&&&&&&&&
\cr
\noalign{\hrule}
height 2mm
&&&&&&&&&&&&&&&&
\cr
& 4           &
& 2.5         &
& 4.12        &
& 3.56        &
& 2.91        &
& 2.50        &
&             &
& 0.011776059 &
\cr
height 2mm
&&&&&&&&&&&&&&&&
\cr
\noalign{\hrule}
height 2mm
&&&&&&&&&&&&&&&&
\cr
& 5           &
& 2.6         &
& 3.58        &
&             &
&             &
&             &
& 0.009405211 &
&             &
\cr
height 2mm
&&&&&&&&&&&&&&&&
\cr
\noalign{\hrule}
height 2mm
&&&&&&&&&&&&&&&&
\cr
& 6           &
& 2.7         &
& 3.47        &
&             &
&             &
&             &
& 0.006558950 &
&             &
\cr
height 2mm
&&&&&&&&&&&&&&&&
\cr
\noalign{\hrule}
height 2mm
&&&&&&&&&&&&&&&&
\cr
& 7           &
& 2.8         &
& 3.34        &
&             &
&             &
&             &
& 0.003536751 &
&             &
\cr
height 2mm
&&&&&&&&&&&&&&&&
\cr
\noalign{\hrule}
height 2mm
&&&&&&&&&&&&&&&&
\cr
& 8           &
& 2.9         &
& 3.19        &
&             &
&             &
&             &
& 0.001056651 &
&             &
\cr
height 2mm
&&&&&&&&&&&&&&&&
\cr
\noalign{\hrule}
height 2mm
&&&&&&&&&&&&&&&&
\cr
& 9           &
& 3.0         &
& 3.00        &
&             &
&             &
&             &
& 0.000000000 &
&             &
\cr
height 2mm
&&&&&&&&&&&&&&&&
\cr
\noalign {\hrule}}}
\cr}$$

\bigskip

The parameters $s_i'$, $\kappa_{1, i}$, $\kappa_{2, i}$ and $\kappa_{3, i}$ are chosen such that
$\Psi_1(s_i)$ or $\Psi_2(s_i)$ is maximal.

We put
$${\bf A} := \pmatrix{
a_{1, 1} & \cdots & a_{1, 9}
\cr
\vdots & & \vdots
\cr
a_{9, 1} & \cdots & a_{9, 9}
\cr},
\qquad
{\bf H}
:= \pmatrix{
H(s_1)
\cr
\vdots
\cr
H(s_9)
\cr},
\qquad
{\bf B}
:= \pmatrix{
\Psi_2(s_1)
\cr
\vdots
\cr
\Psi_2(s_4)
\cr\noalign{\smallskip}
\Psi_1(s_5)
\cr
\vdots
\cr
\Psi_1(s_9)
\cr}.$$
Then (7.1) and (7.2) can be written as
$$({\bf I} - {\bf A} ){\bf H}\ge {\bf B},
\leqno(7.3)$$
where ${\bf I}$ is the unit matrix.

In order to resolve (7.3), we first solve the system of linear equations
$$({\bf I} - {\bf A}) {\bf X} = {\bf B},
\leqno(7.4)$$
by using {\it Maple} and obtain
$${\bf X}
= \pmatrix{
0.0223939 \cdots
\cr\noalign{\vskip 0.5mm}
0.0217196 \cdots
\cr\noalign{\vskip 0.5mm}
0.0202876 \cdots
\cr\noalign{\vskip 0.5mm}
0.0181433 \cdots
\cr\noalign{\vskip 0.5mm}
0.0158644 \cdots
\cr\noalign{\vskip 0.5mm}
0.0129923 \cdots
\cr\noalign{\vskip 0.5mm}
0.0100686 \cdots
\cr\noalign{\vskip 0.5mm}
0.0078162 \cdots
\cr\noalign{\vskip 0.5mm}
0.0072943 \cdots
\cr}.$$
 
From (7.3) and (7.4), we deduce that
$$({\bf I} - {\bf A}) ({\bf H} - {\bf X}) \ge {\bf 0}.$$
Since all elements of $({\bf I} - {\bf A})^{-1}$ are positive, it follows that
$${\bf H}\ge {\bf X}.$$
In particular we have
$$H(2.2)\ge 0.0223939.$$
Now taking $\sigma = \{1\}$ and $s=2.2$ in (3.3) of Lemma 3.2,
we find, for $\delta$ sufficiently small, $N_0$ sufficiently large 
and $N\ge N_0$,
$$\eqalign{D(N)
& \le S\big({\cal A}; {\cal P}(N), N^{(1/2-\delta)/2.2}\big)
= \Phi(N, \{1\}, 2.2)
\cr
& \le \big\{A(2.2) - H_{k, N_0}(2.2)\big\} {4 C_N {\rm li}(N)\over 
\log(N^{1/2-\delta})}
\cr
& \le 8 (1 - 0.0223938) \Theta(N)
\cr\noalign{\smallskip}
& \le 7.82085 \Theta(N).
\cr}$$
This completes the proof of Theorem 1.

\smallskip

{\bf Remark 2.}
(i)
The constant $s_1=2.2$ comes from the fact that $\Psi_2(s)$ attains the maximal value at $s=s_1$ (approximately).
Since $H(s)$ is decreasing on $[1, 10]$, we have $H(2.1)\ge 0.0223939$.
In order to obtain a better lower (which leads to a smaller constant than 7.82085), 
we must look for a new weighted inequality (as in Lemma 4.1 and 4.2) 
such that the corresponding main term $\Psi(2.1)$ has a lager lower bound than 0.015826357.

(ii)
If we divide the interval $[2, 3]$ into more subintervals than 9,
it is certain that we can obtain a better result.
But the improvement is very minuscule.

\vskip 5mm

\noindent{\bf \S\ 8. Proof of Theorem 3}

\medskip

In the case of the twin primes problem, we need to sieve the following sequence
$${\cal B} := \{p+2 : p\le x\}.$$
Thinking to Lemmas 2.7 and 2.9,
we have ${4\over 7}$ for the level of distribution
in place ${1\over 2}$ in the Bombieri--Vinogradov theorem.
Thus we can take $Q := x^{4/7-\delta}$ and $\underline d := Q/d$
in the definitions described in Section 3.
As before, we can prove the corresponding Propositions 3 and 4 with the 
following modification:
In the definition of $\Psi_1(s)$ we add a factor ${7\over 8}$ before $I_1(s)$,
and in the definition of $\Psi_2(s)$
we replace the factor ${2\over 5}$ before the sum by ${7\over 20}$.
When we use the switching principle to treat the terms $\Omega_3$ and $\Gamma_i$ for $5\le i\le 21$,
the related error terms can be estimated by using Lemma 2.9 
which has ${4\over 7}$ for the level of distribution (see [26], page 380).

\bigskip

\centerline{{\bf Table 2.} Choice of parameters}
\vskip -4mm
$$\displaylines{
\vbox{\tabskip = 0pt\offinterlineskip
\halign{
\vrule # & &\hfil$ $ $#$ $ \!\! $ \hfil & \vrule #\cr
\noalign{\hrule}
height 2mm
&&&&&&&&&&&&&&&&
\cr
& \,\,     i              \,\,        &
& \,\,\,   s_i            \,\,\,      &
& \,\,\,\, s_i'           \,\,\,\,    &
& \,\,\,   \kappa_{1,i}   \,\,\,      &
& \,\,\,   \kappa_{2,i}   \,\,\,      &
& \,\,\,   \kappa_{3,i}   \,\,\,      &
& \quad\,\,\, \Psi_1(s_i) \quad\,\,\, &
& \quad\,\,\, \Psi_2(s_i) \quad\,\,\, &
\cr
height 2mm
&&&&&&&&&&&&&&&&
\cr
\noalign{\hrule}
height 2mm
&&&&&&&&&&&&&&&&
\cr
& 1           &
& 2.1         &
& 4.93        &
& 3.62        &
& 2.86        &
& 2.34        &
&             &
& 0.020914508 &
\cr
height 2mm
&&&&&&&&&&&&&&&&
\cr
\noalign{\hrule}
height 2mm
&&&&&&&&&&&&&&&&
\cr
& 2           &
& 2.2         &
& 4.91        &
& 3.62        &
& 2.85        &
& 2.33        &
&             &
& 0.020399717 &
\cr
height 2mm
&&&&&&&&&&&&&&&&
\cr
\noalign{\hrule}
height 2mm
&&&&&&&&&&&&&&&&
\cr
& 3           &
& 2.3         &
& 5.00        &
& 3.63        &
& 2.82        &
& 2.30        &
&             &
& 0.019005124 &
\cr
height 2mm
&&&&&&&&&&&&&&&&
\cr
\noalign{\hrule}
height 2mm
&&&&&&&&&&&&&&&&
\cr
& 4           &
& 2.4         &
& 4.52        &
& 3.64        &
& 2.87        &
& 2.40        &
&             &
& 0.016618139 &
\cr
height 2mm
&&&&&&&&&&&&&&&&
\cr
\noalign{\hrule}
height 2mm
&&&&&&&&&&&&&&&&
\cr
& 5           &
& 2.5         &
& 3.72        &
&             &
&             &
&             &
& 0.013597508 &
&             &
\cr
height 2mm
&&&&&&&&&&&&&&&&
\cr
\noalign{\hrule}
height 2mm
&&&&&&&&&&&&&&&&
\cr
& 6           &
& 2.6         &
& 3.62        &
&             &
&             &
&             &
& 0.010644985 &
&             &
\cr
height 2mm
&&&&&&&&&&&&&&&&
\cr
\noalign{\hrule}
height 2mm
&&&&&&&&&&&&&&&&
\cr
& 7           &
& 2.7         &
& 3.49        &
&             &
&             &
&             &
& 0.007155027 &
&             &
\cr
height 2mm
&&&&&&&&&&&&&&&&
\cr
\noalign{\hrule}
height 2mm
&&&&&&&&&&&&&&&&
\cr
& 8           &
& 2.8         &
& 3.35        &
&             &
&             &
&             &
& 0.003741586 &
&             &
\cr
height 2mm
&&&&&&&&&&&&&&&&
\cr
\noalign{\hrule}
height 2mm
&&&&&&&&&&&&&&&&
\cr
& 9           &
& 2.9         &
& 3.19        &
&             &
&             &
&             &
& 0.001087780 &
&             &
\cr
height 2mm
&&&&&&&&&&&&&&&&
\cr
\noalign{\hrule}
height 2mm
&&&&&&&&&&&&&&&&
\cr
& 10          &
& 3.0         &
& 3.00        &
&             &
&             &
&             &
& 0.000000000 &
&             &
\cr
height 2mm
&&&&&&&&&&&&&&&&
\cr
\noalign {\hrule}}}
\cr}$$

\bigskip

As before we can prove
$$H(2.1)\ge 0.0287118.$$
Thus for $\delta$ sufficiently small, $x_0$ sufficiently large and 
$x\ge x_0$, we have
$$\eqalign{\pi_2(x)
& \le S\big({\cal B}; {\cal P}(2), x^{(1/2-\delta)/2.1}\big)
\cr
& \le 3.5 (1 - 0.0287117) \, \Pi(x)
\cr
& \le 3.39951 \, \Pi(x).
\cr}$$
This completes the proof of Theorem 3.

\goodbreak

\vskip 5mm

\noindent{\bf \S\ 9. Chen's system of weights}

\medskip

Let $${\cal A} := \{N-p : p\le N\},
\qquad
{\cal P}(q) := \{p : (p, q) = 1\}.$$
The inequality (9.1) below appeared in [9] (page 479, (11))
with $(\kappa, \sigma) = ({1\over 12}, {1\over 3.047}), ({1\over 
9.2}, {1\over 3.41})$
without proof. 
Cai [3] gave a proof with an extra assumption 
$3\sigma+\kappa>1$.
Here we present a proof
without Cai's assumption.
This removal is important in our argument.

\proclaim Lemma 9.1.
Let $0<\kappa<\sigma<{1\over 3}$. Then we have
$$D_{1,2}(N)\ge
S({\cal A}; {\cal P}(N), N^\kappa)
- \dm S_1
- S_2
- \dm S_3
+ \dm S_4
+ O(N^{1-\kappa}),
\leqno(9.1)$$
where $S_i = S_i(\kappa, \sigma)$ $(1\le i\le 4)$ are defined by
$$\eqalign{
& S_1 := \sum_{\scriptstyle N^\kappa\le p<N^\sigma\atop\scriptstyle (p, N)=1}
S({\cal A}_p; {\cal P}(N), N^\kappa),
\cr
& S_2 := \mathop{\sum \,\, \sum}_{\scriptstyle N^\sigma\le 
p_1<p_2<(N/p_1)^{1/2}
\atop\scriptstyle (p_1p_2, N)=1}
S({\cal A}_{p_1 p_2}; {\cal P}(Np_1), p_2),
\cr
& S_3 := \mathop{\sum \,\, \sum}_{\scriptstyle N^\kappa\le 
p_1<N^\sigma\le p_2<(N/p_1)^{1/2}
\atop\scriptstyle (p_1p_2, N)=1}
S({\cal A}_{p_1 p_2}; {\cal P}(Np_1), p_2),
\cr
& S_4 := \mathop{\sum \,\, \sum \,\, \sum}_{\scriptstyle N^\kappa\le 
p_1<p_2<p_3<N^\sigma
\atop\scriptstyle (p_1p_2p_3, N)=1}
S({\cal A}_{p_1 p_2 p_3}; {\cal P}(Np_1), p_2).
\cr}$$

\noindent{\sl Proof}.
Clearly the desired inequality (9.1) is equivalent to
$$D_{1,2}(N)\ge
\sum_{a\in {\cal A}, \, (a, P(N^\kappa))=1}
\big(1 - \dm s_1(a) - s_2(a) - \dm s_3(a) + \dm s_4(a)\big)
+ O(N^{1-\kappa}),
\leqno(9.2)$$
where
$$\eqalign{
& s_1(a) := \sum_{\scriptstyle N^\kappa\le 
p<N^\sigma\atop\scriptstyle p\mid a, \, (p, N)=1} 1,
\cr
& s_2(a)
:= \mathop{\sum\,\,\sum}_{{\scriptstyle N^\sigma\le p_1<p_2<(N/p_1)^{1/2}
\atop\scriptstyle p_1p_2\mid a, \, (p_1p_2, N)=1}
\atop\scriptstyle p\mid (a/p_1p_2)\Rightarrow p\ge p_2} 1,
\cr
& s_3(a)
:= \mathop{\sum\,\,\sum}_{{\scriptstyle N^\kappa\le p_1<N^\sigma\le 
p_2<(N/p_1)^{1/2}
\atop\scriptstyle p_1p_2\mid a, \, (p_1p_2, N)=1}
\atop\scriptstyle p\mid (a/p_1p_2)\Rightarrow p\ge p_2} 1,
\cr
& s_4(a)
:= \mathop{\sum\,\,\sum\,\,\sum}_{{\scriptstyle N^\kappa\le 
p_1<p_2<p_3<N^\sigma
\atop\scriptstyle p_1p_2p_3\mid a, \, (p_1p_2p_3, N)=1}
\atop\scriptstyle p\mid (a/p_1p_2p_3)\Rightarrow p\ge p_2} 1.
\cr}$$

Let 
$$\delta^*(a) := \cases{
1 & if $\Omega(a)\le 2$,
\cr\noalign{\smallskip}
0 & otherwise.
\cr}$$
Then it is easy to see
$$D_{1,2}(N)
\ge \sum_{a\in {\cal A}, \, (a, P(N^\kappa))=1} \delta^*(a)
= \sum_{a\in {\cal A}, \, (a,P(N^\kappa))=1} \mu(a)^2 \delta^*(a) + 
O(N^{1/2}),$$
where we have used the fact that
$$\sum_{\scriptstyle a\in {\cal A}\atop{\atop\scriptstyle (a, P(N^\kappa))=1}}
\{1 - \mu(a)^2\} \delta^*(a)
\le \sum_{N^\kappa\le p\le N^{1/2}} 1
\ll N^{1/2}.$$

Similarly if we write
$$\delta(a) := 1 - \dm s_1(a) - s_2(a) - \dm s_3(a) + \dm s_4(a),$$
we can show that
$$\sum_{a\in {\cal A}, \, (a, P(N^\kappa))=1} \delta(a)
= \sum_{a\in {\cal A}, \, (a, P(N^\kappa))=1} \mu(a)^2 \delta(a) + 
O(N^{1-\kappa}).$$

Thus in order to prove (9.2) it suffices to verify that
$$\delta^*(a)\ge \delta(a)
\leqno(9.3)$$
for $a\in {\cal A}$, $\mu(a)^2 = 1$ and $(a, P(N^\kappa)) = (a, N) = 1$.

We first observe that (9.3) is trivial if $\Omega(a)\le 2$,
since $\delta^*(a) = 1$ and $s_4(a)=0$ in this case.
It remains to show that $\delta(a)\le 0$ in all other cases,
which can be verified as follows:

If $\Omega(a) = 3$ and $s_1(a)=0$,
then $s_2(a) = 1$ and $s_3(a) = s_4(a) = 0$.
Thus $\delta(a)=0$.

If $\Omega(a) = 3$ and $s_1(a)=1$,
then $s_3(a) = 1$ and $s_2(a) = s_4(a) = 0$.
Thus $\delta(a)=0$.

If $\Omega(a) = 3$ and $s_1(a)=2$,
then $s_2(a) = s_3(a) = s_4(a) = 0$.
Thus $\delta(a)=0$.

If $\Omega(a) = 3$ and $s_1(a)=3$,
then $s_2(a) = s_3(a) = 0$ and $s_4(a) = 1$.
Thus $\delta(a)=0$.

If $\Omega(a) \ge 4$ and $s_1(a) = 1$, then $s_3(a) = 1$ and $s_2(a) 
= s_4(a) = 0$.
Thus $\delta(a)=0$.

If $\Omega(a) \ge 4$ and $s_1(a) = 2$, then $s_2(a) = s_3(a) = s_4(a) = 0$.
Thus $\delta(a)=0$.

If $\Omega(a) \ge 4$ and $s_1(a)\ge 3$, then $s_2(a) = s_3(a) = 0$ 
and $s_4(a)=s_1(a)-2$.
Thus $\delta(a)=0$.
This completes the proof.
\hfill
$\square$

\smallskip

The main difference between (9.1) and Chen's other weighted inequalities
(see (34) of [7] and page 425 of [8]) is the additional positive term $S_4$.
However a direct application of sieve to $S_4$ leads to zero contribution.
In order to take advantage of $S_4$, Chen used (9.1) with 
two different couples of parameters
$(\kappa, \sigma)$. Then an agreeable application of the Buchstab identity
and switching principle
leads to some compensation.
This idea was also used by Cai \& Lu [4] and Cai [3].
Here we make some modifications of their argument such that
this process is more powerful.

\proclaim Lemma 9.2.
Let $0<\kappa_1<\kappa_2<\rho<\sigma_2<\sigma_1<{1\over 3}$
such that $3\kappa_1+\rho\ge \dm$.
Then we have
$$\leqalignno{\qquad
4D_{1,2}(N)
& \ge
4S({\cal A}; {\cal P}(N), N^{\kappa_1})
- \Upsilon_1
- \Upsilon_2
- \Upsilon_3
+ \Upsilon_4
+ \Upsilon_5
+ \Upsilon_6
- 2 \Upsilon_7
& (9.4)
\cr
& \quad
- 2 \Upsilon_8
- \Upsilon_9
- \Upsilon_{10}
+ \Upsilon_{11}
+ \Upsilon_{12}
- \Upsilon_{13}
- \Upsilon_{14}
+ \Upsilon_{15}
+ O(N^{1-\kappa_1}),
\cr}$$
where
$$\eqalign{
& \Upsilon_1 := \sum_{\scriptstyle N^{\kappa_1}\le 
p<N^{\kappa_2}\atop\scriptstyle (p,N)=1}
S({\cal A}_p; {\cal P}(N), p),
\cr
& \Upsilon_2 := \sum_{\scriptstyle N^{\kappa_1}\le 
p<N^{\sigma_1}\atop\scriptstyle (p, N)=1}
S({\cal A}_p; {\cal P}(N), N^{\kappa_1}),
\cr
& \Upsilon_3 := \sum_{\scriptstyle N^{\kappa_1}\le 
p<N^{\sigma_2}\atop\scriptstyle (p,N)=1}
S({\cal A}_p; {\cal P}(N), N^{\kappa_1}),
\cr
& \Upsilon_4 := \mathop{\sum \,\, \sum }_{\scriptstyle N^{\kappa_1}\le 
p_1<p_2<N^{\kappa_2}
\atop\scriptstyle (p_1p_2,N)=1}
S({\cal A}_{p_1p_2}; {\cal P}(N), N^{\kappa_1}),
\cr
& \Upsilon_5 := \mathop{\sum \,\, \sum}_{\scriptstyle
N^{\kappa_1}\le p_1<N^{\kappa_2}\le 
p_2<N^{1/2-2\kappa_1}/p_1\atop\scriptstyle (p_1p_2,N)=1}
S({\cal A}_{p_1p_2}; {\cal P}(N), N^{\kappa_1}),
\cr
& \Upsilon_6
:= \mathop{\sum \,\,\sum}_{\scriptstyle
N^{\kappa_1}\le p_1<N^{\kappa_2}, N^{\rho}\le 
p_2<N^{\sigma_2}\atop\scriptstyle (p_1p_2,N)=1}
S({\cal A}_{p_1p_2}; {\cal P}(N), p_1),
\cr
& \Upsilon_7 :=
\mathop{\sum \,\, \sum}_{\scriptstyle N^{\sigma_1}\le p_1<p_2<(N/p_1)^{1/2}
\atop\scriptstyle (p_1p_2, N)=1}
S({\cal A}_{p_1 p_2}; {\cal P}(Np_1), p_2),
\cr
& \Upsilon_8 := \mathop{\sum \,\, \sum}_{\scriptstyle N^{\sigma_2}\le 
p_1<p_2<(N/p_1)^{1/2}
\atop{\atop\scriptstyle (p_1p_2, N)=1}}
S({\cal A}_{p_1 p_2}; {\cal P}(Np_1), p_2),
\cr
& \Upsilon_9 :=
\mathop{\sum \,\, \sum}_{\scriptstyle N^{\kappa_1}\le 
p_1<N^{\sigma_1}\le p_2<(N/p_1)^{1/2}
\atop\scriptstyle (p_1p_2, N)=1}
S({\cal A}_{p_1 p_2}; {\cal P}(Np_1), p_2),
\cr
& \Upsilon_{10}
:= \mathop{\sum \,\, \sum}_{\scriptstyle N^{\kappa_2}\le 
p_1<N^{\sigma_2}\le p_2<(N/p_1)^{1/2}
\atop\scriptstyle (p_1p_2, N)=1}
S({\cal A}_{p_1 p_2}; {\cal P}(Np_1), N^{\sigma_1}),
\cr
& \Upsilon_{11}
:= \mathop{\sum \,\, \sum \,\, \sum}_{\scriptstyle N^{\kappa_2}\le 
p_1<p_2<p_3<N^{\sigma_2}
\atop{\atop\scriptstyle (p_1p_2p_3, N)=1}}
S({\cal A}_{p_1 p_2 p_3}; {\cal P}(N), p_2),
\cr
& \Upsilon_{12} := \mathop{\sum \,\, \sum \,\, \sum}_{\scriptstyle 
N^{\kappa_2}\le p_1<N^{\sigma_2}, \,
N^{\sigma_1}\le p_2<p_3<(N/p_1)^{1/2}\atop\scriptstyle (p_1p_2p_3, N)=1}
S({\cal A}_{p_1 p_2 p_3}; {\cal P}(Np_1), p_2),
\cr
& \Upsilon_{13}
:= \mathop{\sum \,\, \sum \,\, \sum \,\, \sum}_{\scriptstyle
N^{\kappa_1}\le p_1<p_2<p_3<p_4<N^{\kappa_2} \atop\scriptstyle 
(p_1p_2p_3p_4,N)=1}
S({\cal A}_{p_1p_2p_3p_4}; {\cal P}(N), p_2),
\cr
& \Upsilon_{14}
:= \mathop{\sum \,\, \sum \,\, \sum \,\, \sum}_{\scriptstyle
N^{\kappa_1}\le p_1<p_2<p_3<N^{\kappa_2}\le p_4<N^{1/2-2\kappa_1}/p_3
\atop\scriptstyle (p_1p_2p_3p_4,N)=1}
S({\cal A}_{p_1p_2p_3p_4}; {\cal P}(N), p_2),
\cr
& \Upsilon_{15} := \mathop{\sum \,\, \sum \,\, \sum \,\, \sum}_{\scriptstyle
N^{\kappa_2}\le p_1<N^{\sigma_2}\le p_2<p_3<p_4<N^{\sigma_1}
\atop\scriptstyle (p_1p_2p_3p_4, N)=1}
S({\cal A}_{p_1p_2p_3p_4}; {\cal P}(N), p_3).
\cr}$$

\noindent{\sl Proof}.
The inequality (9.1) with $(\kappa, \sigma) = (\kappa_2, \sigma_2)$ implies
$$2 D_{1,2}(N)
\ge
2 S({\cal A}; {\cal P}(N), N^{\kappa_2})
- S_1(\kappa_2, \sigma_2)
- 2 \Upsilon_8
- S_3(\kappa_2, \sigma_2)
+ \Upsilon_{11}
+ O(N^{1-\kappa_2}).
\leqno(9.5)$$

Buchstab's identity, when applied three times, gives the equality
$$\eqalign{2S({\cal A}; {\cal P}(N), N^{\kappa_2})
& = 2 S({\cal A}; {\cal P}(N), N^{\kappa_1})
- \Upsilon_1
- \sum_{\scriptstyle N^{\kappa_1}\le p<N^{\kappa_2}\atop\scriptstyle (p,N)=1}
S({\cal A}_p; {\cal P}(N), N^{\kappa_1})
\cr
& \quad
+ \Upsilon_4
- \mathop{\sum \,\, \sum \,\, \sum}_{\scriptstyle N^{\kappa_1}\le 
p_1<p_2<p_3<N^{\kappa_2}
\atop\scriptstyle (p_1p_2p_3,N)=1}
S({\cal A}_{p_1p_2p_3}; {\cal P}(N), p_1).
\cr}$$
Similarly a simple application of Buchstab's identity yields
$$S_1(\kappa_2, \sigma_2)
= \sum_{\scriptstyle N^{\kappa_2}\le p<N^{\sigma_2}\atop\scriptstyle (p,N)=1}
S({\cal A}_p; {\cal P}(N), N^{\kappa_1})
- \mathop{\sum \,\, \sum}_{\scriptstyle N^{\kappa_1}\le 
p_1<N^{\kappa_2}\le p_2<N^{\rho}
\atop\scriptstyle (p_1p_2,N)=1}
S({\cal A}_{p_1p_2}; {\cal P}(N), p_1)
- \Upsilon_6.$$
Clearly $p_1<N^{\sigma_2}$ and $\sigma_2<\sigma_1<{1\over 3}$ imply 
that $N^{\sigma_1}<(N/p_1)^{1/2}$.
Thus by Buchstab's identity, we can write
$$\eqalign{S_3(\kappa_2, \sigma_2)
& = \mathop{\sum \,\, \sum}_{\scriptstyle N^{\kappa_2}\le 
p_1<N^{\sigma_2}\le p_2<N^{\sigma_1}
\atop\scriptstyle (p_1p_2, N)=1}
S({\cal A}_{p_1 p_2}; {\cal P}(Np_1), N^{\sigma_1})
\cr
& \quad
+ \mathop{\sum \,\, \sum \,\, \sum}_{\scriptstyle
N^{\kappa_2}\le p_1<N^{\sigma_2}\le 
p_2<p_3<N^{\sigma_1}\atop\scriptstyle (p_1p_2p_3, N)=1}
S({\cal A}_{p_1 p_2 p_3}; {\cal P}(Np_1), p_3)
+ O(N^{1-\sigma_2})
\cr
& \quad
+ \mathop{\sum \,\, \sum}_{\scriptstyle
N^{\kappa_2}\le p_1<N^{\sigma_2}, \, N^{\sigma_1}\le 
p_2<(N/p_1)^{1/2}\atop\scriptstyle (p_1p_2, N)=1}
S({\cal A}_{p_1 p_2}; {\cal P}(Np_1), N^{\sigma_1})
\cr
& \quad
- \mathop{\sum \,\, \sum \,\, \sum}_{\scriptstyle N^{\kappa_2}\le 
p_1<N^{\sigma_2}, \,
N^{\sigma_1}\le p_2<p_3<(N/p_1)^{1/2}\atop\scriptstyle (p_1p_2p_3, N)=1}
S({\cal A}_{p_1 p_2 p_3}; {\cal P}(Np_1), p_2)
\cr
& = \Upsilon_{10} - \Upsilon_{12}
+ \mathop{\sum \,\, \sum \,\, \sum}_{\scriptstyle
N^{\kappa_2}\le p_1<N^{\sigma_2}\le 
p_2<p_3<N^{\sigma_1}\atop\scriptstyle (p_1p_2p_3, N)=1}
S({\cal A}_{p_1 p_2 p_3}; {\cal P}(Np_1), p_3)
+ O(N^{1-\sigma_2}).
\cr}$$
Inserting these into (9.5), we find that
$$2D_{1,2}(N)
\ge
2S({\cal A}; {\cal P}(N), N^{\kappa_1})
- \Upsilon_1
- \Upsilon_3
+ \Upsilon_4
+ \Upsilon_6
- 2 \Upsilon_8
- \Upsilon_{10}
+ \Upsilon_{11}
+ \Upsilon_{12}
+ \Delta
+ O(N^{1-\kappa_2}),
\leqno(9.6)$$
where
$$\eqalign{\Delta
&  :=
- \mathop{\sum \,\, \sum \,\, \sum}_{\scriptstyle
N^{\kappa_1}\le p_1<p_2<p_3<N^{\kappa_2}\atop\scriptstyle (p_1p_2p_3, N)=1}
S({\cal A}_{p_1 p_2p_3}; {\cal P}(Np_1), p_1)
\cr
& \hskip 4,8mm
- \mathop{\sum \,\, \sum \,\, \sum}_{\scriptstyle
N^{\kappa_2}\le p_1<N^{\sigma_2}\le 
p_2<p_3<N^{\sigma_1}\atop\scriptstyle (p_1p_2p_3, N)=1}
S({\cal A}_{p_1 p_2p_3}; {\cal P}(Np_1), p_3)
\cr
& \hskip 4,8mm
+ \mathop{\sum \,\, \sum}_{\scriptstyle N^{\kappa_1}\le 
p_1<N^{\kappa_2}\le p_2<N^{\rho}
\atop\scriptstyle (p_1p_2,N)=1}
S({\cal A}_{p_1p_2}; {\cal P}(N), p_1).
\cr}$$

Next we shall further decompose $\Delta$.
In view of $3\kappa_1+\rho\ge \dm$,
we have $N^{\rho}\ge N^{1/2-2\kappa_1}/p_1$ provided $p_1\ge N^{\kappa_1}$.
Thus Buchstab's identity allows us to write
$$\eqalign{\mathop{\sum \,\, \sum}_{\scriptstyle N^{\kappa_1}\le 
p_1<N^{\kappa_2}\le p_2<N^{\rho}
\atop\scriptstyle (p_1p_2,N)=1} \!\!\!\!
S({\cal A}_{p_1p_2}; {\cal P}(N), p_1)
& \ge \Upsilon_5
\cr
& - \mathop{\sum \,\, \sum \,\, \sum}_{\scriptstyle
N^{\kappa_1}\le p_1<p_2<N^{\kappa_2}\le p_3<N^{1/2-2\kappa_1}/p_2
\atop\scriptstyle (p_1p_2p_3,N)=1} \!\!\!\!
S({\cal A}_{p_1p_2p_3}; {\cal P}(N), p_1).
\cr}$$
Thus we have
$$\Delta\ge \Upsilon_5 + \Delta_1,$$
where
$$\eqalign{\Delta_1
&  :=
- \mathop{\sum \,\, \sum \,\, \sum}_{\scriptstyle N^{\kappa_1}\le 
p_1<p_2<p_3<N^{\kappa_2}
\atop\scriptstyle (p_1p_2p_3,N)=1}
S({\cal A}_{p_1p_2p_3}; {\cal P}(N), p_1)
\cr
& \hskip 4,8mm
- \mathop{\sum \,\, \sum \,\, \sum}_{\scriptstyle
N^{\kappa_1}\le p_1<p_2<N^{\kappa_2}\le p_3<N^{1/2-2\kappa_1}/p_2
\atop\scriptstyle (p_1p_2p_3,N)=1}
S({\cal A}_{p_1p_2p_3}; {\cal P}(N), p_1)
\cr
& \hskip 4,8mm
- \mathop{\sum \,\, \sum \,\, \sum}_{\scriptstyle
N^{\kappa_2}\le p_1<N^{\sigma_2}\le p_2<p_3<N^{\sigma_1}
\atop\scriptstyle (p_1p_2p_3, N)=1}
S({\cal A}_{p_1 p_2p_3}; {\cal P}(Np_1), p_3).
\cr}$$
Now the inequality (9.6) becomes
$$\leqalignno{2 D_{1,2}(N)
& \ge
2 S({\cal A}; {\cal P}(N), N^{\kappa_1})
- \Upsilon_1
- \Upsilon_3
+ \Upsilon_4
+ \Upsilon_5
+ \Upsilon_6
& (9.7)
\cr\noalign{\smallskip}
& \quad
- 2 \Upsilon_8
- \Upsilon_{10}
+ \Upsilon_{11}
+ \Upsilon_{12}
+ \Delta_1
+ O(N^{1-\kappa_2}).
\cr}$$

The inequality (9.1) with $(\kappa, \sigma) = (\kappa_1, \sigma_1)$ gives us
$$2 D_{1,2}(N)
\ge
2 S({\cal A}; {\cal P}(N), N^{\kappa_1})
- \Upsilon_2
- 2 \Upsilon_7
- \Upsilon_9
+ S_4(\kappa_1, \sigma_1)
+ O(N^{1-\kappa_1}).
\leqno(9.8)$$
Adding (9.7) to (9.8) yields
$$\leqalignno{\qquad
4 D_{1,2}(N)
& \ge
4 S({\cal A}; {\cal P}(N), N^{\kappa_1})
- \Upsilon_1
- \Upsilon_2
- \Upsilon_3
+ \Upsilon_4
+ \Upsilon_5
+ \Upsilon_6
& (9.9)
\cr\noalign{\smallskip}
& \quad
- 2 \Upsilon_7
- 2 \Upsilon_8
- \Upsilon_9
- \Upsilon_{10}
+ \Upsilon_{11}
+ \Upsilon_{12}
+ \Delta_2
+ O(N^{1-\kappa_1}),
\cr}$$
where
$$\Delta_2
:= \Delta_1
+ \mathop{\sum \,\, \sum \,\, \sum}_{\scriptstyle N^{\kappa_1}\le 
p_1<p_2<p_3<N^{\sigma_1}
\atop\scriptstyle (p_1p_2p_3, N)=1}
S({\cal A}_{p_1 p_2 p_3}; {\cal P}(N), p_2).$$
Clearly all domains of summation in the three terms on the right-hand 
side of $\Delta_1$ are distinct
and are contained in the domain of summation of the last triple sums
on the right-hand side of $\Delta_2$ (since 
$3\kappa_1+\sigma_1>3\kappa_1+\rho\ge {1\over 2}$).
Therefore we have
$$\eqalign{\Delta_2
& \ge - \mathop{\sum \,\, \sum \,\, \sum}_{\scriptstyle 
N^{\kappa_1}\le p_1<p_2<p_3<N^{\kappa_2}
\atop\scriptstyle (p_1p_2p_3,N)=1}
\big\{S({\cal A}_{p_1p_2p_3}; {\cal P}(N), p_1) - S({\cal 
A}_{p_1p_2p_3}; {\cal P}(N), p_2)\big\}
\cr
& \hskip 4,8mm
- \mathop{\sum \,\, \sum \,\, \sum}_{\scriptstyle
N^{\kappa_1}\le p_1<p_2<N^{\kappa_2}\le p_3<N^{1/2-2\kappa_1}/p_2
\atop\scriptstyle (p_1p_2p_3,N)=1}
\big\{S({\cal A}_{p_1p_2p_3}; {\cal P}(N), p_1) - S({\cal 
A}_{p_1p_2p_3}; {\cal P}(N), p_2)\big\}
\cr
& \hskip 4,8mm
+ \mathop{\sum \,\, \sum \,\, \sum}_{\scriptstyle
N^{\kappa_2}\le p_1<N^{\sigma_2}\le p_2<p_3<N^{\sigma_1}
\atop\scriptstyle (p_1p_2p_3, N)=1}
\big\{S({\cal A}_{p_1p_2p_3}; {\cal P}(N), p_2) - S({\cal 
A}_{p_1p_2p_3}; {\cal P}(N), p_3)\big\}
\cr
& = - \Upsilon_{13} - \Upsilon_{14} + \Upsilon_{15} + O(N^{1-\kappa_1}).
\cr}$$
Combining this with (9.9), we obtain the required result.
This completes the proof.
\hfill
$\square$

\medskip

{\bf Remark 3.}
Lemmas 9.1 and 9.2 are also valid for
$$\eqalign{
& {\cal A}' := \{p+2 : p\le x\},
\cr
& {\cal A}'' := \{p+2 : x<p\le x+x^\theta\},
\cr
& {\cal A}''' := \{N-p : \alpha N<p\le \alpha N+N^\theta\},
\cr}$$
if we make some suitable modifications.
For example, we have
$$\leqalignno{4 \pi_{1,2}(x)
& \ge
4 S({\cal A}'; {\cal P}(2), x^{\kappa_1})
- \Upsilon_1'
- \Upsilon_2'
- \Upsilon_3'
+ \Upsilon_4'
+ \Upsilon_5'
+ \Upsilon_6'
- 2 \Upsilon_7'
& (9.10)
\cr
& \quad
- 2 \Upsilon_8'
- \Upsilon_9'
- \Upsilon_{10}'
+ \Upsilon_{11}'
+ \Upsilon_{12}'
- \Upsilon_{13}'
- \Upsilon_{14}'
+ \Upsilon_{15}'
+ O(N^{1-\kappa_1}),
\cr}$$
where $\Upsilon_j'$ is similarly defined as $\Upsilon_j$ 
with the difference that
${\cal A}$ is replaced by ${\cal A}'$,
${\cal P}(N)$ by ${\cal P}(2)$,
${\cal P}(Np_1)$ by ${\cal P}(2p_1)$,
$(N/p_1)^{1/2}$ by $(x/p_1)^{1/2}$,
$N^{1/2-2\kappa_1}/p_2$ by $x^{4/7-2\kappa_1}/p_2$ 
(in $\Upsilon_5$ and $\Upsilon_{14}$),
$N^{\rho}$ by $x^{\rho}$,
$N^{\kappa_i}$ by $x^{\kappa_i}$,
$N^{\sigma_i}$ by $x^{\sigma_i}$
and that the conditions
$(p,N)=1$, $(p_1p_2,N)=1$, $(p_1p_2p_3,N)=1$ and $(p_1p_2p_3p_4,N)=1$ 
are eliminated.
The assumption on the parameters is
$$\textstyle 0<\kappa_1<\kappa_2<\rho<\sigma_2<\sigma_1<{1\over 3},
\qquad
3\kappa_1+\rho\ge {4\over 7}.$$
The last condition is necessary 
in the proof of $\Delta\ge \Upsilon_5'+\Delta_1$.

\vskip 5mm

\noindent{\bf \S\ 10. Proofs of Theorems 2 and 5}

\medskip

For simplicity, we write
${\cal L} := \log N$ and
use $B$ to denote a suitable positive constant determined by Lemma 2.3.
We shall estimate all terms $\Upsilon_i$ in the inequality (9.4).
For this we suppose that
$$\textstyle {1\over 12}
= \kappa_1 < \kappa_2\le {1\over 8},
\quad
{1\over 4}
= \rho
<\sigma_2<\sigma_1<{1\over 3},
\quad
3\sigma_1+\kappa_1>1,
\quad
2\sigma_1+\sigma_2+\kappa_2>1.
\leqno(10.1)$$

\smallskip

$1^\circ$
{\it Lower bound of $S({\cal A}; {\cal P}(N), N^{\kappa_1})$}

\smallskip

We apply (2.5) of Lemma 2.2 with
$$X = {\rm li}(N),
\quad
w(p)=\cases{
p/\varphi(p) &  if $p\in {\cal P}(N)$,
\cr\noalign{\smallskip}
0       & otherwise,
\cr}
\quad
z = N^{\kappa_1},
\quad
Q = {\sqrt N\over {\cal L}^B}.$$
Since $|\lambda_l^\pm(q)|\le 1$,
Lemma 2.3 with the choice $f(1) = 1$ and $f(m)=0$ if $m\ge 2$ implies that
$$\Big|\sum_{l<L} \sum_{q\mid P(z)} \lambda_l^\pm(q) r({\cal A}, q)\Big|
\ll_\varepsilon \sum_{q\le \sqrt N/{\cal L}^B} \mu(q)^2
\max_{y\le N} \max_{(a, q)=1}
\bigg|\pi(y; q,a)-{{\rm li}(y)\over \varphi(q)}\bigg|
\ll_\varepsilon {N\over {\cal L}^3}.$$
In view of $V(z)\sim 2 e^{-\gamma} C_N/\log z$ ($\gamma$ is the Euler 
constant) and $C_N\gg 1$,
we can deduce
$$S({\cal A}; {\cal P}(N), N^{\kappa_1})
\ge \{F_0 + O(\varepsilon)\} \, \Theta(N)
\qquad{\rm with}\qquad
F_0 := 2 f(1/2\kappa_1)/\kappa_1 e^\gamma.
\leqno(10.2)$$

$2^\circ$
{\it Upper bounds of $\Upsilon_1$, $\Upsilon_2$ and $\Upsilon_3$}

\smallskip

We apply (2.4) of Lemma 2.2 with
$$X = {{\rm li}(N)\over \varphi(p)},
\quad
w(p)=\cases{
p/\varphi(p)     &  if $p\in {\cal P}(N)$,
\cr\noalign{\smallskip}
0                & otherwise,
\cr}
\quad
z = p,
\quad
Q = {\sqrt N\over p{\cal L}^B}$$
to $S({\cal A}_p; {\cal P}(N), p)$.
The contribution of the error term in (2.4) is
$$\eqalign{
& \ll_\varepsilon \sum_{N^{\kappa_1}\le p<N^{\kappa_2}, \, (p, N)=1}
\sum_{q\le \sqrt N/p {\cal L}^B, \, q\mid P(p)} |r({\cal A}, pq)|
\cr
& \ll_\varepsilon \sum_{d\le \sqrt N/{\cal L}^B}
\mu(d)^2 \max_{y\le N} \max_{(a, d)=1}
\bigg|\pi(y; d, a) - {{\rm li}(y)\over \varphi(d)}\bigg|
\ll_\varepsilon {N\over {\cal L}^3}
\cr}$$
by Lemma 2.3 with the same choice of $f$ as above.
Thus
$$\Upsilon_1 \le {\{1+O(\varepsilon)\} N\over {\cal L}}
\sum_{N^{\kappa_1}\le p<N^{\kappa_2}} {V(p)\over \varphi(p)}
F\bigg({\log(\sqrt N/p)\over \log p}\bigg)
+ O\bigg({N\over {\cal L}^3}\bigg).$$
The standard procedure for replacing sums over primes by integrals yields
$$\Upsilon_1
\le \{F_1 + O(\varepsilon)\} \, \Theta(N),
\leqno(10.3)$$
where
$$F_1
:= {2\over e^\gamma} \int_{\kappa_1}^{\kappa_2} {F(1/2t-1)\over t^2} \d t
= {4\over e^\gamma} \int_{1/2\kappa_2-1}^{1/2\kappa_1-1} F(t) \d t.$$

Similarly we can prove
$$\Upsilon_i
\le \{F_i + O(\varepsilon)\} \, \Theta(N)
\qquad
(i = 2, 3),
\leqno(10.4)$$
where
$$F_2
:= {4\over e^\gamma}
\int_{(1/2-\sigma_1)/\kappa_1}^{(1/2-\kappa_1)/\kappa_1} {F(t)\over 
1-2\kappa_1t}  \d t,
\qquad
F_3
:= {4\over e^\gamma}
\int_{(1/2-\sigma_2)/\kappa_1}^{(1/2-\kappa_1)/\kappa_1} {F(t)\over 
1-2\kappa_1t}  \d t.$$

$3^\circ$
{\it Lower bounds of $\Upsilon_4$ and $\Upsilon_5$}

\smallskip

As before we can deduce, from Lemmas 2.2 and 2.3, that
$$\Upsilon_i
\ge \{F_i + O(\varepsilon)\} \, \Theta(N)
\qquad
(i = 4, 5),
\leqno(10.5)$$
where
$$\eqalign{
F_4
& := {2\over \kappa_1 e^\gamma}
\int_{\kappa_1}^{\kappa_2} {\d t\over t}
\int_{t}^{\kappa_2} f\bigg({1/2-t-u\over \kappa_1}\bigg) {\d u\over u},
\cr
F_5
& := {2\over \kappa_1 e^\gamma}
\int_{\kappa_1}^{\kappa_2} {\d t\over t}
\int_{\kappa_2}^{1/2-2\kappa_1-t} f\bigg({1/2-t-u\over 
\kappa_1}\bigg) {\d u\over u}.
\cr}$$
We have used the following fact to remove the condition $(p_1p_2,N)=1$:
$$\sum_{\scriptstyle p_1|N\atop\scriptstyle p_1\ge N^\kappa}
\sum_{N^\kappa\le p_2<N^{\sigma}} {N\over p_1p_2}
+ \sum_{N^\kappa\le p_1<N^{\sigma}}
\sum_{\scriptstyle p_2|N\atop\scriptstyle p_2\ge N^\kappa} {N\over p_1p_2}
\ll N^{1-\kappa} {\cal L}.$$

$4^\circ$
{\it Upper bounds of $\Upsilon_i$ for $i = 7, 8, 9, 10, 13, 14$}

\smallskip

We shall only majorize $\Upsilon_7$ and the others can be treated similarly.

Since $\sigma_1>{1\over 4}$,
the quantity $\Upsilon_7$ is equal to the number of primes $p\le N$
such that $N-p=p_1p_2p_3$
with $N^{\sigma_1} \le p_1<p_2< (N/p_1)^{1/2}$, $p_3\ge p_2$ and 
$(p_1p_2p_3, N)=1$.
Define
$$\eqalign{
{\cal M}
& := \{m : m = p_1p_2, \,
N^{\sigma_1} \le p_1<p_2< (N/p_1)^{1/2}, \, (p_1p_2, N)=1\},
\cr
{\cal B}
& := \{b : b=N-mp\le N, \, m\in {\cal M},  \, p\le N/m\}.
\cr}$$
It is clear that
$$\Upsilon_7
\le S({\cal B}; {\cal P}(N), N^{1/2}) + O(N^{1/2}).$$
By applying (2.4) of Lemma 2.2 with
$$
X = \sum_{m\in {\cal M}} {\rm li}\Big({N\over m}\Big),
\qquad
w(p) = \cases{
p/\varphi(p)        & if $p\in {\cal P}(N)$,
\cr\noalign{\smallskip}
0                   & otherwise,
\cr}
\qquad
Q = {\sqrt N\over {\cal L}^B},$$
we obtain
$$\Upsilon_7\le {8 \, C_N X\over \log N} \{1 + O(\varepsilon)\}
+ O_\varepsilon\big(\sqrt N + R_3 + R_4\big),
\leqno(10.6)$$
where
$$\eqalign{
& R_3 := \sum_{\scriptstyle q\le \sqrt N/{\cal L}^B 
\atop{\atop\scriptstyle (q,N)=1}}
\mu(q)^2
\bigg|\sum_{\scriptstyle m\in{\cal M} \atop{\atop\scriptstyle (q, m)=1}}
\bigg(\sum_{\scriptstyle mp\le N \atop{\atop\scriptstyle
mp\equiv N ({\rm mod} \, q)}} 1
- {{\rm li}(N/m)\over \varphi(q)}
\bigg)\bigg|,
\cr\
& R_4 := \sum_{q\le \sqrt N/{\cal L}^B, \, (q,N)=1}
{\mu(q)^2\over \varphi(q)}
\sum_{m\in {\cal M}, \, (q, m)>1}
{\rm li} \Big({N\over m}\Big).
\cr}$$
Let $f(m)$ be the characteristic function of ${\cal M}$.
Since $m\le N^{3/4}$ for $m\in {\cal M}$, Lemma 2.3 implies
$$R_3
= \sum_{\scriptstyle q\le \sqrt N/{\cal L}^B \atop{\atop\scriptstyle (q,N)=1}}
\mu(q)^2
\bigg|\sum_{\scriptstyle m\le N^{5/6}\atop{\atop\scriptstyle (q, m)=1}} f(m)
\bigg(\sum_{\scriptstyle mp\le N \atop{\atop\scriptstyle
mp\equiv N ({\rm mod} \, q)}} 1
-{{\rm li}(N/m)\over \varphi(q)}\bigg)\bigg|
\ll {N\over {\cal L}^3}.
\leqno(10.7)$$

Noticing that $(d, m)>1$ implies $(d, m)\ge N^{\sigma_1}$ for $m\in 
{\cal M}$, we have
$$\eqalign{R_4
& \ll {N\over {\cal L}} \sum_{q\le \sqrt N} {\mu(q)^2\over \varphi(q)}
\sum_{m\le N^{3/4}, \, (q, m)\ge N^{\sigma_1}} {1\over m}
\cr
& \ll {N\over {\cal L}} \sum_{q\le N} {\mu(q)^2\over \varphi(q)}
\sum_{d\mid q, \, d\ge N^{\sigma_1}} {1\over d}
\sum_{n\le N^{3/4}/d} {1\over n}
\cr
& \ll N \sum_{q\le \sqrt N} {\mu(q)^2\over \varphi(q)}
\sum_{d\mid q, \, d\ge N^{\sigma_1}} {1\over d}
\cr
& \ll N \sum_{N^{\sigma_1}<d\le N} {1\over d}
\sum_{l\le \sqrt N/d} {\mu(dl)^2\over \varphi(dl)}
\cr
& \ll N \sum_{N^{\sigma_1}<d\le N} {\mu(d)^2\over d\varphi(d)}
\sum_{l\le \sqrt N/d} {\mu(l)^2\over \varphi(l)}.
\cr}$$
Since the function $\mu(n)^2/\varphi(n)$ is multiplicative and
$\mu(p^\nu)^2/\varphi(p^\nu) = 1/(p-1)$ for $\nu=1$ and $=0$ for $\nu\ge 2$,
it is plain to see that
$$\sum_{l\le t} {\mu(n)^2\over \varphi(n)}\asymp \log t.$$
Thus
$$R_4\ll N^{1-\sigma_1} {\cal L}^2.
\leqno(10.8)$$

By the prime number theorem, we obtain
$$\eqalign{X
& = \{1 + o(1)\} \mathop{\sum \,\, \sum}_{N^{\sigma_1}\le 
p_1<p_2<(N/p_1)^{1/2}}
{N\over p_1 p_2\log(N/p_1p_2)}
\cr
& = \{1 + o(1)\} {N\over {\cal L}} \int_{\sigma_1}^{1/3} 
{\log(1/t-2)\over t(1-t)} \d t
\cr
& = \{1 + o(1)\} {N\over {\cal L}} \int_2^{1/\sigma_1-1} 
{\log(t-1)\over t} \d t.
\cr}
\leqno(10.9)$$
Inserting (10.7)--(10.9) into (10.6) yields
$$\Upsilon_7
\le \{F_7 + O(\varepsilon)\} \, \Theta(N),
\leqno(10.10)$$
where
$$F_7 := 8 \int_2^{1/\sigma_1-1} {\log(t-1)\over t} \d t.$$

Similarly we can prove that
$$\Upsilon_i
\le \{F_i + O(\varepsilon)\} \, \Theta(N)
\qquad
(i = 8, 10),
\leqno(10.11)$$
where
$$\eqalign{
F_8
& := 8 \int_2^{1/\sigma_2-1} {\log(t-1)\over t} \d t,
\cr
F_{10}
& := 8
\int_{\kappa_2}^{\sigma_2}{\log(1/\sigma_2-1-t/\sigma_2)\over t (1-t)} \d t.
\cr}$$
[We need to use the assumption $2\sigma_1+\sigma_2+\kappa_2>1$ in 
$\Upsilon_{10}$.]

For the terms $\Upsilon_9, \Upsilon_{13}$ and $\Upsilon_{14}$ with 
$p_1\le N^{1/10}$,
we can apply Lemma 2.6 instead of Lemma 2.3.
A similar argument allows us to show that
$$\Upsilon_i
\le \{F_i + O(\varepsilon)\} \, \Theta(N)
\qquad
(i = 9, 13, 14),
\leqno(10.12)$$
where
$$\eqalign{
F_9
& := {36\over 5} \int_{\kappa_1}^{1/10} 
{\log(1/\sigma_1-1-t/\sigma_1)\over t (1-t)^2} \d t
+ 8 \int_{1/10}^{\sigma_1} {\log(1/\sigma_1-1-t/\sigma_1)\over t (1-t)} \d t,
\cr
F_{13}
& := {36\over 5}
\int_{\kappa_1}^{1/10} {\d t_1\over t_1(1-t_1)}
\int_{t_1}^{\kappa_2} {\d t_2\over t_2^2}
\int_{t_2}^{\kappa_2} {\d t_3\over t_3}
\int_{t_3}^{\kappa_2} \omega\bigg({1-t_1-t_2-t_3-t_4\over t_2}\bigg) 
{\d t_4\over t_4}
\cr
& \quad
+ 8 \int_{1/10}^{\kappa_2} {\d t_1\over t_1}
\int_{t_1}^{\kappa_2} {\d t_2\over t_2^2}
\int_{t_2}^{\kappa_2} {\d t_3\over t_3}
\int_{t_3}^{\kappa_2} \omega\bigg({1-t_1-t_2-t_3-t_4\over t_2}\bigg) 
{\d t_4\over t_4},
\cr
F_{14}
& := {36\over 5}
\int_{\kappa_1}^{1/10} {\d t_1\over t_1(1-t_1)}
\int_{t_1}^{\kappa_2} {\d t_2\over t_2^2}
\int_{t_2}^{\kappa_2} {\d t_3\over t_3}
\int_{\kappa_2}^{1/2-2\kappa_1-t_3} 
\omega\bigg({1-t_1-t_2-t_3-t_4\over t_2}\bigg) {\d t_4\over t_4}
\cr
& \quad
+ 8 \int_{1/10}^{\kappa_2} {\d t_1\over t_1}
\int_{t_1}^{\kappa_2} {\d t_2\over t_2^2}
\int_{t_2}^{\kappa_2} {\d t_3\over t_3}
\int_{\kappa_2}^{1/2-2\kappa_1-t_3} \omega\bigg({1-t_1-t_2-t_3-t_4\over 
t_2}\bigg) {\d t_4\over t_4}.
\cr}$$
[We need to use the assumption $3\sigma_1+\kappa_1>1$ in $\Upsilon_9$
and Lemma 2.10 in $\Upsilon_{13}$ and $\Upsilon_{14}$.]

By inserting (10.2)--(10.5), (10.10)--(10.12) and
by using the trivial lower bounds $\Upsilon_i\ge 0$ ($i = 6, 11, 12, 15$) into (9.4),
we get the following inequality
$$D_{1,2}(N)
\ge \{F(\kappa_1, \kappa_2, \rho, \sigma_2, \sigma_1) + 
O(\varepsilon)\} \, \Theta(N),$$
where
$$\eqalign{F(\kappa_1, \kappa_2, \rho, \sigma_2, \sigma_1)
& := \textstyle {1\over 4} (4F_0
- F_1
- F_2
- F_3
+ F_4
+ F_5
\cr
& \quad
- 2F_7
- 2F_8
- F_9
- F_{10}
- F_{13}
- F_{14}).
\cr}$$
Taking
$\kappa_1 = {1\over 12}$,
$\kappa_2 = {29\over 250}$,
$\rho = {1\over 4}$,
$\sigma_2 = {141\over 500}$
and
$\sigma_1 = {41\over 125}$,
a numerical computation gives us
$$\eqalign{\textstyle F(\kappa_1, \kappa_2, \rho, \sigma_2, \sigma_1)
& \ge \textstyle {1\over 4}(
   4\times 13.473613
- 3.891854
- 20.432098
- 17.327241
\cr
& \quad
+ 0.697375
+ 2.118119
- 2\times 0.004609
- 2\times 0.434368
\cr
& \quad
- 5.161945
- 5.468377
- 0.023310
- 0.182860)
\cr
& > 0.83607.
\cr}$$
[For the integrals $F_{13}$ and $F_{14}$, we make use of 
$\omega(u)\le 0.561522$ for $u\ge 3.5$.]
This completes the proof of Theorem 2.
\hfill
$\square$

\smallskip

Theorem 5 can be proved in the same way.
The only difference is to replace Lemmas 2.3 and 2.6 by Lemma 2.4.
Here, the choice of parameters is
$$(\theta, \kappa_1, \kappa_2, \rho, \sigma_2, \sigma_1)
= \textstyle
(0.971, \, {(2\theta-1)/12}, \, 0.111, \, {(2\theta-1)/4}, \, 0.271, 0.313).
\eqno\square$$

\vskip 5mm

\noindent{\bf \S\ 11. Proof of Theorem 4}

\medskip

The proof of Theorem 4 is very similar to that of Theorem 2.
But we must use Lemmas 2.5, 2.7 and 2.8 in place of Lemmas 2.3 and 2.6.
In order to take the advantage of these lemmas,
we must carry out a more careful and delicate analysis.
Thus the proof will be slightly complicated.

Suppose that the parameters satisfy the following conditions:
$$\textstyle
  {2\over 21}
= \kappa_1
<\kappa_2\le {1\over 7},
\quad
{2\over 7}
= \rho<\sigma_2<{29\over 100}<\sigma_1<{1\over 3},
\quad
3\sigma_1+\kappa_1>1,
\quad
2\sigma_1+\sigma_2+\kappa_2>1.$$

\smallskip

$1^\circ$
{\it Lower bounds of $S({\cal A}'; {\cal P}(2), x^{\kappa_1})$}

\smallskip

By (2.5) of Lemma 2.2 and Lemma 2.7, we can easily prove
$$S({\cal A}'; {\cal P}(2), x^{\kappa_1})
\ge \{G_0 + O(\varepsilon)\} \, \Pi(x)
\quad{\rm with}\quad
G_0 := f(4/7\kappa_1)/\kappa_1 e^\gamma.
\leqno(11.1)$$

\smallskip

$2^\circ$
{\it Upper bounds of $\Upsilon_1'$, $\Upsilon_2'$ and $\Upsilon_3'$}

\smallskip

We divide the interval $[x^{\kappa_1}, x^{\kappa_2}]$ into $O({\cal L})$
subintervals of the form $[P, 2P)$
and apply (2.4) of Lemma 2.2 to $S({\cal A}_p';{\cal P}(2),p)$ for 
$p\in [P, 2P)$.
We have
$$S({\cal A}_p';{\cal P}(2),p)
\le {\{1 + O(\varepsilon)\} x\over {\cal L}} {V(p)\over \varphi(p)} 
F\bigg({\log(Q/P)\over \log p}\bigg)
+ \sum_{l<L} \sum_{q\mid P(p)} \lambda_l^+(q) r({\cal A'}, pq)$$
where $Q=x^{4/7-\varepsilon}$ and $\lambda_l^+(q)$ is well factorable 
of level $Q/P$ and of order 1.

Denote by $\pi_P$ the characteristic function of the primes in the 
interval $[P, 2P)$.
Noticing that $P\le x^{\kappa_2}$ $\Rightarrow$ $P\le Q/P$,
Lemma 2.1 shows that $\pi_P*\lambda_l^+$ is well factorable of level 
$Q$ and of order 2.
Thus Lemma 2.7 allows us to deduce that
$$\sum_{P\le p<2P} \sum_{l<L} \sum_{q\mid P(p)} \lambda_l^+(q) r({\cal A'}, pq)
\ll_\varepsilon x/(\log x)^4$$
and
$$\Upsilon_1'
\le {\{1 + O(\varepsilon)\} x\over \log x}
\sum_{x^{\kappa_1}\le p<x^{\kappa_2}} {V(p)\over \varphi(p)} 
F\bigg({\log(Q/p)\over \log p}\bigg)
+ O\bigg({x\over (\log x)^4}\bigg).$$
Since $V(p)\sim e^{-\gamma} C/\log p$,
the prime number theorem implies that
$$\Upsilon_1'
\le \{G_1 + O(\varepsilon)\} \, \Pi(x)
\qquad{\rm with}\qquad
G_1 := {7\over 4 e^\gamma}
\int_{4/7\kappa_2-1}^{4/7\kappa_1-1} F(t) \d t.
\leqno(11.2)$$

We divide the interval of summation $[x^{\kappa_1}, x^{\sigma_1}]$ of 
$\Upsilon_2'$ into three parts:
$$[x^{\kappa_1}, x^{2/7-\varepsilon}],
\qquad
[x^{2/7-\varepsilon}, x^{29/100}],
\qquad
[x^{29/100}, x^{\sigma_1}],$$
and use (2.4) of Lemma 2.2 to handle each sum.
As before we apply Lemma 2.7, 
the condition (C.2) and (C.3) of Lemma 2.8, respectively,
to control the corresponding error terms.
We find
$$\Upsilon_2'\le \{G_2 + O(\varepsilon)\} \, \Pi(x),
\leqno(11.3)$$
where
$$G_2
:= {1\over \kappa_1 e^\gamma} \bigg\{
\int_{2/7\kappa_1}^{4/7\kappa_1-1} {F(t)\over 4/7\kappa_1-t} \d t
+ \int_{13/50\kappa_1}^{2/7\kappa_1} {F(t)\over 2/\kappa_1-t}  \d t
+ \int_{(11/20-\sigma_1)/\kappa_1}^{13/50\kappa_1} {F(t)\over 
11/20\kappa_1-t} \d t\bigg\}.$$

Similarly
$$\Upsilon_3'\le \{G_3 + O(\varepsilon)\} \, \Pi(x),$$
where
$$G_3
:= {1\over \kappa_1 e^\gamma} \bigg\{
\int_{2/7\kappa_1}^{(4/7-\kappa_1)/\kappa_1} {F(t)\over 4/7\kappa_1-t} \d t
+ \int_{(2-6\sigma_2)/\kappa_1}^{2/7\kappa_1} {F(t)\over 2/\kappa_1-t} \d t
\bigg\}.
\leqno(11.4)$$

\smallskip

$3^\circ$
{\it Lower bounds of $\Upsilon_4'$ and $\Upsilon_5'$}

\smallskip

In view of $3\kappa_2\le {3\over 7}<{4\over 7}$,
a similar argument proving (11.2) implies that
$$\Upsilon_4'
\ge \{G_4 + O(\varepsilon)\} \, \Pi(x),
\leqno(11.5)$$
where
$$G_4
:= {1\over \kappa_1 e^\gamma}
\int_{\kappa_1}^{\kappa_2} {\d t\over t}
\int_{t}^{\kappa_2} f\bigg({4/7-t-u\over \kappa_1}\bigg) {\d u\over u}.$$

Our assumptions on $\kappa_1$ and $\rho$ imply that
$p_1^2 p_2\le x^{4/7-\varepsilon}$ and $p_1^2\le x^{4/7-\varepsilon}$.
As before we can apply (2.5) of Lemma 2.2 and Lemma 2.7 to get
$$\Upsilon_5'
\ge \{G_5 + O(\varepsilon)\} \, \Pi(x),
\leqno(11.6)$$
where
$$\eqalign{
& G_5
:= {1\over \kappa_1 e^\gamma}
\int_{\kappa_1}^{\kappa_2} {\d t\over t}
\int_{\kappa_2}^{4/7-2\kappa_1-t} f\bigg({4/7-t-u\over 
\kappa_1}\bigg) {\d u\over u}.
\cr}$$

\smallskip

$4^\circ$
{\it Upper bounds of $\Upsilon_i'$ for $i=7, 8, 9, 10, 13, 14$}

\smallskip

We shall apply the technique of [12].
Since $\sigma_1>{2\over 7}$,
the quantity $\Upsilon_7'$ is equal to the number of primes $p\le x$
such that $p+2=p_1p_2p_3$
with $x^{\sigma_1} \le p_1<p_2< (x/p_1)^{1/2}$ and $p_3\ge p_2$.

Introduce the set
$${\cal B} := \{b-2 : b = p_1p_2p_3\le x, \; 
x^{\sigma_1}\le p_1<p_2<p_3\}.$$
Then we have
$$\Upsilon_7' = S({\cal B};{\cal P}(2), x^{1/2}) + O(x^{1/2}).$$
Let $\Delta := 1 + {\cal L}^{-4}$.
We cover the set ${\cal B}$ by cuboids
$${\cal B}(t_1, t_2, t_3)
:= \big\{b-2 : b= p_1p_2p_3\le x, \;
p_i\in [\Delta^{t_i}, \Delta^{t_i+1}) \,\,\, \hbox{for} \,\,\, 1\le 
i\le 3\big\}$$
where $t_i$ are integers satisfying
$x^{\sigma_1}\le \Delta^{t_1}\le \Delta^{t_2}\le \Delta^{t_3}$ and 
$\Delta^{t_1+t_2+t_3+3}\le x$.
In view of $x^{2/7}\le p_2\le x^{(1-\sigma_2)/2}\le x^{2/5}$,
Lemma 2.5 with the choice
$$\alpha_m
= \cases{
1 & if $m=p_1p_3$
\cr\noalign{\smallskip}
0 & otherwise
\cr},
\qquad
\beta_n
= \cases{
1 & if $n=p_2$
\cr\noalign{\smallskip}
0 & otherwise
\cr}$$
implies the inequality
$$\sum_{(q, 2)=1} \lambda_l^+(q)
\bigg(|{\cal B}(t_1, t_2, t_3)_q|-{|{\cal B}(t_1, t_2, t_3)|\over 
\varphi(q)}\bigg)
\ll {x\over (\log x)^{18}},$$
where $\lambda_l^+(q)$ is well factorable of order 1 and of level 
$Q=x^{\theta(t_2)}$
with $\theta(t_2)=(2+t_2)/4$.

Thus we find by (2.4) of Lemma 2.2,
$$S({\cal B}(t_1, t_2, t_3);{\cal P}(2), x^{1/2})
\le {2C \{1+O(\varepsilon)\}\over \theta(t_2) {\cal L}} 
|{\cal B}(t_1, t_2, t_3)|
+ O\bigg({x\over (\log x)^{18}}\bigg).$$
Since the number of cuboids ${\cal B}(t_1, t_2, t_3)$ is 
$O\big((\log x)^{15}\big)$, we have
$$\eqalign{\sum_{(t_1,t_2,t_3)} {|{\cal B}(t_1, t_2, t_3)|\over \theta(t_2)}
& = \mathop{\sum\,\,\,\sum}_{x^{\sigma_1}\le p_1<p_2\le (x/p_1)^{1/2}}
{4x\over p_1p_2\log(x/p_1p_2) (2+\log p_2/\log x)}
+ O\bigg({x\over (\log x)^2}\bigg)
\cr
& = {4x \{1+O(\varepsilon)\}\over \log x}
\mathop{\int\int}_{\sigma_1\le t\le u\le (1-t)/2}
{\d t \d u\over t u (1-t-u)(2+u)}.
\cr}$$
Combining these estimates, we obtain
$$\Upsilon_7'
\le \{G_7 + O(\varepsilon)\} \, \Pi(x),
\leqno(11.7)$$
where
$$G_7 := 8 \mathop{\int\int}_{\sigma_1\le t\le u\le (1-t)/2}
{\d t \d u\over t u (1-t-u)(2+u)}.$$

Analogously we have
$$\Upsilon_8'
\le \{G_8 + O(\varepsilon)\} \, \Pi(x),
\leqno(11.8)$$
where
$$G_8 := 8 \mathop{\int\int}_{\sigma_2\le t\le u\le (1-t)/2}
{\d t \d u\over t u (1-t-u)(2+u)}.$$

For $\Upsilon_9'$, the assumption $3\sigma_1+\kappa_1>1$ allows us to write
$$\Upsilon_9' = S({\cal B}'; {\cal P}(2), x^{1/2}) + O(x^{1/2})$$
with
$${\cal B}' := \{b-2 : b = p_1p_2p_3\le x, \;
x^{\kappa_1}\le p_1<x^{\sigma_1}\le p_2<p_3\}.
\leqno(11.9)$$
We decompose ${\cal B}' = {\cal B}_1'\cup \cdots\cup {\cal B}_6'$,
where ${\cal B}_1'$, $\dots$, ${\cal B}_6'$ are defined as in (11.9)
but we add respectively the extra conditions 
$$\eqalign{
& \hbox{$p_1\le x^{1/10}$ in ${\cal B}_1'$};
\cr\noalign{\smallskip}
& \hbox{$p_1>x^{1/10}$ and $p_1p_2\le x^{1/2}$ in ${\cal B}_2'$};
\cr\noalign{\smallskip}
& \hbox{$p_1>x^{1/10}$ and $p_1^{-2} p_2^8>x^3$ in ${\cal B}_3'$};
\cr\noalign{\smallskip}
& \hbox{$p_2>x^{2/5}$ and $p_1^{-2} p_2^8\le x^3$ in ${\cal B}_4'$};
\cr\noalign{\smallskip}
& \hbox{$x^{(1-\sigma_2)/2}<p_2\le x^{2/5}$ and $p_1p_2>x^{1/2}$ in 
${\cal B}_5'$};
\cr\noalign{\smallskip}
& \hbox{$p_2\le x^{(1-\sigma_2)/2}$ and $p_1p_2>x^{1/2}$ in ${\cal B}_6'$}.
\cr}$$
Again we can use Lemma 2.5
with
$$\eqalign{
& \hbox{$\theta(\log p_1/\log x) = (1+2\log p_1/\log x)/2$ for ${\cal B}_1'$};
\cr\noalign{\smallskip}
& \hbox{$\theta(\log p_1/\log x) = (5-2\log p_1/\log x)/8$ for ${\cal 
B}_2'$ and ${\cal B}_3'$};
\cr\noalign{\smallskip}
& \hbox{$\theta(\log p_2/\log x) = 1-\log p_2/\log x$ for ${\cal B}_4'$};
\cr\noalign{\smallskip}
& \hbox{$\theta(\log p_2/\log x) = (2+\log p_2/\log x)/4$ for ${\cal 
B}_5'$ and ${\cal B}_6'$}.
\cr}$$
Then we have
$$\Upsilon_9'
\le \{G_9 + O(\varepsilon)\} \, \Pi(x),
\leqno(11.10)$$
where
$$\eqalign{G_9
& = 4 \mathop{\int\int}_{\scriptstyle\kappa_1\le t\le \sigma_1\le u\le (1-t)/2
\atop\scriptstyle t\le 1/10}
{\d t\d u\over t u (1-t-u) (1+2t)}
\cr
& \quad
+ 16 \mathop{\int\int}_{\scriptstyle\kappa_1\le t\le \sigma_1\le u\le (1-t)/2
\atop\scriptstyle t\ge 1/10, \, t+u\le 1/2}
{\d t\d u\over t u (1-t-u) (5-2t)}
\cr
& \quad
+ 16 \mathop{\int\int}_{\scriptstyle\kappa_1\le t\le \sigma_1\le u\le (1-t)/2
\atop\scriptstyle t\ge 1/10, \, 8u\ge 2t+3}
{\d t\d u\over t u (1-t-u) (5-2t)}
\cr
& \quad
+ 2 \mathop{\int\int}_{\scriptstyle\kappa_1\le t\le \sigma_1\le u\le (1-t)/2
\atop\scriptstyle u\ge 2/5, \, 8u\le 2t+3}
{\d t\d u\over t u (1-t-u) (1-u)}
\cr
& \quad
+ 8 \mathop{\int\int}_{\scriptstyle\kappa_1\le t\le \sigma_1\le u\le (1-t)/2
\atop\scriptstyle (1-\sigma_1)/2\le u\le 2/5, \, t+u\ge 1/2}
{\d t\d u\over t u (1-t-u) (2+u)}
\cr
& \quad
+ 8 \mathop{\int\int}_{\scriptstyle\kappa_1\le t\le \sigma_1\le u\le (1-t)/2
\atop\scriptstyle u\le (1-\sigma_1)/2, \, t+u\ge 1/2}
{\d t\d u\over t u (1-t-u) (2+u)}.
\cr}$$

Similarly in view of the assumption $2\sigma_1+\sigma_2+\kappa_2>1$, 
we can prove
$$\Upsilon_{10}'
\le \{G_{10} + O(\varepsilon)\} \, \Pi(x),
\leqno(11.11)$$
where
$$\eqalign{G_{10}
& = 16 \mathop{\int\int}_{\scriptstyle\kappa_2\le t\le \sigma_2\le u\le (1-t)/2
\atop\scriptstyle t+u\le 1/2}
{\d t\d u\over t u (1-t-u) (5-2t)}
\cr
& \quad
+ 16 \mathop{\int\int}_{\scriptstyle\kappa_2\le t\le \sigma_2\le u\le (1-t)/2
\atop\scriptstyle 8u\ge 2t+3}
{\d t\d u\over t u (1-t-u) (5-2t)}
\cr
& \quad
+ 2 \mathop{\int\int}_{\scriptstyle\kappa_2\le t\le \sigma_2\le u\le (1-t)/2
\atop\scriptstyle 2/5\le u\le (2\kappa_2+3)/8}
{\d t\d u\over t u (1-t-u) (1-u)}
\cr
& \quad
+ 2 \mathop{\int\int}_{\scriptstyle\kappa_2\le t\le \sigma_2\le u\le (1-t)/2
\atop\scriptstyle u\ge (2\kappa_2+3)/8, \, 8u\le 2t+3}
{\d t\d u\over t u (1-t-u) (1-u)}
\cr
& \quad
+ 8 \mathop{\int\int}_{\scriptstyle\kappa_2\le t\le \sigma_2\le u\le (1-t)/2
\atop\scriptstyle 1/2-\kappa_2\le u\le 2/5, \, t+u\ge 1/2}
{\d t\d u\over t u (1-t-u) (2+u)}
\cr
& \quad
+ 8 \mathop{\int\int}_{\scriptstyle\kappa_2\le t\le \sigma_2\le u\le (1-t)/2
\atop\scriptstyle (1-\sigma_2)/2\le u\le 1/2-\kappa_2, \, t+u\ge 1/2}
{\d t\d u\over t u (1-t-u) (2+u)}
\cr
& \quad
+ 8 \mathop{\int\int}_{\scriptstyle\kappa_2\le t\le \sigma_2\le u\le (1-t)/2
\atop\scriptstyle u\le (1-\sigma_2)/2, \, t+u\ge 1/2}
{\d t\d u\over t u (1-t-u) (2+u)}.
\cr}$$

More easily we can prove that
$$\Upsilon_i
\le \{G_i + O(\varepsilon)\} \, \Theta(N)
\qquad
(i = 13, 14),
\leqno(11.12)$$
where
$$\eqalign{
G_{13}
& := 4
\int_{\kappa_1}^{1/10} {\d t_1\over t_1(1+2t_1)}
\int_{t_1}^{\kappa_2} {\d t_2\over t_2^2}
\int_{t_2}^{\kappa_2} {\d t_3\over t_3}
\int_{t_3}^{\kappa_2} \omega\bigg({1-t_1-t_2-t_3-t_4\over t_2}\bigg) 
{\d t_4\over t_4}
\cr
& \quad
+ 16
\int_{1/10}^{\kappa_2} {\d t_1\over t_1(5-2t_1)}
\int_{t_1}^{\kappa_2} {\d t_2\over t_2^2}
\int_{t_2}^{\kappa_2} {\d t_3\over t_3}
\int_{t_3}^{\kappa_2} \omega\bigg({1-t_1-t_2-t_3-t_4\over t_2}\bigg) 
{\d t_4\over t_4},
\cr\noalign{\medskip}
G_{14}
& := 4
\int_{\kappa_1}^{1/10} {\d t_1\over t_1(1+2t_1)}
\int_{t_1}^{\kappa_2} {\d t_2\over t_2^2}
\int_{t_2}^{\kappa_2} {\d t_3\over t_3}
\int_{\kappa_2}^{4/7-2\kappa_1-t_3} 
\omega\bigg({1-t_1-t_2-t_3-t_4\over t_2}\bigg) {\d t_4\over t_4}
\cr
& \quad
+ 16
\int_{1/10}^{\kappa_2} {\d t_1\over t_1(5-2t_1)}
\int_{t_1}^{\kappa_2} {\d t_2\over t_2^2}
\int_{t_2}^{\kappa_2} {\d t_3\over t_3}
\int_{\kappa_2}^{4/7-2\kappa_1-t_3} 
\omega\bigg({1-t_1-t_2-t_3-t_4\over t_2}\bigg) {\d t_4\over t_4}.
\cr}$$

Inserting these estimations and the trivial lower bounds 
$\Upsilon_i'\ge 0$ $(i = 6, 11, 12, 15)$ into (9.10),
we obtain
$$\pi_{1,2}(x)
\ge \{G(\kappa_1, \kappa_2, \rho, \sigma_2, \sigma_1) + 
O(\varepsilon)\} \, \Pi(x),$$
where
$$\eqalign{G(\kappa_1, \kappa_2, \rho, \sigma_2, \sigma_1)
& := \textstyle {1\over 4} (4G_0
- G_1
- G_2
- G_3
+ G_4
+ G_5
\cr
& \quad
- 2G_7
- 2G_8
- G_9
- G_{10}
- G_{13}
- G_{14}).
\cr}$$
Taking
$\kappa_1 = {2\over 21}$,
$\kappa_2 = {13\over 100}$,
$\rho = {2\over 7}$,
$\sigma_2 = {36\over 125}$
and
$\sigma_1 = {332\over 1000}$,
a numerical computation gives us
$$\eqalign{\textstyle G(\kappa_1, \kappa_2, \rho, \sigma_2, \sigma_1)
& \ge \textstyle {1\over 4}(
   4\times 5.894705
- 1.611441
- 7.921437
- 6.736885
\cr
& \quad
+ 0.270916
+ 0.913995
- 2\times 0.000124
- 2\times 0.145114
\cr
& \quad
- 1.790090
- 1.930545
- 0.006814
- 0.059690)
\cr
& > 1.10409.
\cr}$$
[For the integrals $F_{13}$ and $F_{14}$, we make use of 
$\omega(u)\le 0.561522$ for $u\ge 3.5$
and $\omega(u)\le 0.567144$ for $u\ge 2$.]
This completes the proof of Theorem 4.
\hfill
$\square$

\vskip 10mm

\centerline{\bf References}

\bigskip

\item{[1]}
{\author E. Bombieri \& H. Davenport},
Small differences between prime numbers,
{\it Proc. Roy. Soc}. Ser. A {\bf 239} (1966), 1--18.

\item{[2]}
{\author E. Bombieri, J.B. Friedlander \& H. Iwaniec},
Primes in arithmetic progressions to large moduli,
{\it Acta Math}. {\bf 156} (1986), 203--251.

\item{[3]}
{\author Y.C. Cai},
A remark on Chen's theorem,
{\it Acta Arith.} {\bf 102} (2002), 339--352.

\item{[4]}
{\author Y.C. Cai \& M.G. Lu},
Chen's theorem in short intervals,
{\it Acta Arith.} {\bf 91} (1999), 311--323.

\item{[5]}
{\author Y.C. Cai \& M.G. Lu},
On Chen's theorem,
in: {\it Analytic number theory} (Beijing/Kyoto, 1999), 99--119, 
Dev. Math. {\bf 6}, Kluwer Acad. Publ., Dordrecht, 2002.
	
\item{[6]}
{\author Y.C. Cai \& M.G. Lu},
On the upper bound for $\pi_2(x)$,
{\it Acta Arith.}, to appear.

\item{[7]}
{\author J.R. Chen},
On the representation of a large even integer
as the sum of a prime and the product of at most two primes,
{\it Sci. Sinica} {\bf 16} (1973), 157--176.

\item{[8]}
{\author J.R. Chen},
On the representation of a large even integer
as the sum of a prime and the product of at most two primes (II),
{\it Sci. Sinica} {\bf 21} (1978), 421--430.

\item{[9]}
{\author J.R. Chen},
Further improvement on the constant in the proposition `1+2':
On the representation of a large even integer
as the sum of a prime and the product of at most two primes (II) (in Chinese),
{\it Sci. Sinica} {\bf 21} (1978), 477--494.

\item{[10]}
{\author J.R. Chen},
On the Goldbach's problem and the sieve methods,
{\it Sci. Sinica} {\bf 21} (1978), 701--739.

\item{[11]}
{\author E. Fouvry},
Autour du th\'eor\`eme de Bombieri--Vinogradov,
{\it Ann. scient. Ec. Norm. Sup.} {\bf 20} (1987), 617--640.

\item{[12]}
{\author E. Fouvry \& F. Grupp},
On the switching principle in sieve theory,
{\it J. reine angew. Math.} {\bf 370} (1986), 101--125.

\item{[13]}
{\author E. Fouvry \& F. Grupp},
Weighted sieves and twin prime type equations,
{\it Duke J. Math.} {\bf 58} (1989), 731--748.

\item{[14]}
{\author H. Halberstam \& H.-E. Richert},
{\it Sieve Methods},
Academic Press, London, 1974.

\item{[15]}
{\author G.H. Hardy \& J.E. Littlewood},
Some problems of `partitio numerorum' III :
On the expression of a number as a sum of primes,
{\it Acta Math.} {\bf 44} (1923), 1--70.

\item{[16]}
{\author H. Iwaniec},
Primes of the type $\phi(x,\,y)+A$ where $\phi$ is a quadratic form,
{\it Acta Arith.} {\bf 21} (1972), 203--234.

\item{[17]}
{\author H. Iwaniec},
A new form of the error term in the linear sieve,
{\it Acta Arith.} {\bf 37} (1980), 307--320.

\item{[18]}
{\author H. Iwaniec \& J. Pomykala},
Sums and difference of quartic norms,
{\it Mathematik} {\bf 40} (1993), 233--245.

\item{[19]}
{\author C.B. Pan},
On the upper bound of the number of ways to represent an even 
integer as a sum of two primes,
{\it Sci. Sinica} {\bf 23} (1980), 1367--1377.

\item{[20]}
{\author C.D. Pan},
A new application of the Yu.V. Linnik large sieve method,
{\it Chinese Math. Acta} {\bf 5} (1964), 642--652.

\item{[21]}
{\author C.G. Pan \& X.X. Ding},
A new mean value theorem,
{\it Sci. Sinica} (1979), Special Issue II on Math., 149--161.

\item{[22]}
{\author C.D. Pan \& C.B. Pan},
{\it Goldbach Conjecture},
Science Press, Beijing, China, 1992.

\item{[23]}
{\author P.M. Ross},
A short intervals result in additive prime number theory,
{\it J. London Math. Soc}. (2), {\bf 17} (1978), 219--227.

\item{[24]}
{\author S. Salerno \& A. Vitolo},
$p+2=P_2$ in short intervals,
{\it Note Mat.} {\bf 13} (1993), 309--328.

\item{[25]}
{\author A. Selberg},
On elementary methods in prime number theory and their limitations,
{\it Den 11te Skandinaviske Matemtikerkongnes, Trondheim}, 1949,
Johan Grundt Tanums Forlag, Oslo, 1952, 13--22.

\item{[26]}
{\author J. Wu},
Sur la suite des nombres premiers jumeaux,
{\it Acta Arith.} {\bf 55} (1990), 365--394.

\item{[27]}
{\author J. Wu},
Th\'eor\`emes g\'en\'eralis\'es de Bombieri--Vinogradov dans les 
petits intervalles,
{\it Quart. J. Math. Oxford} (2), {\bf 44} (1993), 109--128.

\item{[28]}
{\author J. Wu},
Sur l'\'equation $p + 2 = P_2$ dans les petits intervalles,
{\it J. London Math. Soc}. (2), {\bf 49} (1994), 61--72.

\bigskip

Institut Elie Cartan

UMR 7502 UHP-CNRS-INRIA

Universit\'e Henri Poincar\'e (Nancy 1)

54506 Vand\oe uvre--l\`es--Nancy

FRANCE

e--mail: wujie@iecn.u-nancy.fr

\end
\bye